\documentclass[a4paper, 12pt]{article}

\usepackage{amsmath} 
\usepackage{theorem}
\usepackage{amssymb}
\usepackage{amsfonts}
\usepackage{latexsym}
\usepackage{amscd}
\usepackage{diagramb}
\usepackage{yfonts}
\usepackage{hyperref}
\usepackage{graphicx}
\usepackage{enumerate}
\usepackage{cancel}
\usepackage{extarrows}
\usepackage{upgreek}

\usepackage{pxfonts} 

\usepackage{stmaryrd}

\input xy
\xyoption{all}

\DeclareMathAlphabet{\mathpzc}{OT1}{pzc}{m}{it}

\usepackage{xspace,colortbl}
\usepackage{color}

\definecolor{verde}{rgb}{0.,0.7,0.}
\definecolor{indigo}{rgb}{.18, .34, .78}
\definecolor{indigo1}{rgb}{.18, .24, .78}
\definecolor{indigo2}{rgb}{.18, .14, .78}
\definecolor{indigo3}{rgb}{.18, 0., .78}
\definecolor{rojo}{rgb}{1,0,0}
\definecolor{negro}{rgb}{0,0,0}
\definecolor{lila}{rgb}{.46, .16, .78}
\definecolor{lila1}{rgb}{.46, .16, .86}
\definecolor{lila2}{rgb}{.56, .16, .86}
	\definecolor{lila3}{rgb}{.63, .16, .78}
\definecolor{lila4}{rgb}{.7, .16, .78}
\definecolor{lila5}{rgb}{.78, .26, .78}
\definecolor{lila6}{rgb}{.6, 0., .78}

\theoremstyle{plain}
\theoremheaderfont{\normalfont\bfseries}

\newtheorem{thm}{Theorem}[section]

\newtheorem{lma}[thm]{Lemma}
\newtheorem{cor}[thm]{Corollary}
\newtheorem{defn}[thm]{Definition}

\theorembodyfont{\rmfamily}

\newtheorem{rem}[thm]{Remark}
\newtheorem{prop}[thm]{Proposition}

\newcommand{\qed}{\hfill\quad\fbox{\rule[0mm]{0,0cm}{0,0mm}}  \par\bigskip}

\newcommand{\x}{\mbox{-}}

\newcommand{\s}{\hspace{0,06cm}}
\newcommand{\Cat}{\operatorname {Cat}}
\newcommand{\Mat}{\operatorname {Mat}}

\newcommand{\Mnd}{{\rm Mnd}}

\newcommand{\Dbl}{{\rm Dbl}}
\newcommand{\Ps}{{\rm Ps}}

\newcommand{\Fun}{{\rm Fun}}

\newcommand{\Comp}{{\rm Comp}}
\newcommand{\Span}{\operatorname{Span}}

\newcommand{\Aa}{{\mathbb A}}
\newcommand{\Bb}{{\mathbb B}}
\newcommand{\Cc}{{\mathbb C}}
\newcommand{\Dd}{{\mathbb D}}
\newcommand{\I}{{\mathbb I}}
\newcommand{\aaa}{{\textfrak{a}}}

\newcommand{\Del}{\boxtimes}
\newcommand{\comp}{\circ}
\newcommand{\iso}{\cong}
\newcommand{\ot}{\otimes}

\newcommand{\C}{{\mathcal C}}

\newcommand{\D}{{\mathcal D}}
\newcommand{\F}{{\mathcal F}}
\newcommand{\G}{{\mathcal G}}

\newcommand{\HH}{{\mathcal H}}

\newcommand{\A}{{\mathcal A}}
\newcommand{\B}{{\mathcal B}}

\def\ul{\underline}

\newcommand{\crta}{\overline}

\newcommand{\Id}{\operatorname {Id}}

\def\K{{\mathcal K}}  

\def\Dd{{\mathbb D}}

\newcommand{\Mod}{\operatorname{Mod}}

\newcommand{\Lax}{\operatorname{Lax}}
\newcommand{\Dist}{\operatorname{Dist}}

\newcommand{\cref}[1]{C.~\ref{c:#1}}

\newcommand{\lelabel}[1]{\label{le:#1}}
\newcommand{\leref}[1]{Lemma~\ref{le:#1}}
\newcommand{\eqlabel}[1]{\label{eq:#1}}
\newcommand{\equref}[1]{(\ref{eq:#1})}
\newcommand{\thlabel}[1]{\label{th:#1}}
\newcommand{\thref}[1]{Theorem~\ref{th:#1}}
\newcommand{\delabel}[1]{\label{de:#1}}
\newcommand{\deref}[1]{Definition~\ref{de:#1}}
\newcommand{\prlabel}[1]{\label{pr:#1}}
\newcommand{\prref}[1]{Proposition~\ref{pr:#1}}
\newcommand{\colabel}[1]{\label{co:#1}}
\newcommand{\coref}[1]{Corollary~\ref{co:#1}}

\newcommand{\selabel}[1]{\label{se:#1}}
\newcommand{\seref}[1]{Section~\ref{se:#1}}
\newcommand{\sslabel}[1]{\label{ss:#1}}
\newcommand{\ssref}[1]{Subsection~\ref{ss:#1}}

\newdir{:=}{{}}
\newcommand{\fit}[3]{\ar@{:=}@/{#3}/[#1] |{\Downarrow #2} }

\begin{document}

\title{Bifunctor Theorem and \\ strictification tensor product \\ for double categories with lax double functors}

\author{Bojana Femi\'c \vspace{6pt} \\
{\small Mathematical Institite of  \vspace{-2pt}}\\
{\small Serbian Academy of Sciences and Arts } \vspace{-2pt}\\
{\small Kneza Mihaila 36,} \vspace{-2pt}\\
{\small 11 000 Belgrade, Serbia} \vspace{-2pt}\\
{\small femicenelsur@gmail.com}}

\date{}
\maketitle

\begin{abstract}
We introduce a candidate for the inner hom for $Dbl_{lx}$, the category of double categories and lax double functors, and 
characterize a lax double functor into it obtaining a lax double quasi-functor. The latter consists of a pair of lax double functors 
with four 2-cells resembling distributive laws. We extend this characterization to a double category isomorphism  
$q\x\Lax_{hop}(\Aa\times\Bb,\Cc)  \iso \Lax_{hop}(\Aa, \llbracket\Bb,\Cc\rrbracket)$. We show that instead of a Gray monoidal product 
in $Dbl_{lx}$ we obtain a product that in a sense strictifies lax double quasi-functors. 
We explain why laxity of double functors hinders $\llbracket-,-\rrbracket$ from making $Dbl_{lx}$ a closed and enriched category 
over 2-categories (or double categories). We prove a bifunctor theorem by which 
certain type of lax double quasi-functors give rise to lax double functors on the Cartesian product, extend it to a 
double functor $q\x\Lax_{hop}^{st}(\Aa\times\Bb,\Cc)\to\Lax_{hop}(\Aa\times\Bb,\Cc)$ and show how it restricts to a double equivalence. The 
(un)currying double functors are studied. We prove that a lax double functor 
from the trivial double category is a monad in the codomain double category, and show that our above double functor in the form 
$q\x\Lax_{hop}(*\times *,\Dd)\to\Lax_{hop}(*,\Dd)$ recovers the specification $\Comp(\Dd):\Mnd\Mnd(\Dd)\to\Mnd(\Dd)$ 
of the natural transformation $\Comp$. 
\end{abstract}


\section{Introduction}

In recent years the importance of double categories, and more generally of internal categories, has been increasingly recognized in the literature. 
It was observed by various authors ({\em e.g.} \cite{Shul, Shul:Fr, Doug, FGK}) that it is often more convenient to work in the internal {\em i.e.} double 
categorical setting, than in the bicategorical one. In $\Mod$-type bicategories the 1-cells are (also) ``objects'' but of different nature 
than the 0-cells, and they do not present real maps between 0-cells, so the latter are missing in the picture. 
This also happens in the 2-category $\Mnd(\K)$ of monads in a 2-category $\K$, 
a fact which gave rise to the introduction of the (pseudo) double category of monads in a (pseudo) double category in \cite{FGK}. 
Namely, it is known that various algebraic structures can be expressed as monads in suitable bicategories, but that 
the corresponding morphisms are not morphisms of monads, considered as 1-cells in $\Mnd(\K)$ (up to weakening of 2-categories). Rather than being 
real morphisms of monads, 1-cells in $\Mnd(\K)$ are distributive laws between them. In 
order to include morphisms of monads, vertical 1-cells among monads are introduced in \cite{FGK} as well as the corresponding (pseudo) double category.

\smallskip

For pseudo double categories the Strictification Theorem is proved in \cite[Section 7.5]{GP:L}. 
One has that the category of pseudo double categories and pseudo double functors is equivalent 
to the category $Dbl_{st}$ of double categories and strict double functors. However, 
the lax double functors can not be ``strictified'', so the category of double categories and lax double functors $Dbl_{lx}$ is properly more general 
than $Dbl_{st}$. Apart from the Cartesian monoidal product known in the literature for both categories, a Gray type monoidal product $\ot$ was introduced 
in \cite{Gabi} for $Dbl_{st}$. In \cite{Fem} we have described the monoidal category structure of $\Aa\ot\Bb$ for double categories $\Aa$ and $\Bb$. 
In the present paper we show that although one can construct natural candidates $\Aa\ot\Bb$ and $\llbracket\Aa,\Bb\rrbracket$ for the tensor product, 
respectively inner hom, for the category $Dbl_{lx}$, it turns out that $\Aa\ot\Bb$ does not satisfy the expected universal property, and that 
$\llbracket\Aa,\Bb\rrbracket$ is not a bifunctor. Instead of a Gray type monoidal product for $Dbl_{lx}$, we prove that $\Aa\ot\Bb$ 
satisfies a universal property by which lax double functors $\Aa\to\llbracket\Bb,\Cc\rrbracket$ bijectively correspond to strict double functors $\Aa\ot\Bb\to\Cc$. 
\smallskip

Recent results on 2-categories of \cite{FMS} naturally inspired us to study the analogous properties in double categories. 
Namely, in {\em loc. cit.} conditions were studied for two families of lax functors with a common codomain 2-category to collate into a lax bifunctor, 
{\em i.e.} a lax functor on the Cartesian product 2-category. (This question corresponds to a 2-category analogue of the first Proposition in \cite{McL}, 
page 37.) 
The authors proved a version of a bifunctor theorem for lax functors, which even extends to a 2-functor $K:\Dist(\Aa,\Bb, \Cc)
\to\Lax_{op}(\Aa\times\Bb,\Cc)$ into the corresponding 2-category of bifunctors. The 2-functor $K$ is proved to restrict to a 2-equivalence on certain 
sub-2-categories. 

We noticed that the conditions 
found by the authors to fulfill the above-mentioned goal are the weak (lax) version of the 2-categorical part of the data of a cubical double functor, 
that we introduced in \cite[Definition 2.2]{Fem}. 
Namely, starting from the Gray type closed monoidal structure on the category $(Dbl_{st}, \ot)$ 
constructed in \cite{Gabi}, we characterized in \cite[Proposition 2.1]{Fem} a strict double functor $F:\Aa\to \llbracket\Bb,\Cc\rrbracket$ 
with the codomain the inner hom object. We obtained that $F$ corresponds to two families of double functors with codomain $\Cc$, satisfying a 
longer list of conditions. The latter pair of families we called a cubical double functor, in analogy to \cite[Section 4.2]{GPS}. 

Our above-mentioned observation led us to conjecture, and to prove it became one of our main goals in this paper, 
that by weakening our characterization in \cite[Proposition 2.1]{Fem} to double categories and 
{\em lax} double functors, 
would lead to a double functor into the corresponding double category of lax double bifunctors, 
generalizing the above 2-functor $K$ to a double categorical setting. 
We attained this goal in Section 4 (concretely, we constructed the double functor $\F$ in \prref{F}), where we also identify a double equivalence functor 
which is a restriction of $\F$ (\thref{equiv-fun}). 
We present this and the rest of our results in more details in the continuation.

\medskip

We start by introducing the double category $\Lax_{hop}(\Aa, \Bb)$ of lax double functors of double categories $\Aa\to\Bb$, horizontal oplax 
transformations as 1h-cells (horizontal 1-cells), vertical lax transformations as 1v-cells (vertical 1-cells), and modifications. 
It is a generalization of 
inner-homs in $Dbl_{st}$ from \cite[Section 2.2]{Gabi} to lax double functors. 
In \seref{mon str} we explore $\llbracket\Aa,\Bb\rrbracket:=\Lax_{hop}(\Aa, \Bb)$ as a candidate for inner-hom in $Dbl_{lx}$ and we show why 
it fails to be one. 
We then characterize a lax double functor $\F:\Aa\to\llbracket\Bb,\Cc\rrbracket$ as a pair of two families of lax double functors into $\Cc$, 
satisfying a list of properties, that we call a {\em lax double quasi-functor}. We use this characterization to describe a double 
category $\Aa\ot\Bb$ in \ssref{Gray prod}. We finish Section 3 by showing why laxity of double functors prevents $-\ot-$ from being a monoidal product 
and also $Dbl_{lx}$ from being a category enriched over 2-categories (and then also over double categories). 

In Section 4 we introduce the double category $q\x\Lax_{hop}(\Aa\times\Bb,\Cc)$ of lax double quasi-functors, 
horizontal oplax transformations as 1h-cells, vertical lax transformations as 1v-cells, and modifications. 
In this double categorical context, the 1-cells in $q\x\Lax_{hop}(\Aa\times\Bb,\Cc)$ 
have four defining axioms $HOT^q_1-HOT^q_4$, whereas in the analogous 2-categorical situation, the 1-cells of the 2-category $\Dist(\A,\B, \C)$ for 
2-categories $\A,\B,\C$ from \cite{FMS} have a single axiom, corresponding to our $HOT^q_1$, called a Yang-Baxter equation therein. 
(The 2-category $\Dist(\A,\B, \C)$, in turn, is a lax version on 0-cells and an oplax version on 1-cells of the 2-category 
$q\x\Fun(\A\times\B, \C)$ from \cite[Section I.4]{Gray}. Namely, the 0-cells in $\Dist(\A,\B, \C)$ are pairs of families of {\em lax} functors of 
2-categories together with their distributive law, whereas the 0-cells in $q\x\Fun(\A\times\B, \C)$ are ``quasi-functors of two variables'' defined 
in terms of pairs of families of {\em strict 2-functors}. Morphisms of distributive laws of lax functors, {\em i.e.} 1-cells in $\Dist(\A,\B, \C)$ 
are oplax natural transformations, while the 1-cells in $q\x\Fun(\A\times\B, \C)$ are quasi-natural transformations, which are lax 
(see I.4.1 and I.3.3 of \cite{Gray}).) In \ssref{2-cells} we prove that the double categories $q\x\Lax_{hop}(\Aa\times\Bb,\Cc)$ and 
$\Lax_{hop}(\Aa, \llbracket\Bb,\Cc\rrbracket)$ are isomorphic. 


The objective in Section 5 is to find a description of lax double quasi-functors in terms of ordinary lax double functors on the Cartesian product $\Aa\times\Bb$. 
In order to obtain this description we find it necessary to require that the vertical lax transformations in $\llbracket\Bb,\Cc\rrbracket$ be strict, 
in which case we obtain a double category $\llbracket\Bb,\Cc\rrbracket^{st}$. Concretely, it is on the isomorphic counterpart 
$q\x\Lax_{hop}^{st}(\Aa\times\Bb,\Cc)$ of $\Lax_{hop}(\Aa, \llbracket\Bb,\Cc\rrbracket^{st})$ (in the isomorphism from the previous paragraph) that we 
managed to construct the double functor $\F: q\x\Lax_{hop}^{st}(\Aa\times\Bb,\Cc) \to \Lax_{hop}(\Aa\times\Bb,\Cc)$ in \prref{F}. Restricting to certain double 
subcategories we obtain a double equivalence $\F': q\x\Lax_{hop}^{st-u}(\Aa\times\Bb,\Cc) \to \Lax_{hop}^{ud}(\Aa\times\Bb,\Cc)$ in \thref{equiv-fun}. In terms 
of pseudo double functors it comes down to a double equivalence $\F'': q\x\Ps_{hop}^{st}(\Aa\times\Bb,\Cc) \to \Ps_{hop}(\Aa\times\Bb,\Cc).$

The double category isomorphism from Section 4, the double functor $\F$ above and the double equivalences $\F'$ and $\F''$ are generalizations to double categories 
of the corresponding results in \cite{FMS}. In Section 6 we show applications of these results of ours in three different contexts. 
In \equref{uncurry} we obtain a double categorical version $\Lax_{hop}(\Aa, \llbracket\Bb,\Cc\rrbracket^{st})\to \Lax_{hop}(\Aa\times\Bb,\Cc)$ of the 
``uncurrying'' 2-functor $J$ from \cite[Section 4]{FMS} and establish a ``currying'' functor, {\em i.e.} a 2-equivalence 
$\Lax_{hop}^{ud}(\Aa\times\Bb,\Cc)\simeq\Lax_{hop}^u(\Aa, \llbracket\Bb,\Cc\rrbracket^{st-u})$ in a double categorical setting.

In \ssref{Gray} we show a universal property of $\Aa\ot\Bb$ which extends to a double category isomorphism. The universal property that $\Aa\ot\Bb$ satisfies 
is that for every lax double quasi-functor $H:\Aa\times\Bb\to\Cc$ there is a unique strict double functor $\crta H: \Aa\ot\Bb\to\Cc$ such that $H=\crta{H}J$, 
where $J:\Aa\times\Bb\to\Aa\ot\Bb$ is a naturally obtained lax double quasi-functor.

The final Subsection is devoted to applications to monads in double categories. 
In it we show that a monad in a double category $\Dd$, as defined in \cite[Definition 2.4]{FGK}, is a lax double functor $*\to\Dd$ from the trivial double 
category. Moreover, we obtain isomorphisms of double categories $\Lax_{hop}(*,\Dd)\iso\Mnd(\Dd)$ and $q\x\Lax_{hop}(*\times *,\Dd)\iso\Mnd(\Mnd(\Dd))$. 
We argue that a version of our double functor $\F$ from above, $q\x\Lax_{hop}(*\times *,\Dd)\to\Lax_{hop}(*,\Dd)$, corresponds via the above isomorphisms to the natural transformation $\Comp:\Mnd\Mnd\to\Mnd$ evaluated at the double category $\Dd$. Some prospects of further research are indicated. 

\medskip

The reader is assumed to be familiar with the notion of double categories, 
for the reference we recommend \cite{GP:L, MG}. All double categories in this paper will be strict.

\section{The double category $\Lax_{hop}(\Aa, \Bb)$ and more} \selabel{definitions}

The double category $\Lax_{hop}(\Aa, \Bb)$ consists of lax double functors $\Aa\to\Bb$, horizontal oplax transformations of the 
latter lax double functors as 1h-cells, vertical lax transformations as 1v-cells, 
and modifications with respect to the two types of transformations as 2-cells. For reader's convenience we give the explicit definitions of all 
these notions in this Section. Moreover, we will give the definitions of the notions that we will be using in \seref{mon str}.

Let us first fix the notation in a double category $\Dd$. Objects we denote by $A,B, \dots$, horizontal 1-cells we will call 
for brevity 1h-cells and denote them by $f, f', g, \dots$ (and by $(K,k), (K',k'), (L,l), \dots$ in the Cartesian product of double categories 
$\Aa\times\Bb$), vertical 1-cells we will call 1v-cells and denote by $u,v, U, \dots$, 
and squares we will call just 2-cells and denote them by $\omega, \zeta, \dots$. We denote 
the horizontal identity 1-cell by $1_A$, vertical identity 1-cell by $1^A$ for an object $A\in\Dd$, 
horizontal identity 2-cell $\I_A$ 
on the 1h-cell $1_A$, 
horizontal identity 2-cell on a 1v-cell $u$ by $Id^u$, and vertical identity 2-cell on a 1h-cell $f$ by $Id_f$ (with subindexes we denote 
those identity 1- and 2-cells which come from the horizontal 2-category lying in $\Dd$). 
For a (vertically) globular 2-cell $\alpha$, that is, one whose 1v-cells are identities, we will write $\alpha: f\Rightarrow g$ 
as in bicategories. 
The composition of 1h-cells as well as the horizontal composition of 2-cells we will denote by juxtaposition, while the composition 
of 1v-cells as well the vertical composition of 2-cells we will denote by fractions $\frac{\bullet}{\bullet}$. When combining 
horizontal and vertical composition of 2-cells we will also use the notation: $[\alpha\vert\beta]:=\beta\alpha$ for the horizontal composition.

When dealing with pseudodouble categories, we use the convention that the horizontal direction is weak and the vertical one strict. 
For this reason our lax double functors will be lax in the horizontal direction. We stress this fact only in the definition that follows, 
and will not repeat it afterwards.

\begin{defn}
A {\em (horizontally) lax double functor} $F:\Cc\to\Dd$ between double categories is given by: 1) the data: images on objects, 1h-, 1v- and 2-cells of $\Cc$, 
globular 2-cells: 
{\em compositor} $F_{gf}: F(g)F(f)\Rightarrow F(gf)$ and {\em unitor} $F_A: 1_{F(A)}\Rightarrow F(1_A)$ in $\Dd$, 
and 2) rules (in $\Dd$): 
\begin{itemize}
\item ($\text{vertical functor}$)
$$\text{{\em (lx.f.v1)}} \label{lx.f.v1} \quad\qquad  \frac{F(u)}{F(u')}=F(\frac{u}{u'}), \qquad\qquad 
\text{{\em (lx.f.v2)}} \label{lx.f.v2}\quad\qquad F(1^A)=1^{F(A)};$$
\item $\text{(horizontal functor)}$
$$\text{{\em (lx.f.h1)}} \label{lx.f.h1} \quad\qquad  F(\frac{\omega}{\zeta})=\frac{F(\omega)}{F(\zeta)}, \qquad\qquad 
\text{{\em (lx.f.h2)}} \label{lx.f.h2} \quad\qquad F(Id_f)=Id_{F(f)};$$
\item $\text{(coherence with compositors)}$
$$\text{{\em (lx.f.hex)}} \label{lx.f.hex} \qquad \frac{[F_{gf}\vert Id_{F(h)}]}{F_{h,gf}}=\frac{[Id_{F(f)}\vert F_{hg}]}{F_{hg,f}} \quad\qquad 
\text{{\em (lx.f.u)}} \label{lx.f.u} \qquad \frac{[F_A\vert Id_{F(f)}]}{F_{f1_A}}=\Id_{F(f)}=\frac{[Id_{F(f)}\vert F_B]}{F_{1_Bf}};$$
\item $\text{(naturality of the compositor)}$

$$
\text{{\em (lx.f.c-nat)}} \label{lx.f.c-nat} \qquad
\scalebox{0.86}{  
\bfig
\putmorphism(-150,500)(1,0)[F(A)`F(B)`F(f)]{600}1a
 \putmorphism(450,500)(1,0)[\phantom{F(A)}`F(C) `F(g)]{620}1a

 \putmorphism(-150,50)(1,0)[F(A)`F(B)`F(f')]{600}1a
 \putmorphism(450,50)(1,0)[\phantom{F(A)}`F(C) `F(g')]{620}1a

\putmorphism(-180,500)(0,-1)[\phantom{Y_2}``u]{450}1l
\putmorphism(450,500)(0,-1)[\phantom{Y_2}``]{450}1r
\putmorphism(300,500)(0,-1)[\phantom{Y_2}``v]{450}0r
\putmorphism(1080,500)(0,-1)[\phantom{Y_2}``w]{450}1r
\put(30,280){\fbox{$F(\alpha)$}}
\put(660,280){\fbox{$F(\beta)$}}

\putmorphism(-150,-400)(1,0)[F(A)`F(C) `F(g'f')]{1200}1a

\putmorphism(-180,50)(0,-1)[\phantom{Y_2}``=]{450}1l
\putmorphism(1080,50)(0,-1)[\phantom{Y_3}``=]{450}1r
\put(320,-180){\fbox{$F_{g'f'}$}}

\efig}
\quad=\quad
\scalebox{0.86}{
\bfig
\putmorphism(-150,500)(1,0)[F(A)`F(B)`F(f)]{600}1a
 \putmorphism(450,500)(1,0)[\phantom{F(A)}`F(C) `F(g)]{620}1a
 \putmorphism(-150,50)(1,0)[F(A)`F(C)`F(gf)]{1220}1a

\putmorphism(-180,500)(0,-1)[\phantom{Y_2}``=]{450}1r
\putmorphism(1080,500)(0,-1)[\phantom{Y_2}``=]{450}1r
\put(350,260){\fbox{$F_{gf}$}}

\putmorphism(-150,-400)(1,0)[F(A)`F(C) `F(g'f')]{1200}1a

\putmorphism(-180,50)(0,-1)[\phantom{Y_2}``u]{450}1l
\putmorphism(1120,50)(0,-1)[\phantom{Y_3}``w]{450}1r
\put(320,-180){\fbox{$F(\beta\alpha)$}} 
\efig};$$
\item $\text{(naturality of the unitor)}$

$$\text{{\em (lx.f.u-nat)}} 
\label{lx.f.u-nat} \qquad\quad 
\scalebox{0.86}{
\bfig
\putmorphism(-250,500)(1,0)[F(A)`F(A)` =]{550}1a
 \putmorphism(-250,50)(1,0)[F(A')`F(A')` =]{550}1a
 \putmorphism(-250,-400)(1,0)[F(A')`F(A')` F(1_{A'})]{550}1a

\putmorphism(-280,500)(0,-1)[\phantom{Y_2}``F(u)]{450}1l
 \putmorphism(-280,70)(0,-1)[\phantom{F(A)}` `=]{450}1l

\putmorphism(300,500)(0,-1)[\phantom{Y_2}``F(u)]{450}1r
\putmorphism(300,70)(0,-1)[\phantom{Y_2}``=]{450}1r
\put(-120,290){\fbox{$Id^{F(u)}$}}
\put(-100,-150){\fbox{$F_{A'}$}}
\efig}
=
\scalebox{0.86}{
\bfig
\putmorphism(-250,500)(1,0)[F(A)`F(A)` =]{550}1a
 \putmorphism(-250,50)(1,0)[F(A)`F(A)` F(1_A)]{550}1a
 \putmorphism(-250,-400)(1,0)[F(A')`F(A')` F(1_{A'})]{550}1a

\putmorphism(-280,500)(0,-1)[\phantom{Y_2}``= ]{450}1l
 \putmorphism(-280,70)(0,-1)[\phantom{F(A)}` `F(u)]{450}1l

\putmorphism(300,500)(0,-1)[\phantom{Y_2}``=]{450}1r
\putmorphism(300,70)(0,-1)[\phantom{Y_2}``F(u)]{450}1r
\put(-120,290){\fbox{$F_A$}}
\put(-120,-150){\fbox{$F(Id^u)$}}
\efig},
\qquad
$$
\end{itemize}
where $u,u'$ are composable 1v-cells, $\omega, \zeta$ vertically composable 2-cells, and $f,g,h$ composable 1h-cells. 

A {\em pseudodouble functor} is a lax double functor whose compositor and unitor 2-cells are invertible. 
\end{defn}


We now define horizontal oplax and vertical lax transformations between lax double functors and their modifications, and their respective compositions.  

\begin{defn} \delabel{hor nat tr}
A {\em horizontal oplax transformation} $\alpha$ between lax double functors $F,G: \Aa\to\Bb$ consists of the following:
\begin{enumerate}
\item for every 0-cell $A$ in $\Aa$ a 1h-cell $\alpha(A):F(A)\to G(A)$ in $\Bb$,
\item for every 1v-cell $u:A\to A'$ in $\Aa$ a 2-cell in $\Bb$:
$$
\scalebox{0.86}{
\bfig
\putmorphism(-150,50)(1,0)[F(A)`G(A)`\alpha(A)]{560}1a
\putmorphism(-150,-320)(1,0)[F(A')`G(A')`\alpha(A')]{600}1a
\putmorphism(-180,50)(0,-1)[\phantom{Y_2}``F(u)]{370}1l
\putmorphism(410,50)(0,-1)[\phantom{Y_2}``G(u)]{370}1r
\put(30,-110){\fbox{$\alpha^u$}}
\efig}
$$
\item 
for every 1h-cell $f:A\to B$  in $\Aa$ there is a 2-cell in $\Bb$:
$$
\scalebox{0.86}{
\bfig
 \putmorphism(-170,500)(1,0)[F(A)`F(B)`F(f)]{540}1a
 \putmorphism(360,500)(1,0)[\phantom{F(f)}`G(B) `\alpha(B)]{560}1a
 \putmorphism(-170,120)(1,0)[F(A)`G(A)`\alpha(A)]{540}1a
 \putmorphism(360,120)(1,0)[\phantom{G(B)}`G(A) `G(f)]{560}1a
\putmorphism(-180,500)(0,-1)[\phantom{Y_2}``=]{380}1r
\putmorphism(940,500)(0,-1)[\phantom{Y_2}``=]{380}1r
\put(280,310){\fbox{$\delta_{\alpha,f}$}}
\efig}
$$
\end{enumerate}
so that the following are satisfied:
\begin{itemize}
\item (coherence with composition for $\delta_{\alpha,-}$): for any composable 1h-cells $f$ and $g$ in $\Aa$ the 2-cell 
$\delta_{\alpha,gf}$ satisfies: \\
{\em (h.o.t.-1)} \label{h.o.t.-1} 
$$\scalebox{0.86}{
\bfig
  \putmorphism(-750,500)(1,0)[F(A)`\phantom{F(B)}`F(f)]{600}1a
\putmorphism(-130,500)(1,0)[F(A)`F(C)`F(g)]{580}1a
\putmorphism(-730,500)(0,-1)[\phantom{Y_2}``=]{400}1r
\putmorphism(420,500)(0,-1)[\phantom{Y_2}``=]{400}1r
\putmorphism(-730,90)(0,-1)[\phantom{Y_2}``=]{400}1r
\putmorphism(1030,90)(0,-1)[\phantom{Y_2}``=]{400}1r
 \putmorphism(450,90)(1,0)[F(C)`G(C) `\alpha(C)]{580}1a
  \putmorphism(-750,90)(1,0)[F(A)`\phantom{F(B)}`F(gf)]{1200}1a
\put(-270,320){\fbox{$F_{gf}$}}
 \putmorphism(-750,-315)(1,0)[F(A)`G(A)`\alpha(A)]{620}1a
 \putmorphism(-120,-315)(1,0)[\phantom{F(B)}`G(C) `G(gf)]{1170}1a
\put(-250,-120){\fbox{$\delta_{\alpha,gf}$}}
\efig}= 
\scalebox{0.86}{
\bfig
 \putmorphism(450,450)(1,0)[F(B)`F(C) `F(g)]{680}1a
 \putmorphism(1120,450)(1,0)[\phantom{F(B)}`G(C) `\alpha(C)]{600}1a
\put(1000,200){\fbox{$\delta_{\alpha,g}$}}

  \putmorphism(-150,0)(1,0)[F(A)` F(B) `F(f)]{600}1a
\putmorphism(450,0)(1,0)[\phantom{F(A)}` G(B) `\alpha(B)]{680}1a
 \putmorphism(1120,0)(1,0)[\phantom{F(A)}`G(C) ` G(g)]{620}1a

\putmorphism(450,450)(0,-1)[\phantom{Y_2}``=]{450}1l
\putmorphism(1710,450)(0,-1)[\phantom{Y_2}``=]{450}1r

 \putmorphism(-150,-450)(1,0)[F(A)`G(A)`\alpha(A)]{600}1a
 \putmorphism(450,-450)(1,0)[\phantom{F(B)}`G(B) `G(f)]{680}1a
 \putmorphism(1120,-450)(1,0)[\phantom{F(B)}`G(C) `G(g)]{620}1a

\putmorphism(-180,0)(0,-1)[\phantom{Y_2}``=]{450}1r
\putmorphism(1040,0)(0,-1)[\phantom{Y_2}``=]{450}1r
\put(350,-240){\fbox{$\delta_{\alpha,f}$}}
\put(1000,-660){\fbox{$G_{gf}$}}

 \putmorphism(450,-900)(1,0)[G(A)` G(C) `G(gf)]{1300}1a

\putmorphism(450,-450)(0,-1)[\phantom{Y_2}``=]{450}1l
\putmorphism(1750,-450)(0,-1)[\phantom{Y_2}``=]{450}1r
\efig}
$$ 
(coherence with identity for $\delta_{\alpha,-}$):
$$\text{{\em (h.o.t.-2)}} \label{h.o.t.-2} \qquad\quad
\scalebox{0.86}{
\bfig
 \putmorphism(-150,420)(1,0)[F(A)`F(A)`=]{500}1a
\putmorphism(-180,420)(0,-1)[\phantom{Y_2}``=]{370}1l
\putmorphism(320,420)(0,-1)[\phantom{Y_2}``=]{370}1r
 \putmorphism(-150,50)(1,0)[F(A)`F(A)`F(1_A)]{500}1a
 \put(-80,250){\fbox{$F_A$}} 
\putmorphism(330,50)(1,0)[\phantom{F(A)}`G(A) `\alpha(A)]{560}1a
 \putmorphism(-170,-350)(1,0)[F(A)`G(A)`\alpha(A)]{520}1a
 \putmorphism(350,-350)(1,0)[\phantom{F(A)}`G(A) `G(1_A)]{560}1a

\putmorphism(-180,50)(0,-1)[\phantom{Y_2}``=]{400}1l
\putmorphism(910,50)(0,-1)[\phantom{Y_2}``=]{400}1r
\put(240,-150){\fbox{$\delta_{\alpha,1_A}$}}
\efig}
=
\scalebox{0.86}{
\bfig
 \putmorphism(-150,420)(1,0)[F(A)`G(A)`\alpha(A)]{500}1a
\putmorphism(-180,420)(0,-1)[\phantom{Y_2}``=]{370}1l
\putmorphism(320,420)(0,-1)[\phantom{Y_2}``=]{370}1r
  \put(-100,230){\fbox{$\Id_{\alpha(A)}$}} 
\putmorphism(-150,50)(1,0)[F(A)` \phantom{Y_2} `\alpha(A)]{450}1a

\putmorphism(350,50)(1,0)[G(A)` G(A) `=]{470}1a
\putmorphism(330,-300)(1,0)[G(A)` G(A) `G(1_A)]{480}1b
\putmorphism(330,50)(0,-1)[\phantom{Y_2}``=]{350}1l
\putmorphism(800,50)(0,-1)[\phantom{Y_2}``=]{350}1r
\put(470,-150){\fbox{$G_A$}}
\efig}
$$

\item (coherence with composition and identity for $\alpha_\bullet$): for any composable 1v-cells $u$ and $v$ in $\Aa$:
$$\text{{\em (h.o.t.-3)}} \label{h.o.t.-3} \qquad\alpha^{\frac{u}{v}}=\frac{\alpha^u}{\alpha^v}\quad\qquad\text{ and}\quad\qquad
\text{{\em (h.o.t.-4)}} \label{h.o.t.-4} \qquad\alpha^{1^A}=\Id_{\alpha(A)};$$

\item (oplax naturality of 2-cells):
for every 2-cell in $\Aa$
$\scalebox{0.86}{
\bfig
\putmorphism(-150,50)(1,0)[A` B`f]{400}1a
\putmorphism(-150,-270)(1,0)[A'`B' `g]{400}1b
\putmorphism(-170,50)(0,-1)[\phantom{Y_2}``u]{320}1l
\putmorphism(250,50)(0,-1)[\phantom{Y_2}``v]{320}1r
\put(0,-140){\fbox{$a$}}
\efig}$ 
the following identity in $\Bb$ must hold:\\
$\text{{\em (h.o.t.-5)}}$ \label{h.o.t.-5} 
$$
\scalebox{0.86}{
\bfig
\putmorphism(-150,500)(1,0)[F(A)`F(B)`F(f)]{600}1a
 \putmorphism(450,500)(1,0)[\phantom{F(A)}`G(B) `\alpha(B)]{640}1a

 \putmorphism(-150,50)(1,0)[F(A')`F(B')`F(g)]{600}1a
 \putmorphism(450,50)(1,0)[\phantom{F(A)}`G(B') `\alpha(B')]{640}1a

\putmorphism(-180,500)(0,-1)[\phantom{Y_2}``F(u)]{450}1l
\putmorphism(450,500)(0,-1)[\phantom{Y_2}``]{450}1r
\putmorphism(300,500)(0,-1)[\phantom{Y_2}``F(v)]{450}0r
\putmorphism(1100,500)(0,-1)[\phantom{Y_2}``G(v)]{450}1r
\put(0,260){\fbox{$F(a)$}}
\put(700,270){\fbox{$\alpha^v$}}

\putmorphism(-150,-400)(1,0)[F(A')`G(A') `\alpha(A')]{640}1a
 \putmorphism(450,-400)(1,0)[\phantom{A'\ot B'}` G(B') `G(g)]{680}1a

\putmorphism(-180,50)(0,-1)[\phantom{Y_2}``=]{450}1l
\putmorphism(1120,50)(0,-1)[\phantom{Y_3}``=]{450}1r
\put(320,-200){\fbox{$\delta_{\alpha,g}$}}

\efig}
\quad=\quad
\scalebox{0.86}{
\bfig
\putmorphism(-150,500)(1,0)[F(A)`F(B)`F(f)]{600}1a
 \putmorphism(450,500)(1,0)[\phantom{F(A)}`G(B) `\alpha(B)]{680}1a
 \putmorphism(-150,50)(1,0)[F(A)`G(A)`\alpha(A)]{600}1a
 \putmorphism(450,50)(1,0)[\phantom{F(A)}`G(B) `G(f)]{680}1a

\putmorphism(-180,500)(0,-1)[\phantom{Y_2}``=]{450}1r
\putmorphism(1100,500)(0,-1)[\phantom{Y_2}``=]{450}1r
\put(350,260){\fbox{$\delta_{\alpha,f}$}}
\put(650,-180){\fbox{$G(a)$}}

\putmorphism(-150,-400)(1,0)[F(A')`G(A') `\alpha(A')]{640}1a
 \putmorphism(490,-400)(1,0)[\phantom{F(A')}` G(B'). `G(g)]{640}1a

\putmorphism(-180,50)(0,-1)[\phantom{Y_2}``F(u)]{450}1l
\putmorphism(450,50)(0,-1)[\phantom{Y_2}``]{450}1l
\putmorphism(610,50)(0,-1)[\phantom{Y_2}``G(u)]{450}0l 
\putmorphism(1120,50)(0,-1)[\phantom{Y_3}``G(v)]{450}1r
\put(40,-180){\fbox{$\alpha^u$}} 
\efig}
$$
\end{itemize}
A {\em horizontal strict transformation} is a  horizontal oplax transformation for which the 2-cells $\delta_{\alpha,f}$ in item 3. are identities. 
\end{defn}


The lax version of the above Definition we will need in \coref{obtained oplax trans}. 
A {\em horizontal lax transformation} $\alpha$ between lax double functors $F,G: \Aa\to\Bb$ differs from its oplax counterpart 
in that the globular 2-cells $\delta_{\alpha,f}$ for any 1h-cell $f$ in $\Aa$ goes in the other direction, namely 
$$
\scalebox{0.86}{
\bfig
 \putmorphism(-170,500)(1,0)[F(A)`G(A)`\alpha(A)]{540}1a
 \putmorphism(360,500)(1,0)[\phantom{G(B)}`G(A) `G(f)]{560}1a
 \putmorphism(-170,120)(1,0)[F(A)`F(B)`F(f)]{540}1a
 \putmorphism(360,120)(1,0)[\phantom{F(f)}`G(B) `\alpha(B)]{560}1a
\putmorphism(-180,500)(0,-1)[\phantom{Y_2}``=]{380}1r
\putmorphism(940,500)(0,-1)[\phantom{Y_2}``=]{380}1r
\put(280,310){\fbox{$\sigma_{\alpha,f}$}}
\efig}
$$
and the axioms (h.o.t.-1)-(h.o.t.-5) are accordingly changed by the analogous axioms that we will refer to as to 
(h.l.t.-1)-(h.l.t.-5). Indeed, note that only the two axioms (h.o.t.-1), (h.o.t.-2) and (h.o.t.-5) are changed into (h.l.t.-1), 
(h.l.t.-2) and (h.l.t.-5).

The above two Definitions are ``oplax, respectively lax, and horizontal'' versions of a ``strong vertical transformation'' from 
\cite[Section 7.4]{GP:L} for strict (rather then pseudo) double categories. Similarly, the following is a horizontal version of a 
``strong modification'' from {\em loc. cit.} with $H$ and $K$ being identities.

The composition of 1h-cells in $\Lax_{hop}(\Aa, \Bb)$, that is of horizontal oplax transformations $\alpha$ and 
$\beta$ acting between lax double functors $F,G,H:\Aa\to\Bb$, is given by the vertical composition of transformations, 
that we make explicit here:

\begin{lma} \lelabel{vert comp hor.ps.tr.}
Vertical composition of two horizontal oplax transformations 
$F\stackrel{\alpha}{\Rightarrow} G \stackrel{\beta}{\Rightarrow}H$ between lax functors $F, G, H:\Aa\to\Bb$, 
denoted by $\frac{\alpha}{\beta}$, is well-given by: 
\begin{enumerate}
\item for every 0-cell $A$ in $\Aa$ a 1h-cell in $\Bb$:
$$(\frac{\alpha}{\beta})(A)=\big( F(A)\stackrel{\alpha(A)}{\longrightarrow}G(A) \stackrel{\beta(A)}{\longrightarrow} H(A) \big),$$ 
\item for every 1v-cell $u:A\to A'$ in $\Aa$ a 2-cell in $\Bb$:
$$(\frac{\alpha}{\beta})^u=
\bfig
\putmorphism(-150,50)(1,0)[F(A)`G(A)`\alpha(A)]{560}1a
 \putmorphism(430,50)(1,0)[\phantom{F(A)}`H(A) `\beta(A)]{620}1a

\putmorphism(-150,-320)(1,0)[F(A')`G(A')`\alpha(A')]{600}1a
\putmorphism(-180,50)(0,-1)[\phantom{Y_2}``F(u)]{370}1l
\putmorphism(380,50)(0,-1)[\phantom{Y_2}``G(u)]{370}1r
\put(-30,-130){\fbox{$\alpha^u$}}

 \putmorphism(460,-320)(1,0)[\phantom{F(A)}`H(A') `\beta(A')]{630}1a

\putmorphism(1060,50)(0,-1)[\phantom{Y_2}``H(u)]{370}1r
\put(640,-110){\fbox{$\beta^u$}}
\efig
$$
\item for every 1h-cell $f:A\to B$  in $\Aa$ a 2-cell in $\Bb$:
$$\delta_{\frac{\alpha}{\beta},f}=
\bfig

 \putmorphism(-200,-50)(1,0)[F(A)`F(B)` F(f)]{650}1a
 \putmorphism(430,-50)(1,0)[\phantom{F(B)}`G(B) `\alpha(B)]{700}1a

 \putmorphism(-200,-450)(1,0)[F(A)`G(A)`\alpha(A)]{650}1a
 \putmorphism(430,-450)(1,0)[\phantom{A\ot B}`G(B) `G(f)]{700}1a
 \putmorphism(1050,-450)(1,0)[\phantom{A'\ot B'}` H(B) `\beta(B)]{700}1a

\putmorphism(-230,-50)(0,-1)[\phantom{Y_2}``=]{400}1r
\putmorphism(1050,-50)(0,-1)[\phantom{Y_2}``=]{400}1r
\put(300,-240){\fbox{$ \delta_{\alpha,f}  $}}
\put(1000,-660){\fbox{$\delta_{\beta,f}$}}

 \putmorphism(450,-850)(1,0)[G(A)` H(A) `\beta(A)]{700}1a
 \putmorphism(1080,-850)(1,0)[\phantom{A''\ot B'}`G(B). ` H(f)]{700}1a

\putmorphism(450,-450)(0,-1)[\phantom{Y_2}``=]{400}1l
\putmorphism(1750,-450)(0,-1)[\phantom{Y_2}``=]{400}1r
\efig
$$
\end{enumerate}
\end{lma}

\begin{proof}
In \cite[Lemma 3.7]{Fem} we proved that the vertical composition of two horizontal {\em pseudonatural} transformations 
between {\em double pseudo} functors is given in the same way as in the statement of the present Lemma. 
For the purpose of the present setting, for horizontal {\em oplax} transformations between {\em lax} functors, 
we have checked that the same holds, in the exactly same way: the proof does not rely on the nature of the coherence 2-cells 
of the transformations, nor of the double functors in question. 
\qed\end{proof}

Since the vertical composition of horizontal oplax transformations is defined in terms of the horizontal composition in $\Bb$, 
it is strictly associative. From here one also sees that $\Lax_{hop}(\Aa,\Bb)$ as a to-be-constricted double category is strict.

\smallskip

\begin{defn}
A {\em vertical lax transformation} $\alpha_0$ between lax double functors $F,G: \Aa\to\Bb$ consists of: 
\begin{enumerate}
\item a 1v-cell $\alpha_0(A):F(A)\to G(A)$ in $\Bb$ for every 0-cell $A$ in $\Aa$; 
\item 
for every 1h-cell $f:A\to B$ in $\Aa$ a 2-cell in $\Bb$:
$$ 
\bfig
\putmorphism(-150,180)(1,0)[F(A)`F(B)`F(f)]{560}1a
\putmorphism(-150,-190)(1,0)[G(A)`G(B)`G(f)]{600}1a
\putmorphism(-180,180)(0,-1)[\phantom{Y_2}``\alpha_0(A)]{370}1l
\putmorphism(410,180)(0,-1)[\phantom{Y_2}``\alpha_0(B)]{370}1r
\put(-10,20){\fbox{$(\alpha_0)_f$}}
\efig
$$
\item for every 1v-cell $u:A\to A'$ in $\Aa$ a 2-cell in $\Bb$:
$$
\bfig
 \putmorphism(-90,500)(1,0)[F(A)`F(A) `=]{540}1a
\putmorphism(450,500)(0,-1)[\phantom{Y_2}`F(\tilde A) `F(u)]{400}1r
\putmorphism(-90,-300)(1,0)[G(\tilde A)`G(\tilde A) `=]{540}1a
\putmorphism(450,100)(0,-1)[\phantom{Y_2}``\alpha_0(\tilde A)]{400}1r
\putmorphism(-120,100)(0,-1)[\phantom{Y_2}``G(u)]{400}1l
\putmorphism(-120,500)(0,-1)[\phantom{Y_2}`G(A) `\alpha_0(A)]{400}1l
\put(60,50){\fbox{$\alpha_0^u$}}
\efig
$$
\end{enumerate}
which need to satisfy: 
\begin{itemize}
\item (coherence with composition for $(\alpha_0)_\bullet$):  \\
$$\text{{\em (v.l.t.\x 1)}} \label{v.l.t.-1} \qquad
\scalebox{0.86}{
\bfig
\putmorphism(-150,500)(1,0)[F(A)`F(B)`F(f)]{600}1a
 \putmorphism(450,500)(1,0)[\phantom{F(A)}`F(C) `F(g)]{640}1a

 \putmorphism(-150,50)(1,0)[G(A)`G(B)`G(f)]{600}1a
 \putmorphism(450,50)(1,0)[\phantom{F(A)}`G(C) `G(g)]{640}1a

\putmorphism(-180,500)(0,-1)[\phantom{Y_2}``\alpha_0(A)]{450}1l
\putmorphism(450,500)(0,-1)[\phantom{Y_2}``]{450}1r
\putmorphism(300,500)(0,-1)[\phantom{Y_2}``\alpha_0(B)]{450}0r
\putmorphism(1100,500)(0,-1)[\phantom{Y_2}``\alpha_0(C)]{450}1r
\put(0,260){\fbox{$(\alpha_0)_f$}}
\put(700,270){\fbox{$(\alpha_0)_g$}}

\putmorphism(-150,-400)(1,0)[G(A)`G(C) `G(gf)]{1260}1a

\putmorphism(-180,50)(0,-1)[\phantom{Y_2}``=]{450}1l
\putmorphism(1120,50)(0,-1)[\phantom{Y_3}``=]{450}1r
\put(320,-170){\fbox{$G_{gf}$}}

\efig}
=
\scalebox{0.86}{
\bfig
\putmorphism(-150,500)(1,0)[F(A)`F(B)`F(f)]{600}1a
 \putmorphism(450,500)(1,0)[\phantom{F(A)}`F(C) `F(g)]{640}1a

 \putmorphism(-150,50)(1,0)[F(A)`F(C)`F(gf)]{1260}1a

\putmorphism(-180,500)(0,-1)[\phantom{Y_2}``=]{450}1r
\putmorphism(1100,500)(0,-1)[\phantom{Y_2}``=]{450}1r
\put(350,260){\fbox{$F_{gf}$}}
\put(350,-180){\fbox{$(\alpha_0)_{gf}$}}

\putmorphism(-180,50)(0,-1)[\phantom{Y_2}``\alpha_0(A)]{450}1l
\putmorphism(1120,50)(0,-1)[\phantom{Y_3}``\alpha_0(C)]{450}1r
\putmorphism(-150,-400)(1,0)[G(A)`G(C) `G(gf)]{1260}1a
\efig}
$$

(coherence with identity for $(\alpha_0)_\bullet$):
$$\text{{\em (v.l.t.\x 2)}} \label{v.l.t.-2} \qquad\qquad 
\scalebox{0.86}{
\bfig
 \putmorphism(-170,420)(1,0)[F(A)`F(A)`=]{500}1a
\putmorphism(-180,420)(0,-1)[\phantom{Y_2}``=]{370}1l
\putmorphism(280,420)(0,-1)[\phantom{Y_2}``=]{370}1r
  \put(-40,250){\fbox{$F_A$}} 
\putmorphism(-170,50)(1,0)[F(A)` F(A) `F(1_A)]{450}1a

\putmorphism(-170,50)(0,-1)[\phantom{Y_2}``\alpha_0(A)]{350}1l
\putmorphism(280,50)(0,-1)[\phantom{Y_2}``\alpha_0(A)]{350}1r
\putmorphism(-150,-300)(1,0)[G(A)` G(A) `G(1_A)]{440}1b
\put(-70,-140){\fbox{$(\alpha_0)_{1_A}$}}
\efig}
=
\scalebox{0.86}{
\bfig
 \putmorphism(-170,420)(1,0)[F(A)`F(A)`=]{500}1a
\putmorphism(-180,420)(0,-1)[\phantom{Y_2}``\alpha_0(A)]{370}1l
\putmorphism(280,420)(0,-1)[\phantom{Y_2}``\alpha_0(A)]{370}1r
  \put(-80,220){\fbox{$\Id^{\alpha_0(A)}$}} 
\putmorphism(-170,50)(1,0)[G(A)` G(A) `=]{450}1a

\putmorphism(-170,50)(0,-1)[\phantom{Y_2}``=]{350}1l
\putmorphism(280,50)(0,-1)[\phantom{Y_2}``=]{350}1r
\putmorphism(-150,-300)(1,0)[G(A)` G(A) `G(1_A)]{440}1b
\put(-40,-140){\fbox{$G_A$}}
\efig}
$$
\item (coherence with composition for $\alpha_0^\bullet$):
$$ \text{{\em (v.l.t.\x 3)}} \label{v.l.t.-3} \qquad\quad 
\bfig 
 \putmorphism(-150,500)(1,0)[F(A)`F(A) `=]{500}1a
\putmorphism(-130,500)(0,-1)[\phantom{Y_2}`G(A) `\alpha_0(A)]{400}1l
\put(0,250){\fbox{$\alpha_0^u$}}
\putmorphism(-150,-300)(1,0)[G(\tilde A)`G(\tilde A) `=]{460}1a
\putmorphism(-130,110)(0,-1)[\phantom{Y_2}``G(u)]{400}1l
\putmorphism(380,500)(0,-1)[\phantom{Y_2}` `F(u)]{400}1r
\putmorphism(380,100)(0,-1)[F(\tilde A)` `\alpha_0(\tilde A)]{400}1l
\putmorphism(480,110)(1,0)[`F(\tilde A)`=]{460}1a
\putmorphism(390,-700)(1,0)[\phantom{G(A)}`G(\tilde{\tilde A})`=]{570}1a
\putmorphism(920,100)(0,-1)[\phantom{(B, \tilde A')}``F(u')]{400}1r
\putmorphism(920,-300)(0,-1)[F(\tilde{\tilde A})`` \alpha_0(\tilde{\tilde A})]{400}1r
\putmorphism(400,-300)(0,-1)[\phantom{(B, \tilde A)}`G(\tilde{\tilde A}) `G(u')]{400}1l
\put(530,-190){\fbox{$\alpha_0^{u'}$}}
\efig
=
\alpha_0^{\frac{u}{u'}}
$$
(coherence with identity for $\alpha^\bullet$):\\
$$\text{{\em (v.l.t.\x 4)}} \label{v.l.t.-4} \qquad\qquad \alpha_0^{1^A}=\Id^{\alpha_0(A)}  \hspace{4cm}$$
\item (lax naturality of 2-cells):
for every 2-cell in $\Aa$
$\scalebox{0.86}{
\bfig
\putmorphism(-150,50)(1,0)[A` B`f]{400}1a
\putmorphism(-150,-270)(1,0)[A'`B' `g]{400}1b
\putmorphism(-170,50)(0,-1)[\phantom{Y_2}``u]{320}1l
\putmorphism(250,50)(0,-1)[\phantom{Y_2}``v]{320}1r
\put(0,-140){\fbox{$a$}}
\efig}$ 
the following identity in $\Bb$ must hold:\\
$$\text{{\em (v.l.t.\x 5)}} \label{v.l.t.-5} \qquad
\scalebox{0.86}{
\bfig
 \putmorphism(-130,500)(1,0)[F(A)`F(A) `=]{560}1a
 \putmorphism(530,500)(1,0)[` `F(f)]{410}1a
\putmorphism(-130,500)(0,-1)[\phantom{Y_2}`G(A) `\alpha_0(A)]{450}1l
\put(40,50){\fbox{$\alpha_0^u$}}
\putmorphism(-130,-400)(1,0)[G(\tilde A)` `=]{500}1a
\putmorphism(-150,50)(0,-1)[\phantom{Y_2}``G(u)]{450}1l
\putmorphism(450,50)(0,-1)[\phantom{Y_2}`G(\tilde A)`\alpha_0(\tilde A)]{450}1l
\putmorphism(450,500)(0,-1)[\phantom{Y_2}`F(\tilde A) `F(u)]{450}1l
\put(600,260){\fbox{$F(a)$}}
\putmorphism(420,50)(1,0)[\phantom{(B, \tilde A)}``F(g)]{500}1a
\putmorphism(1030,50)(0,-1)[\phantom{(B, A')}`G(\tilde A)`\alpha_0(\tilde B)]{450}1r
\putmorphism(1030,500)(0,-1)[F(B)`F(\tilde B)`F(v)]{450}1r
\putmorphism(420,-400)(1,0)[\phantom{(B, \tilde A)}``G(g)]{500}1a
\put(600,-170){\fbox{$(\alpha_0)_g$}}
\efig}=
\scalebox{0.86}{
\bfig
 \putmorphism(-130,500)(1,0)[F(A)`F(B) `F(f)]{560}1a
 \putmorphism(550,500)(1,0)[` `=]{410}1a
\putmorphism(-180,500)(0,-1)[\phantom{Y_2}`G(A) `\alpha_0(A)]{450}1l
\put(650,50){\fbox{$\alpha_0^v$}}
\putmorphism(-130,-400)(1,0)[G(\tilde A)` `G(g)]{500}1a
\putmorphism(-180,50)(0,-1)[\phantom{Y_2}``G(u)]{450}1l
\putmorphism(450,50)(0,-1)[\phantom{Y_2}`G(\tilde B)`G(v)]{450}1r
\putmorphism(450,500)(0,-1)[\phantom{Y_2}`G(B) `\alpha_0(B)]{450}1r
\put(40,280){\fbox{$(\alpha_0)_f$}}
\putmorphism(-150,50)(1,0)[\phantom{(B, \tilde A)}``G(f)]{500}1a
\putmorphism(1050,50)(0,-1)[\phantom{(B, A')}`G(\tilde B).`\alpha_0(\tilde B)]{450}1r
\putmorphism(1050,500)(0,-1)[F(B)`F(\tilde B)`F(v)]{450}1r
\putmorphism(450,-400)(1,0)[\phantom{(B, \tilde A)}``=]{500}1b
\put(0,-170){\fbox{$G(a)$}}
\efig}
$$
\end{itemize}
\end{defn}

The definition of vertical composition of two vertical lax transformations of lax double functors, that we give in the next Lemma, 
is the same - up to the orientation of the coherence 2-cells - as that in \cite[Lemma 3.8]{Fem} for vertical pseudonatural transformations 
of double pseudofunctors. Identically as in the proof of \leref{vert comp hor.ps.tr.}, the proof of 
well-definedness is direct and does not depend on the coherence structures of double functors and their transformations.

\begin{lma} \lelabel{vert comp vert. lx tr.}
Vertical composition of two vertical lax transformations $\alpha_0: F\Rightarrow G: \Aa\to\Bb$ and 
$\beta_0: G\Rightarrow H:\Aa\to\Bb$ between lax double functors, denoted by $\frac{\alpha_0}{\beta_0}$, is well-given by: 
\begin{itemize}
\item for every 0-cell $A$ in $\Aa$ a 1v-cell on the left below, and for every 1h-cell $f:A\to B$ in $\Aa$ a 2-cell on the right below, both in $\Bb$:
$$(\frac{\alpha_0}{\beta_0})(A)=
\bfig
\putmorphism(-280,500)(0,-1)[F(A)`G(A) `\alpha_0(A)]{450}1l
 \putmorphism(-280,70)(0,-1)[\phantom{F(A)}`H(A) `\beta_0(A)]{450}1l
\efig \qquad
(\frac{\alpha_0}{\beta_0})(f)= 
\bfig
\putmorphism(-250,500)(1,0)[F(A)`F(B)` F(f)]{550}1a
 \putmorphism(-250,50)(1,0)[G(A)`G(B)` G(f)]{550}1a
 \putmorphism(-250,-400)(1,0)[H(A)`H(B)` H(f)]{550}1a

\putmorphism(-280,500)(0,-1)[\phantom{Y_2}``\alpha_0(A)]{450}1l
 \putmorphism(-280,70)(0,-1)[\phantom{F(A)}` `\beta_0(A)]{450}1l

\putmorphism(300,500)(0,-1)[\phantom{Y_2}``\alpha_0(B)]{450}1r
\putmorphism(300,70)(0,-1)[\phantom{Y_2}``\beta_0(B)]{450}1r
\put(-120,290){\fbox{$(\alpha_0)_f$}}
\put(-120,-150){\fbox{$(\beta_0)_f$}}
\efig
$$ 
\item for every 1v-cell $u:A\to A'$ in $\Aa$ a 2-cell in $\Bb$: 
$$(\frac{\alpha_0}{\beta_0})^u=
\bfig
 \putmorphism(920,450)(1,0)[F(A)`F(A) `=]{460}1a
\putmorphism(920,450)(0,-1)[\phantom{Y_2}`F(A') `\alpha_0(A)]{400}1l
\put(1060,250){\fbox{$\alpha_0^u$}}
\putmorphism(920,-400)(1,0)[G(A')`G(A') `=]{460}1a
\putmorphism(1380,450)(0,-1)[\phantom{Y_2}` `F(u)]{400}1r
\putmorphism(1380,50)(0,-1)[F(\tilde A)``\alpha_0(\tilde A)]{450}1r
\putmorphism(440,50)(0,-1)[F(A')`F(A'') `\beta_0(A)]{450}1l
\putmorphism(530,60)(1,0)[``=]{300}1a
\putmorphism(560,-850)(1,0)[`G(A'')`=]{430}1a
\putmorphism(920,50)(0,-1)[\phantom{(B, \tilde A')}``G(u)]{450}1r
\putmorphism(920,-400)(0,-1)[`` \beta_0(\tilde A)]{450}1r
\putmorphism(440,-400)(0,-1)[\phantom{(B, \tilde A)}`G(A'') `H(u)]{450}1l
\put(560,-190){\fbox{$\beta_0^u$}}
\efig
$$
\end{itemize}
\end{lma}

We finally define 2-cells for the double category $\Lax_{hop}(\Aa, \Bb)$.

\begin{defn} \delabel{modif-hv}
A modification $\Theta$ between two horizontal oplax transformations $\alpha$ and $\beta$ and two vertical lax transformations $\alpha_0$ and $\beta_0$ 
depicted below on the left, where the lax double functors $F, G, F\s', G'$ act between $\Aa\to\Bb$, is given 
by a collection of 2-cells in $\Bb$ depicted below on the right:
\begin{equation} \eqlabel{modification cells}
\scalebox{0.86}{
\bfig
\putmorphism(-150,50)(1,0)[F` G`\alpha]{400}1a
\putmorphism(-150,-270)(1,0)[F'`G' `\beta]{400}1b
\putmorphism(-170,50)(0,-1)[\phantom{Y_2}``\alpha_0]{320}1l
\putmorphism(250,50)(0,-1)[\phantom{Y_2}``\beta_0]{320}1r
\put(-30,-140){\fbox{$\Theta$}}
\efig}
\qquad\qquad
\scalebox{0.86}{
\bfig
\putmorphism(-180,50)(1,0)[F(A)` G(A)`\alpha(A)]{550}1a
\putmorphism(-180,-270)(1,0)[F\s'(A)`G'(A) `\beta(A)]{550}1b
\putmorphism(-170,50)(0,-1)[\phantom{Y_2}``\alpha_0(A)]{320}1l
\putmorphism(350,50)(0,-1)[\phantom{Y_2}``\beta_0(A)]{320}1r
\put(0,-140){\fbox{$\Theta_A$}}
\efig}
\end{equation}
which satisfy the following rules: 

\medskip

{\em (m.ho-vl.-1)} \label{m.ho-vl.-1} for every 1h-cell $f$  
$$
\scalebox{0.86}{
\bfig
\putmorphism(-150,500)(1,0)[F(A)`F(B)`F(f)]{600}1a
 \putmorphism(450,500)(1,0)[\phantom{F(A)}`G(B) `\alpha(B)]{620}1a

 \putmorphism(-150,50)(1,0)[F\s'(A)`F\s'(B)`F\s'(f)]{600}1a
 \putmorphism(450,50)(1,0)[\phantom{F(A)}`G'(B) `\beta(B)]{620}1a

\putmorphism(-180,500)(0,-1)[\phantom{Y_2}``\alpha_0(A)]{450}1l
\putmorphism(450,500)(0,-1)[\phantom{Y_2}``]{450}1r
\putmorphism(300,500)(0,-1)[\phantom{Y_2}``\alpha_0(B)]{450}0r
\putmorphism(1080,500)(0,-1)[\phantom{Y_2}``\beta_0(B)]{450}1r
\put(0,280){\fbox{$(\alpha_0)_f$}}
\put(670,280){\fbox{$\Theta_B$}}

\putmorphism(-150,-400)(1,0)[F\s'(A)`G'(A) `\beta(A)]{600}1a
\putmorphism(510,-400)(1,0)[\phantom{Y_2}`G'(B) `G'(f)]{580}1a

\putmorphism(-180,50)(0,-1)[\phantom{Y_2}``=]{450}1l
\putmorphism(1080,50)(0,-1)[\phantom{Y_3}``=]{450}1r
\put(320,-180){\fbox{$\delta_{\beta,f}$}}

\efig}
\quad=\quad
\scalebox{0.86}{
\bfig
\putmorphism(-150,500)(1,0)[F(A)`F(B)`F(f)]{600}1a
 \putmorphism(450,500)(1,0)[\phantom{F(A)}`G(B) `\alpha(B)]{620}1a
\putmorphism(-150,50)(1,0)[F(A)`G(A) `\alpha(A)]{600}1a
\putmorphism(510,50)(1,0)[\phantom{Y_2}`G(B) `G(f)]{580}1a

\putmorphism(-180,500)(0,-1)[\phantom{Y_2}``=]{450}1r
\putmorphism(1080,500)(0,-1)[\phantom{Y_2}``=]{450}1r
\put(350,280){\fbox{$\delta_{\alpha,f}$}}

\putmorphism(-180,50)(0,-1)[\phantom{Y_2}``\alpha_0(A)]{450}1l
\putmorphism(1080,50)(0,-1)[\phantom{Y_3}``\beta_0(B)]{450}1r
\put(20,-180){\fbox{$\Theta_A$}}
\put(670,-180){\fbox{$(\beta_0)_f$}}

\putmorphism(450,50)(0,-1)[\phantom{Y_2}``]{450}1r
\putmorphism(300,50)(0,-1)[\phantom{Y_2}``\beta_0(A)]{450}0r

\putmorphism(-150,-400)(1,0)[F\s'(A)`G'(A) `\beta(A)]{600}1a
\putmorphism(510,-400)(1,0)[\phantom{Y_2}`G'(B) `G'(f)]{580}1a
\efig}
$$
and

{\em (m.ho-vl.-2)} \label{m.ho-vl.-2} for every 1v-cell $u$  
$$
\scalebox{0.86}{
\bfig
 \putmorphism(-150,500)(1,0)[F(A)`F(A) `=]{600}1a
 \putmorphism(550,500)(1,0)[` `\alpha(A)]{400}1a
\putmorphism(-180,500)(0,-1)[\phantom{Y_2}`F\s'(A) `\alpha_0(A)]{450}1l
\put(30,50){\fbox{$\alpha_0^u$}}
\putmorphism(-150,-400)(1,0)[F\s'(\tilde A)` `=]{500}1a
\putmorphism(-180,50)(0,-1)[\phantom{Y_2}``F\s'(u)]{450}1l
\putmorphism(450,50)(0,-1)[\phantom{Y_2}`F\s'(\tilde A)` \alpha_0(\tilde A)]{450}1l
\putmorphism(450,500)(0,-1)[\phantom{Y_2}`F(\tilde A) `F(u)]{450}1l
\put(660,280){\fbox{$\alpha^u$}}
\putmorphism(450,50)(1,0)[\phantom{(B, \tilde A)}``\alpha(\tilde A)]{500}1a
\putmorphism(1070,50)(0,-1)[\phantom{(B, A')}`G'(\tilde A)`\beta_0(\tilde A)]{450}1r
\putmorphism(1070,500)(0,-1)[G(A)`G(\tilde A)`G(u)]{450}1r
\putmorphism(450,-400)(1,0)[\phantom{(B, \tilde A)}``\beta(\tilde A)]{500}1a
\put(640,-170){\fbox{$ \Theta_{\tilde A}$ } } 
\efig}=
\scalebox{0.86}{
\bfig
 \putmorphism(-150,500)(1,0)[F(A)`G(A) `\alpha(A)]{600}1a
 \putmorphism(450,500)(1,0)[\phantom{(B,A)}` `=]{450}1a
\putmorphism(-180,500)(0,-1)[\phantom{Y_2}`F\s'(A) `\alpha_0(A)]{450}1l
\put(620,50){\fbox{$\beta_0^u$}}
\putmorphism(-150,-400)(1,0)[F\s'(\tilde A)` `\beta(\tilde A)]{500}1a
\putmorphism(-180,50)(0,-1)[\phantom{Y_2}``F\s'(u)]{450}1l
\putmorphism(450,50)(0,-1)[\phantom{Y_2}`G'(\tilde A)`G'(u)]{450}1r
\putmorphism(450,500)(0,-1)[\phantom{Y_2}`G'(A) `\beta_0(A)]{450}1r
\put(0,260){\fbox{$\Theta_A$}}
\putmorphism(-150,50)(1,0)[\phantom{(B, \tilde A)}``\beta(A)]{500}1a
\putmorphism(1030,50)(0,-1)[\phantom{(B, A')}` G'(\tilde A). ` \beta_0(\tilde A)]{450}1r
\putmorphism(1030,500)(0,-1)[G(A)`G(\tilde A)` G(u)]{450}1r
\putmorphism(450,-400)(1,0)[\phantom{(B, \tilde A)}``=]{500}1b
\put(70,-170){\fbox{$\beta^u$}}
\efig}
$$
\end{defn}

If we think of double categories as vertices and of lax double functors as horizontal and vertical arrows, then we may think of  
horizontal oplax transformations as (vertically globular) faces between horizontal arrows, and of 
vertical lax transformations as (horizontally globular) faces between vertical arrows.  
This is the reason why we call the composition in \leref{vert comp hor.ps.tr.} ``vertical''. 
The horizontal composition of 2-cells in $\Lax_{hop}(\Aa, \Bb)$ is induced on components by the horizontal composition of the corresponding 2-cells: 
$$(\frac{\Theta}{\Theta' })(A)=
\bfig
\putmorphism(-100,250)(1,0)[F(A)`G(A)`\alpha(A)]{550}1a
 \putmorphism(430,250)(1,0)[\phantom{F(A)}`H(A) `\alpha'(A)]{580}1a

 \putmorphism(-100,-200)(1,0)[F(A)`G(A)`\beta(A)]{550}1a
 \putmorphism(450,-200)(1,0)[\phantom{F(A)}`H(A). `\beta'(A)]{580}1a

\putmorphism(-100,250)(0,-1)[\phantom{Y_2}``\alpha_0(A)]{450}1l
\putmorphism(420,250)(0,-1)[\phantom{Y_2}``]{450}1r
 \putmorphism(400,250)(0,-1)[\phantom{Y_2}``\beta_0(A)]{450}0r
\putmorphism(1020,250)(0,-1)[\phantom{Y_2}``\beta' _0(A)]{450}1r
\put(40,20){\fbox{$\Theta_A$}}
\put(700,20){\fbox{$\Theta_A' $}}
\efig
$$

The vertical composition of modifications is induced on components by the vertical composition of the corresponding 2-cells:
$$(\Sigma\cdot \Theta)(A)=
\bfig
 \putmorphism(-150,470)(1,0)[F(A)`G(A)  `\alpha(A)]{440}1a
\putmorphism(-160,470)(0,-1)[\phantom{Y_2}`F(A) `\alpha_0(A)]{400}1l
\putmorphism(-160,80)(0,-1)[\phantom{Y_2}`F(A)`\alpha_0'(A)]{400}1l
\putmorphism(300,80)(0,-1)[\phantom{Y_2}`G(A).`\gamma_0' (A)]{400}1r
\putmorphism(300,470)(0,-1)[\phantom{Y_2}`G(A) `\beta_0' (A)]{400}1r
\put(-20,280){\fbox{$\Theta_A$}}

\putmorphism(-160,70)(1,0)[\phantom{F(A)}``\beta(A)]{380}1a
\putmorphism(-100,-300)(1,0)[\phantom{Y_2}` `\gamma(A)]{300}1b
\put(-20,-130){\fbox{$\Sigma_A$}}
\efig
$$ 
It is clear that the associativity and unitality of 2-cells in both horizontal and vertical direction hold strictly.

\bigskip

Taking so to say a horizontal and a vertical restriction of modifications in \deref{modif-hv}, we obtain the definitions of:
\begin{itemize} 
\item {\em modifications between horizontal oplax transformations} given by families of (vertically globular) 2-cells 
\begin{equation} \eqlabel{m-hor}
\scalebox{0.86}{
\bfig
\putmorphism(-180,50)(1,0)[F(A)` G(A)`\alpha(A)]{550}1a
\putmorphism(-180,-270)(1,0)[F\s'(A)`G'(A) `\beta(A)]{550}1b
\putmorphism(-170,50)(0,-1)[\phantom{Y_2}``=]{320}1l
\putmorphism(350,50)(0,-1)[\phantom{Y_2}``=]{320}1r
\put(0,-140){\fbox{$\Theta_A$}}
\efig}
\end{equation} 
and axioms (m.ho.-1) and (m.ho.-2) obtained from (m.ho-vl.-1) and (m.ho-vl.-2) by ignoring the 2-cells $(\alpha_0)_f, (\beta_0)_f, 
\alpha_0^u$ and $\beta_0^u$, and 
\item {\em modifications between vertical lax transformations} given by families of (horizontally globular) 2-cells 
\begin{equation} \eqlabel{m-vert}
\scalebox{0.86}{
\bfig
\putmorphism(-180,50)(1,0)[F(A)` G(A)`=]{550}1a
\putmorphism(-180,-270)(1,0)[F\s'(A)`G'(A) `=]{550}1b
\putmorphism(-170,50)(0,-1)[\phantom{Y_2}``\alpha_0(A)]{320}1l
\putmorphism(350,50)(0,-1)[\phantom{Y_2}``\beta_0(A)]{320}1r
\put(0,-140){\fbox{$\Theta_A$}}
\efig}
\end{equation} 
and axioms (m.vl.-1) and (m.vl.-2) obtained from (m.ho-vl.-1) and (m.ho-vl.-2) by ignoring the 2-cells $\delta_{\alpha,f}, \delta_{\beta,f}, 
\alpha^u$ and $\beta^u$.  
\end{itemize}
The above two types of modifications will be used in \ssref{1-cells}.

The oplax version of the vertical transformations we will need in \coref{obtained oplax trans}. 
A {\em vertical oplax transformation} $\alpha_0$ between lax double functors $F,G: \Aa\to\Bb$ differs from its lax counterpart 
in that the (horizontally) globular 2-cell $\alpha_0^u$ for any 1v-cell $u$ in $\Aa$ goes in the other direction, namely 
$$
\scalebox{0.86}{
\bfig
 \putmorphism(-90,500)(1,0)[F(A)`F(A) `=]{540}1a
\putmorphism(-90,-300)(1,0)[G(\tilde A)`G(\tilde A) `=]{540}1a

\putmorphism(-120,500)(0,-1)[\phantom{Y_2}`F(\tilde A) `F(u)]{400}1r
\putmorphism(-120,100)(0,-1)[\phantom{Y_2}``\alpha_0(\tilde A)]{400}1r
\putmorphism(450,100)(0,-1)[\phantom{Y_2}``G(u)]{400}1l
\putmorphism(450,500)(0,-1)[\phantom{Y_2}`G(A) `\alpha_0(A)]{400}1l
\put(60,50){\fbox{$\alpha_0^u$}}
\efig}
$$
and the axioms (v.l.t.-1)-(v.l.t.-5) are accordingly changed by the analogous axioms that we will refer to as to 
(v.o.t.-1)-(v.o.t.-5). Indeed, note that only the two axioms (v.l.t.-3) and (v.l.t.-5) are changed into (v.o.t.-3) and 
(v.o.t.-5). 

Analogously to \leref{vert comp vert. lx tr.}, vertical composition of two vertical oplax transformations 
is given so that the (horizontally) globular 2-cell $(\frac{\alpha_0}{\beta_0})^u$ accordingly changes. 

\bigskip

To the modifications in $\Lax_{hop}(\Aa,\Bb)$ we will refer to as to ``modifications with respect to 
horizontally oplax and vertically lax transformations''. This is the motive for the abbreviations ``m.ho-vl.'' in their 
axioms. Apart from them we will consider:

\begin{defn} \delabel{mod-hl-vo}
Let $\alpha: F\to G, \beta: F\s'\to G'$ be horizontal lax transformations, and 
$\alpha_0: F\to F\s', \beta_0: G\to G'$ vertical oplax transformations. 
A {\em modification with respect to horizontally lax and vertically oplax transformations} $\alpha, \beta, \alpha_0, \beta_0$ 
has the (same) form $\Theta$ and is given by a collection of 2-cells $\Theta_A$ in $\Bb$ of the same form as in \equref{modification cells} 
which satisfy the rules: 

{\em (m.hl-vo.-1)} \label{m.hl-vo.-1} for every 1h-cell $f$  
$$
\scalebox{0.86}{
\bfig
\putmorphism(-150,500)(1,0)[F(A)`G(A)`\alpha(A)]{600}1a
 \putmorphism(450,500)(1,0)[\phantom{F(A)}`G(B) `G(f)]{620}1a
\putmorphism(-150,50)(1,0)[F(A)`F(B) `F(f)]{600}1a
\putmorphism(510,50)(1,0)[\phantom{Y_2}`G(B) `\alpha(B)]{580}1a

\putmorphism(-180,500)(0,-1)[\phantom{Y_2}``=]{450}1r
\putmorphism(1080,500)(0,-1)[\phantom{Y_2}``=]{450}1r
\put(350,280){\fbox{$\sigma_{\alpha,f}$}}

\putmorphism(-180,50)(0,-1)[\phantom{Y_2}``\alpha_0(A)]{450}1l
\putmorphism(1080,50)(0,-1)[\phantom{Y_3}``\beta_0(B)]{450}1r
\put(20,-180){\fbox{$(\alpha_0)_f$}}
\put(670,-180){\fbox{$\Theta_B$}}

\putmorphism(450,50)(0,-1)[\phantom{Y_2}``]{450}1r
\putmorphism(300,50)(0,-1)[\phantom{Y_2}``\alpha_0(B)]{450}0r

\putmorphism(-150,-400)(1,0)[F\s'(A)`F\s'(B) `F\s'(f)]{600}1a
\putmorphism(510,-400)(1,0)[\phantom{Y_2}`G'(B) `\beta(B)]{580}1a
\efig}
\quad=\quad
\scalebox{0.86}{
\bfig
\putmorphism(-150,500)(1,0)[F(A)`G(A)`\alpha(A)]{600}1a
 \putmorphism(450,500)(1,0)[\phantom{F(A)}`G(B) `G(f)]{620}1a

\putmorphism(-150,50)(1,0)[F\s'(A)`G'(A) `\beta(A)]{600}1a
\putmorphism(450,50)(1,0)[\phantom{Y_2}`G'(B) `G'(f)]{580}1a

\putmorphism(-180,500)(0,-1)[\phantom{Y_2}``\alpha_0(A)]{450}1l
\putmorphism(450,500)(0,-1)[\phantom{Y_2}``]{450}1r
\putmorphism(300,500)(0,-1)[\phantom{Y_2}``\beta_0(A)]{450}0r
\putmorphism(1080,500)(0,-1)[\phantom{Y_2}``\beta_0(B)]{450}1r
\put(0,280){\fbox{$\Theta_A$}} 
\put(670,280){\fbox{$(\beta_0)_f$}}

 \putmorphism(-150,-400)(1,0)[F\s'(A)`F\s'(B)`F\s'(f)]{600}1a
 \putmorphism(510,-400)(1,0)[\phantom{F(A)}`G'(B) `\beta(B)]{620}1a

\putmorphism(-180,50)(0,-1)[\phantom{Y_2}``=]{450}1l
\putmorphism(1080,50)(0,-1)[\phantom{Y_3}``=]{450}1r
\put(320,-180){\fbox{$\delta_{\beta,f}$}}

\efig}
$$
and 

\noindent 
{\em (m.hl-vo.-2)} \label{m.hl-vo.-2} for every 1v-cell $u$  
$$
\scalebox{0.86}{
\bfig
 \putmorphism(-170,500)(1,0)[F(A)`G(A) `\alpha(A)]{600}1a
 \putmorphism(450,500)(1,0)[\phantom{(B,A)}` `=]{450}1a
\putmorphism(-180,500)(0,-1)[\phantom{Y_2}`F(\tilde A) `F(u)]{450}1l
\put(650,50){\fbox{$\beta_0^u$}}
\putmorphism(-150,-400)(1,0)[F\s'(\tilde A)` `\beta(\tilde A)]{500}1a
\putmorphism(-180,50)(0,-1)[\phantom{Y_2}``\alpha_0(\tilde A)]{450}1l
\putmorphism(450,50)(0,-1)[\phantom{Y_2}`G'(\tilde A)`\beta_0(\tilde A)]{450}1r
\putmorphism(450,500)(0,-1)[\phantom{Y_2}`G(\tilde A) `G(u)]{450}1r
\put(30,280){\fbox{$\alpha^u$}}
\putmorphism(-170,50)(1,0)[\phantom{(B, \tilde A)}``\alpha(\tilde A)]{500}1a
\putmorphism(1030,50)(0,-1)[\phantom{(B, A')}` G'(\tilde A) ` G'(u)]{450}1r
\putmorphism(1030,500)(0,-1)[G(A)`G(\tilde A)` \beta_0(A)]{450}1r
\putmorphism(450,-400)(1,0)[\phantom{(B, \tilde A)}``=]{500}1b
\put(50,-150){\fbox{$\Theta_{\tilde A}$}}
\efig}
=
\scalebox{0.86}{
\bfig
 \putmorphism(-150,500)(1,0)[F(A)`F(A) `=]{600}1a
 \putmorphism(550,500)(1,0)[` `\alpha(A)]{400}1a
\putmorphism(-180,500)(0,-1)[\phantom{Y_2}`F(\tilde A) `F(u)]{450}1l
\put(30,50){\fbox{$\alpha_0^u$}}
\putmorphism(-150,-400)(1,0)[F\s'(\tilde A)` `=]{500}1a
\putmorphism(-180,50)(0,-1)[\phantom{Y_2}``\alpha_0(\tilde A)]{450}1l
\putmorphism(450,50)(0,-1)[\phantom{Y_2}`F\s'(\tilde A)` F\s'(u)]{450}1l
\putmorphism(450,500)(0,-1)[\phantom{Y_2}`F\s'(A) `\alpha_0(A)]{450}1l
\put(660,280){\fbox{$\Theta_A$}}
\putmorphism(450,50)(1,0)[\phantom{(B, \tilde A)}``\beta(A)]{500}1a
\putmorphism(1070,50)(0,-1)[\phantom{(B, A')}`G'(\tilde A).`G'(u)]{450}1r
\putmorphism(1070,500)(0,-1)[G(A)`G'(A)`\beta_0(A)]{450}1r
\putmorphism(450,-400)(1,0)[\phantom{(B, \tilde A)}``\beta(\tilde A)]{500}1a
\put(640,-170){\fbox{$\beta^u$ } } 
\efig}
$$
\end{defn}

Observe that analogously to the double category $\Lax_{hop}(\Aa,\Bb)$,  
lax double functors as 0-cells, horizontal {\em lax} transformations as 1h-cells, vertical {\em oplax} transformations as 1v-cells and 
modifications from \deref{mod-hl-vo} form another double category $\Lax_{hop}^*(\Aa,\Bb)$.

\section{Candidates for inner-hom and lax Gray type monoidal product in $Dbl_{lx}$, and lax double quasi-functors} \selabel{mon str}

Let $Dbl_{lx}$ denote the category of double categories and lax double functors. 
For the purpose of exploring the existence and properties of an inner-hom in $Dbl_{lx}$, we will denote 
$\llbracket\Aa,\Bb\rrbracket:=\Lax_{hop}(\Aa,\Bb)$ for two double categories $\Aa, \Bb$. 

Observe that the double category $\llbracket\Aa,\Bb\rrbracket$ is analogous to one in \cite[Section 2.2]{Gabi}. 
The double category $\llbracket\Aa,\Bb\rrbracket$ from \cite[Section 2.2]{Gabi} consists of: 0-cells are strict double functors, 
1h-cells are horizontal pseudo-transformations, 1v-cells are vertical pseudo-transformations and 2-cells are modifications among the latter two. 
It is the strictness of double functors that allows that $\llbracket-,-\rrbracket$ be a bifunctor $(Dbl_{st})^{op}\times Dbl_{st}\to Dbl_{st}$, 
where $Dbl_{st}$ is the category of double categories and strict double functors. Furthermore, the author constructs a Gray type monoidal product 
$-\ot-:Dbl_{st}\times Dbl_{st}\to Dbl_{st}$, so that there is an adjunction $(-\ot\Dd, \llbracket\Dd,-\rrbracket)$ for every double category $\Dd$ 
and $(Dbl_{st},\ot)$ is a closed monoidal category. 

In contrast to the case where the double functors are strict or pseudo, our  
double category $\llbracket\Aa,\Bb\rrbracket$ will not induce a bifunctor $\llbracket-,-\rrbracket: (Dbl_{lx})^{op}\times Dbl_{lx}\to Dbl_{lx}$. 
We will explain this in \ssref{not-inner}. 
As can be appreciated from the previous Section, all the cells in our double category $\llbracket\Aa,\Bb\rrbracket$ 
are more general than in \cite[Section 2.2]{Gabi}. Though, the price we pay is that we loose closedness for $Dbl_{lx}$. 

After \ssref{not-inner} we will characterize a lax double functor $F:\Aa\to\llbracket\Bb,\Cc\rrbracket$ for another double category $\Cc$ in terms 
of a bifunctor from the Cartesian product $\Aa\times\Bb\to \Cc$ of double categories. 
Setting $\Cc=\Aa\times\Bb$ and reading off 
the structure of the image double category $F(\Aa)(\Bb)$, we will obtain a full description of a new structure on the underlying 
 double category $\Aa\times\Bb$. Thus obtained structure we will denote by $\Aa\ot\Bb$. 
In \ssref{not Gray} we will comment why this does not provide a Gray type monoidal product on the category $Dbl_{lx}$, although  
the analogous construction works well for the case when pseudo (or strict) double functors are used instead of lax double functors. 
We will though prove a universal property that it satisfies in \ssref{Gray}.

\subsection{Why $\llbracket\Aa,\Bb\rrbracket$ is not an inner-hom} \sslabel{not-inner}

The double category $\llbracket \Aa,\Bb\rrbracket$ does not induce a functor 
$\llbracket -,-\rrbracket: (Dbl_{lx})^{op}\times Dbl_{lx}\to Dbl_{lx}$, 
and so it can not play the role of the inner hom in $Dbl_{lx}$. 
The thing is that at the level of morphisms, given lax double functors $F:\Aa\to\Aa'$ and $G:\Bb\to\Bb'$, what should be a lax double 
functor $\llbracket F,G\rrbracket: \llbracket \Aa',\Bb\rrbracket \to \llbracket \Aa,\Bb'\rrbracket$ can not be defined on 1h-cells. 
Namely, the components of their images should be naturally defined in an analogous way as it has been done in similar constructions, 
that is, as $Gx^{Fh}$ in \cite[Section 2.3]{Gabi} (strict double functor case), or as $\delta_{H(\alpha),f}$ in \cite[Lemma 3.4]{Fem} 
({\em double} pseudofunctor case). Namely, $H:=\llbracket F,G\rrbracket$ on a 1h-cell, that is a horizontal oplax transformation $\alpha: F\to G$, 
should give another horizontal oplax transformation $H(\alpha)$. Its 2-cell component at a 1h-cell $f$ 
in \cite[Lemma 3.4]{Fem} 
was defined via
$$
\bfig
 \putmorphism(-150,250)(1,0)[HF(A)`HF(B)`HF(f)]{600}1a
 \putmorphism(450,250)(1,0)[\phantom{F(A)}`HG(B) `H(\alpha(B))]{640}1a
\putmorphism(-150,-200)(1,0)[HF(A)`HG(A) `H(\alpha(A))]{640}1a
 \putmorphism(450,-200)(1,0)[\phantom{A'\ot B'}` HG(B) `HG(f)]{680}1a

\putmorphism(-180,250)(0,-1)[\phantom{Y_2}``=]{450}1l
\putmorphism(1120,250)(0,-1)[\phantom{Y_3}``=]{450}1r
\put(320,40){\fbox{$\delta_{H(\alpha),f}$}}

\efig
\quad=\quad
\bfig
 \putmorphism(-150,650)(1,0)[HF(A)`HF(B)`HF(f)]{600}1a
 \putmorphism(450,650)(1,0)[\phantom{F(A)}`HG(B) `H(\alpha(B))]{640}1a

\putmorphism(-100,650)(0,-1)[\phantom{Y_2}``=]{450}1l
\putmorphism(1050,650)(0,-1)[\phantom{Y_2}``=]{450}1r
\put(280,460){\fbox{$H_{\alpha(B)F(f)}$}}

\putmorphism(-580,200)(0,-1)[\phantom{Y_2}``=]{450}1l
\putmorphism(1520,200)(0,-1)[\phantom{Y_2}``=]{450}1r
\putmorphism(-550,200)(1,0)[HF(A)`\phantom{HF(A)}`=]{550}1a
\putmorphism(-550,-250)(1,0)[HF(A)`\phantom{HF(A)}`=]{550}1a
\put(-500,-40){\fbox{$H^{F(A)}$}}
\put(1060,-40){\fbox{$(H^{G(B)})^{-1}$}}

 \putmorphism(980,200)(1,0)[\phantom{HF(A)}`HG(B) `=]{550}1a
 \putmorphism(1000,-250)(1,0)[\phantom{HF(A)}`HG(B) `=]{550}1a

\putmorphism(-20,200)(1,0)[HF(A)`HG(B)`H(\alpha(B)F(f))]{1000}1a
 \putmorphism(-20,-250)(1,0)[HF(A)`HG(B)`H\big(G(f)\alpha(A)\big)]{1000}1a

\putmorphism(-80,200)(0,-1)[\phantom{Y_2}``]{450}1r
\putmorphism(-100,200)(0,-1)[\phantom{Y_2}``H(1^A)]{450}0r

\putmorphism(1020,200)(0,-1)[\phantom{Y_2}``]{450}1l
\putmorphism(1050,200)(0,-1)[\phantom{Y_2}``H(1^A)]{450}0l

\put(300,10){\fbox{$H(\delta_{\alpha,f})$}}
\put(300,-480){\fbox{$H_{G(f)\alpha(A)}^{-1}$}}
\putmorphism(-150,-700)(1,0)[HF(A)`HG(A) `H(\alpha(A))]{600}1a
 \putmorphism(450,-700)(1,0)[\phantom{HF(A)}` HG(B), `HG(f)]{600}1a

\putmorphism(-80,-250)(0,-1)[\phantom{Y_2}``=]{450}1l
\putmorphism(1020,-250)(0,-1)[\phantom{Y_3}``=]{450}1r

\efig
$$
which made sense there, as $H$ was a double pseudofunctor. The definition in \cite[Section 2.3]{Gabi} is 
a strict double functor version of $H$. So, in our present context a colax double functor structure is missing. 
In conclusion, 0-cells of $\llbracket \Aa,\Bb\rrbracket$ should be double functors which are strict or pseudo in the horizontal direction 
in order to be able to define $\llbracket F,G\rrbracket(\alpha)$.

\subsection{Generating lax double quasi-functors}

Having in mind the definition of a lax double functor and of $\llbracket\Bb,\Cc\rrbracket$, when writing out the list of 
the data and relations 
that determine a lax double functor $\F:\Aa\to\llbracket\Bb,\Cc\rrbracket$, one gets the following characterization of it:

\begin{prop} \prlabel{char df}
A lax double functor $\F:\Aa\to\llbracket\Bb,\Cc\rrbracket$ of double categories consists of the following: \\
1. lax double functors 
$$(-,A):\Bb\to\Cc\quad\text{ and}\quad (B,-):\Aa\to\Cc$$ 
such that $(-,A)\vert_B=(B,-)\vert_A=(B,A)$, 
for objects $A\in\Aa, B\in\Bb$, \\
2. 2-cells
$$
\scalebox{0.86}{
\bfig
 \putmorphism(-150,50)(1,0)[(B,A)`(B', A)`(k, A)]{600}1a
 \putmorphism(450,50)(1,0)[\phantom{A\ot B}`(B', A') `(B', K)]{680}1a
\putmorphism(-180,50)(0,-1)[\phantom{Y_2}``=]{450}1r
\putmorphism(1100,50)(0,-1)[\phantom{Y_2}``=]{450}1r
\put(350,-190){\fbox{$(k,K)$}}
 \putmorphism(-150,-400)(1,0)[(B,A)`(B,A')`(B,K)]{600}1a
 \putmorphism(450,-400)(1,0)[\phantom{A\ot B}`(B', A') `(k, A')]{680}1a
\efig}
$$

$$
\scalebox{0.86}{
\bfig
\putmorphism(-150,50)(1,0)[(B,A)`(B,A')`(B,K)]{600}1a
\putmorphism(-150,-400)(1,0)[(\tilde B, A)`(\tilde B,A') `(\tilde B,K)]{640}1a
\putmorphism(-180,50)(0,-1)[\phantom{Y_2}``(u,A)]{450}1l
\putmorphism(450,50)(0,-1)[\phantom{Y_2}``(u,A')]{450}1r
\put(0,-180){\fbox{$(u, K)$}}
\efig}
\quad
\scalebox{0.86}{
\bfig
\putmorphism(-150,50)(1,0)[(B,A)`(B',A)`(k,A)]{600}1a
\putmorphism(-150,-400)(1,0)[(B, \tilde A)`(B', \tilde A) `(k,\tilde A)]{640}1a
\putmorphism(-180,50)(0,-1)[\phantom{Y_2}``(B,U)]{450}1l
\putmorphism(450,50)(0,-1)[\phantom{Y_2}``(B',U)]{450}1r
\put(0,-180){\fbox{$(k,U)$}}
\efig}
$$

$$
\scalebox{0.86}{
\bfig
 \putmorphism(-150,500)(1,0)[(B,A)`(B,A) `=]{600}1a
\putmorphism(-180,500)(0,-1)[\phantom{Y_2}`(B, \tilde A) `(B,U)]{450}1l
\put(0,50){\fbox{$(u,U)$}}
\putmorphism(-150,-400)(1,0)[(\tilde B, \tilde A)`(\tilde B, \tilde A) `=]{640}1a
\putmorphism(-180,50)(0,-1)[\phantom{Y_2}``(u,\tilde A)]{450}1l
\putmorphism(450,50)(0,-1)[\phantom{Y_2}``(\tilde B, U)]{450}1r
\putmorphism(450,500)(0,-1)[\phantom{Y_2}`(\tilde B, A) `(u,A)]{450}1r
\efig}
$$ 
in $\Cc$ for every 1h-cells $A\stackrel{F}{\to} A'$ and $B\stackrel{f}{\to} B'$ and 1v-cells $A\stackrel{U}{\to} \tilde A$ and 
$B\stackrel{u}{\to} \tilde B$
which satisfy: 

\noindent $\bullet$ \quad ($(1_B,K)$) \label{1B,K}  
$$
\scalebox{0.86}{
\bfig
 \putmorphism(-170,420)(1,0)[(B,A)`(B,A)`=]{500}1a
\putmorphism(-170,50)(1,0)[(B,A)`(B,A) `(1_B,A)]{520}1a
\putmorphism(360,50)(1,0)[\phantom{F(A)}`(B,A') `(B,K)]{500}1a
\putmorphism(350,400)(1,0)[\phantom{F(A)}` (B,A') `(B,K)]{500}1a

\putmorphism(-170,420)(0,-1)[\phantom{Y_2}``=]{350}1l
\putmorphism(320,420)(0,-1)[\phantom{Y_2}``=]{370}1r
\put(-80,250){\fbox{$(-,A)_B$}}

\putmorphism(830,420)(0,-1)[\phantom{Y_2}``=]{350}1r
\put(470,240){\fbox{$\Id_{(B,K)}$}}
\putmorphism(-180,50)(0,-1)[\phantom{Y_2}``=]{350}1r
\putmorphism(830,50)(0,-1)[\phantom{Y_2}``=]{350}1r
 \putmorphism(-150,-300)(1,0)[(B,A)`(B,A')`(B,K)]{500}1a
 \putmorphism(350,-300)(1,0)[\phantom{A\ot B}`(B, A') `(1_B, A')]{560}1a
\put(170,-130){\fbox{$(1_B,K)$}}
\efig}
\scalebox{0.86}{
\quad=\quad
\bfig
 \putmorphism(-210,420)(1,0)[(B,A)`\phantom{F(A)} `(B,K)]{500}1a
\putmorphism(330,420)(1,0)[(B,A')`(B,A')`=]{500}1a
\putmorphism(-210,420)(0,-1)[\phantom{Y_2}``=]{370}1l
\putmorphism(320,420)(0,-1)[\phantom{Y_2}``=]{370}1l
\putmorphism(910,420)(0,-1)[\phantom{Y_2}``=]{370}1r
 \putmorphism(-210,50)(1,0)[(B,A)`\phantom{Y_2}`(B,K)]{470}1a
 \put(420,250){\fbox{$(-,A')_{B}$}} 
\putmorphism(360,50)(1,0)[(B,A')`(B,A') `(1_B,A')]{580}1a
\put(-130,240){\fbox{$\Id_{(B,K)}$}}
\efig}
$$

\noindent $\bullet$ \quad ($(k,1_A)$) \label{k,1A}
$$
\scalebox{0.86}{
\bfig
 \putmorphism(-150,420)(1,0)[(B,A)`(B,A)`=]{500}1a
\putmorphism(360,420)(1,0)[\phantom{F(A)}`(B,A') `(k,A)]{500}1a
\putmorphism(830,420)(0,-1)[\phantom{Y_2}``=]{350}1r
\put(460,250){\fbox{$\Id_{(k,A)}$}}

\putmorphism(-180,420)(0,-1)[\phantom{Y_2}``=]{370}1l
\putmorphism(320,420)(0,-1)[\phantom{Y_2}``=]{370}1r
 \putmorphism(-150,50)(1,0)[(B,A)`(B,A)`(B,1_A)]{500}1a
 \put(-80,250){\fbox{$(B,-)_A$}} 
\putmorphism(350,50)(1,0)[\phantom{F(A)}`(B',A) `(k,A)]{560}1a
\efig}
\quad
=
\quad
\scalebox{0.86}{
\bfig
\putmorphism(-150,420)(1,0)[(B,A)` \phantom{Y_2}`(k,A)]{450}1a
\putmorphism(370,420)(1,0)[(B',A)` (B',A) `=]{500}1a
\putmorphism(-170,420)(0,-1)[\phantom{Y_2}``=]{350}1l
\putmorphism(350,420)(0,-1)[\phantom{Y_2}``=]{350}1l
\putmorphism(860,420)(0,-1)[\phantom{Y_2}``=]{350}1r
 \putmorphism(420,50)(1,0)[\phantom{Y_2}`(B',A)`(B',1_A)]{520}1a
\putmorphism(-150,50)(1,0)[(B,A)` (B',A) `(k,A)]{500}1a
 \put(430,250){\fbox{$(B',-)_A$}} 
\put(-120,250){\fbox{$\Id_{(k,A)}$}}

\putmorphism(-180,50)(0,-1)[\phantom{Y_2}``=]{350}1r
\putmorphism(860,50)(0,-1)[\phantom{Y_2}``=]{350}1r
\putmorphism(-150,-300)(1,0)[(B,A)`(B,A)`(B,1_A)]{500}1a
 \putmorphism(350,-300)(1,0)[\phantom{A\ot B}`(B', A) `(k, A)]{560}1a
\put(200,-110){\fbox{$(k,1_A)$}}
\efig}
$$
where the 2-cells $(-,A)_{B}$ and $(B,-)_A$ come from laxity of the lax double functors $(-,A)$ and $(B,-)$ 

\noindent $\bullet$ \quad ($(u,1_A)$) \label{u,1A} 
$$
\scalebox{0.86}{
\bfig
\putmorphism(-250,500)(1,0)[(B,A)`(B,A)` =]{550}1a
 \putmorphism(-250,50)(1,0)[(\tilde B,A)`(\tilde B,A)` =]{550}1a
 \putmorphism(-250,-400)(1,0)[(\tilde B,A)`(\tilde B,A)` (\tilde B,1_A)]{550}1a

\putmorphism(-280,500)(0,-1)[\phantom{Y_2}``(u,A)]{450}1l
 \putmorphism(-280,70)(0,-1)[\phantom{F(A)}` `=]{450}1l

\putmorphism(300,500)(0,-1)[\phantom{Y_2}``(u,A)]{450}1r
\putmorphism(300,70)(0,-1)[\phantom{Y_2}``=]{450}1r
\put(-120,270){\fbox{$Id^{(u,A)}$}}
\put(-120,-150){\fbox{$(\tilde B,-)_A$}}
\efig}
=
\scalebox{0.86}{
\bfig
\putmorphism(-250,500)(1,0)[(B,A)`(B,A)` =]{550}1a
 \putmorphism(-250,50)(1,0)[(B,A)`(B,A)` (B,1_A)]{550}1a
 \putmorphism(-250,-400)(1,0)[(\tilde B,A)`(\tilde B,A)` (\tilde B,1_A)]{550}1a

\putmorphism(-280,500)(0,-1)[\phantom{Y_2}``= ]{450}1l
 \putmorphism(-280,70)(0,-1)[\phantom{F(A)}` `(u,A)]{450}1l

\putmorphism(300,500)(0,-1)[\phantom{Y_2}``=]{450}1r
\putmorphism(300,70)(0,-1)[\phantom{Y_2}``(u,A)]{450}1r
\put(-120,290){\fbox{$(B,-)_A$}}
\put(-140,-150){\fbox{$(u,1_A)$}}
\efig}
$$

\noindent $\bullet$ \quad ($(1_B,U)$) \label{1BU} 
$$\scalebox{0.86}{
\bfig
\putmorphism(-250,500)(1,0)[(B,A)`(B,A)` =]{550}1a
 \putmorphism(-250,50)(1,0)[(B,A)`(B,A)` (1_B,A)]{550}1a
 \putmorphism(-250,-400)(1,0)[(B,\tilde A)`(B,\tilde A)` (1_B,\tilde A)]{550}1a

\putmorphism(-280,500)(0,-1)[\phantom{Y_2}``= ]{450}1l
 \putmorphism(-280,70)(0,-1)[\phantom{F(A)}` `(B,U)]{450}1l

\putmorphism(300,500)(0,-1)[\phantom{Y_2}``=]{450}1r
\putmorphism(300,70)(0,-1)[\phantom{Y_2}``(B,U)]{450}1r
\put(-150,290){\fbox{$(-,A)_B$}}
\put(-150,-160){\fbox{$(1_B,U)$}}
\efig}\quad
=
\scalebox{0.86}{
\bfig
\putmorphism(-250,500)(1,0)[(B,A)`(B,A)` =]{550}1a
 \putmorphism(-250,50)(1,0)[(B,\tilde A)`(B,\tilde A)` =]{550}1a
 \putmorphism(-250,-400)(1,0)[(B,\tilde A)`(B,\tilde A)` (\tilde B,1_{\tilde A})]{550}1a

\putmorphism(-280,500)(0,-1)[\phantom{Y_2}``(B,U)]{450}1l
 \putmorphism(-280,70)(0,-1)[\phantom{F(A)}` `=]{450}1l

\putmorphism(300,500)(0,-1)[\phantom{Y_2}``(B,U)]{450}1r
\putmorphism(300,70)(0,-1)[\phantom{Y_2}``=]{450}1r
\put(-140,270){\fbox{$Id^{(B,U)}$}}
\put(-140,-160){\fbox{$(-, \tilde A)_B$}}
\efig}
$$

\noindent $\bullet$ \qquad ($(1^B,K)$) \label{1uB,K} \quad $(1^B,K)=Id_{(B,K)}$ \qquad\text{and}\qquad 
$\bullet$ \quad ($(k,1^A)$) \label{k,1uA} \qquad $(k,1^A)=Id_{(k,A)}$ 

\noindent $\bullet$ \qquad ($(1^B,U)$) \label{1uB,U} \quad $(1^B,U)=Id^{(B,U)}$ \qquad\text{and}\qquad \hspace{-0,22cm} 
$\bullet$ \quad ($(u,1^A)$) \label{u,1uA} \qquad $(u,1^A)=Id^{(u,A)}$

\noindent $\bullet$ \quad ($(k'k,K)$) \label{k'k,K} 
$$
\scalebox{0.86}{
\bfig
 \putmorphism(450,500)(1,0)[(B', A) `(B'', A) `(k', A)]{680}1a
 \putmorphism(1140,500)(1,0)[\phantom{A\ot B}`(B'', A') ` (B'', K)]{680}1a

 \putmorphism(-150,50)(1,0)[(B, A) `(B', A)`(k, A)]{600}1a
 \putmorphism(450,50)(1,0)[\phantom{A\ot B}`(B', A') `(B', K)]{680}1a
 \putmorphism(1130,50)(1,0)[\phantom{A\ot B}`(B'', A') ` (k', A')]{680}1a

\putmorphism(450,500)(0,-1)[\phantom{Y_2}``=]{450}1r
\putmorphism(1750,500)(0,-1)[\phantom{Y_2}``=]{450}1r
\put(1000,270){\fbox{$ (k',K)$}}

 \putmorphism(-150,-400)(1,0)[(B, A)`(B, A') `(B,K)]{640}1a
 \putmorphism(480,-400)(1,0)[\phantom{A'\ot B'}`(B', A') `(k, A')]{680}1a

\putmorphism(-180,50)(0,-1)[\phantom{Y_2}``=]{450}1l
\putmorphism(1120,50)(0,-1)[\phantom{Y_3}``=]{450}1r
\put(310,-200){\fbox{$ (k,K)$}}

 \putmorphism(1170,-400)(1,0)[\phantom{A\ot B}`(B'', A') ` (k', A')]{650}1a
\putmorphism(450,-400)(0,-1)[\phantom{Y_2}``=]{450}1r
\putmorphism(1750,-400)(0,-1)[\phantom{Y_2}``=]{450}1r

 \putmorphism(480,-850)(1,0)[(B, A') `(B'', A'') `(k'k, A')]{1320}1a
\put(920,-620){\fbox{$ (-,A')_{k'k}$}}
\efig}
\quad=\quad
\scalebox{0.86}{
\bfig
 \putmorphism(-150,500)(1,0)[(B, A) `(B', A)`(k,A)]{600}1a
 \putmorphism(450,500)(1,0)[(B', A) `(B'', A) `(k',A)]{680}1a
\putmorphism(-180,500)(0,-1)[\phantom{Y_2}``=]{450}1r
\put(310,270){\fbox{$ (-,A)_{k'k}$}}

 \putmorphism(-150,50)(1,0)[(B,A)` `(k'k,A)]{1130}1a
 \putmorphism(1130,50)(1,0)[(B'',A)`(B'', A') ` (B'',K)]{680}1a

\putmorphism(1050,500)(0,-1)[\phantom{Y_2}``=]{450}1r

 \putmorphism(-150,-400)(1,0)[(B, A)`(B, A') `(B,K)]{640}1a
 \putmorphism(530,-400)(1,0)[\phantom{Y_2X}`(B'', A') `(k'k,A')]{1270}1a
\put(620,-200){\fbox{$ (k'k,K)$}}

\putmorphism(-180,50)(0,-1)[\phantom{Y_2}``=]{450}1l
\putmorphism(1750,50)(0,-1)[\phantom{Y_3}``=]{450}1r
\efig}
$$ 
where $(-,A)_{k'k}$ is the 2-cell from the laxity of $(-,A)$ 

\noindent $\bullet$ \quad ($(k,K'K)$) \label{k,K'K} 
$$
\scalebox{0.86}{
\bfig
 \putmorphism(450,500)(1,0)[(B', A) `(B', A') `(B', K)]{680}1a
 \putmorphism(1140,500)(1,0)[\phantom{A\ot B}`(B', A'') ` (B', K')]{680}1a

 \putmorphism(-150,50)(1,0)[(B, A) `(B', A)`(k,A)]{600}1a
 \putmorphism(450,50)(1,0)[\phantom{A\ot B}`(B', A'') `(B', K'K)]{1350}1a

\putmorphism(450,500)(0,-1)[\phantom{Y_2}``=]{450}1r
\putmorphism(1750,500)(0,-1)[\phantom{Y_2}``=]{450}1r
\put(880,290){\fbox{$ (B',-)_{K'K}$}}

 \putmorphism(-150,-400)(1,0)[(B, A)`(B, A'') `(B, K'K)]{980}1a
 \putmorphism(780,-400)(1,0)[\phantom{A'\ot B'}`(B', A'') `(k,A'')]{980}1a

\putmorphism(-180,50)(0,-1)[\phantom{Y_2}``=]{450}1l
\putmorphism(1750,50)(0,-1)[\phantom{Y_3}``=]{450}1r
\put(560,-200){\fbox{$ (k,K'K)$}}

\efig}
=
\scalebox{0.86}{
\bfig
 \putmorphism(-150,450)(1,0)[(B,A)`(B',A)`(k,A)]{600}1a
 \putmorphism(450,450)(1,0)[\phantom{A\ot B}`(B', A') `(B',K)]{680}1a

 \putmorphism(-150,0)(1,0)[(B,A)`(B,A')`(B,K)]{600}1a
 \putmorphism(450,0)(1,0)[\phantom{A\ot B}`(B', A') `(k,A')]{680}1a
 \putmorphism(1100,0)(1,0)[\phantom{A'\ot B'}`(B', A'') `(B', K')]{660}1a

\putmorphism(-180,450)(0,-1)[\phantom{Y_2}``=]{450}1r
\putmorphism(1100,450)(0,-1)[\phantom{Y_2}``=]{450}1r
\put(350,210){\fbox{$(k,K)$}}
\put(1000,-250){\fbox{$(k,K')$}}

 \putmorphism(450,-450)(1,0)[\phantom{A''\ot B'}` (B, A'') `(B, K')]{680}1a
 \putmorphism(1100,-450)(1,0)[\phantom{A''\ot B'}`(B', A'') ` (k, A'')]{660}1a

\putmorphism(450,0)(0,-1)[\phantom{Y_2}``=]{450}1l
\putmorphism(1750,0)(0,-1)[\phantom{Y_2}``=]{450}1r
 \putmorphism(-150,-450)(1,0)[(B,A)`(B,A')`(B,K)]{600}1a
\putmorphism(-180,-450)(0,-1)[\phantom{Y_2}``=]{450}1r
\putmorphism(1100,-450)(0,-1)[\phantom{Y_2}``=]{450}1r
 \putmorphism(-150,-900)(1,0)[(B,A)`(B,A'')`(B, K'K)]{1280}1a
\put(260,-670){\fbox{$(B,-)_{K'K}$}}
\efig}
$$
where $(B,-)_{K'K}$ is the 2-cell from the laxity of $(B,-)$ 

\noindent $\bullet$ \quad ($(u, K'K)$) \label{u,K'K} 
$$
\scalebox{0.86}{
\bfig
\putmorphism(-150,500)(1,0)[(B,A)`(B,A')`(B,K)]{600}1a
 \putmorphism(470,500)(1,0)[\phantom{F(A)}`(B,A'') `(B, K')]{600}1a
 \putmorphism(-150,50)(1,0)[(B,A)`(B,A'')`(B,K'K)]{1220}1a

\putmorphism(-180,500)(0,-1)[\phantom{Y_2}``=]{450}1r
\putmorphism(1080,500)(0,-1)[\phantom{Y_2}``=]{450}1r
\put(300,290){\fbox{$(B,-)_{K'K}$}}

\putmorphism(-150,-400)(1,0)[(\tilde B,A)`(\tilde B,A'') `(\tilde B,K'K)]{1200}1a

\putmorphism(-180,50)(0,-1)[\phantom{Y_2}``(u,A)]{450}1l
\putmorphism(1080,50)(0,-1)[\phantom{Y_3}``(u,A'')]{450}1r
\put(300,-180){\fbox{$(u,K'K)$}} 
\efig}
=
\scalebox{0.86}{
\bfig
\putmorphism(-150,500)(1,0)[(B,A)`(B,A')`(B,K)]{600}1a
 \putmorphism(470,500)(1,0)[\phantom{F(A)}`(B,A'') `(B, K')]{600}1a

 \putmorphism(-150,50)(1,0)[(\tilde B, A)`(\tilde B,A')`(\tilde B,K)]{600}1a
 \putmorphism(470,50)(1,0)[\phantom{F(A)}`(B'',\tilde A) `(\tilde B,K')]{620}1a

\putmorphism(-180,500)(0,-1)[\phantom{Y_2}``(u,A)]{450}1l
\putmorphism(450,500)(0,-1)[\phantom{Y_2}``]{450}1r
\putmorphism(300,500)(0,-1)[\phantom{Y_2}``(u,A')]{450}0r
\putmorphism(1080,500)(0,-1)[\phantom{Y_2}``(u,A'')]{450}1r
\put(-40,280){\fbox{$(u,K)$}}
\put(620,280){\fbox{$(u,K')$}}

\putmorphism(-150,-400)(1,0)[(\tilde B, A)`(\tilde B,A'') `(\tilde B,K'K)]{1200}1a

\putmorphism(-180,50)(0,-1)[\phantom{Y_2}``=]{450}1l
\putmorphism(1080,50)(0,-1)[\phantom{Y_3}``=]{450}1r
\put(260,-160){\fbox{$(\tilde B,-)_{K'K}$}}

\efig}
$$

\noindent $\bullet$ \quad ($(k'k, U)$) \label{k'k,U} 
$$
\scalebox{0.86}{
\bfig
\putmorphism(-150,500)(1,0)[(B,A)`(B',A)`(k,A)]{600}1a
 \putmorphism(450,500)(1,0)[\phantom{F(A)}`(B'',A) `(k',A)]{620}1a

 \putmorphism(-150,50)(1,0)[(B,\tilde A)`(B',\tilde A)`(k,\tilde A)]{600}1a
 \putmorphism(470,50)(1,0)[\phantom{F(A)}`(B'',\tilde A) `(k',\tilde A)]{620}1a

\putmorphism(-180,500)(0,-1)[\phantom{Y_2}``(B,U)]{450}1l
\putmorphism(450,500)(0,-1)[\phantom{Y_2}``]{450}1r
\putmorphism(300,500)(0,-1)[\phantom{Y_2}``(B',U)]{450}0r
\putmorphism(1080,500)(0,-1)[\phantom{Y_2}``(B'',U)]{450}1r
\put(-40,280){\fbox{$(k,U)$}}
\put(620,280){\fbox{$(k',U)$}}

\putmorphism(-150,-400)(1,0)[(B,\tilde A)`(B'',\tilde A) `(k'k,\tilde A)]{1200}1a

\putmorphism(-180,50)(0,-1)[\phantom{Y_2}``=]{450}1l
\putmorphism(1080,50)(0,-1)[\phantom{Y_3}``=]{450}1r
\put(260,-160){\fbox{$(-,\tilde A)_{k'k}$}}
\efig}
=
\scalebox{0.86}{
\bfig
\putmorphism(-150,500)(1,0)[(B,A)`(B',A)`(k,A)]{600}1a
 \putmorphism(450,500)(1,0)[\phantom{F(A)}`(B'',A) `(k',A)]{620}1a
 \putmorphism(-150,50)(1,0)[(B,A)`(B'',A)`(k'k, A)]{1220}1a

\putmorphism(-180,500)(0,-1)[\phantom{Y_2}``=]{450}1r
\putmorphism(1080,500)(0,-1)[\phantom{Y_2}``=]{450}1r
\put(300,290){\fbox{$(-,A)_{k'k}$}}

\putmorphism(-150,-400)(1,0)[(B,\tilde A)`(B'',\tilde A) `(k'k,\tilde A)]{1200}1a

\putmorphism(-180,50)(0,-1)[\phantom{Y_2}``(B,U)]{450}1l
\putmorphism(1080,50)(0,-1)[\phantom{Y_3}``(B'',U)]{450}1r
\put(300,-180){\fbox{$(k'k,U)$}} 
\efig}
$$

\noindent $\bullet$ \quad ($(\frac{u}{u'}, K)$) \label{uu',K} \qquad $(\frac{u}{u'}, K)=\frac{(u,K)}{(u', K)}$  \qquad\text{and}\qquad 
$\bullet$ \quad ($(k,\frac{U}{U'})$) \label{k,UU'}  \qquad $(k,\frac{U}{U'})=\frac{(k,U)}{(k,U')}$ 

\noindent $\bullet$ \quad ($(u,\frac{U}{U'})$) \label{u,UU'} 
$$(u,\frac{U}{U'})=
\scalebox{0.86}{
\bfig
 \putmorphism(-150,500)(1,0)[(B,A)`(B,A) `=]{600}1a
\putmorphism(-180,500)(0,-1)[\phantom{Y_2}`(B, \tilde A) `(B,U)]{450}1l
\put(0,50){\fbox{$(u,U)$}}
\putmorphism(-150,-400)(1,0)[(\tilde B, \tilde A)`(\tilde B, \tilde A) `=]{640}1a
\putmorphism(-180,50)(0,-1)[\phantom{Y_2}``(u,\tilde A)]{450}1l
\putmorphism(450,50)(0,-1)[\phantom{Y_2}``(\tilde B, U)]{450}1r
\putmorphism(450,500)(0,-1)[\phantom{Y_2}`(\tilde B, A) `(u,A)]{450}1r
\putmorphism(-820,50)(1,0)[(B, \tilde A)``=]{520}1a
\putmorphism(-820,50)(0,-1)[\phantom{(B, \tilde A')}``(B,U')]{450}1l
\putmorphism(-820,-400)(0,-1)[(B, \tilde A')`(\tilde B, \tilde A')`(u,\tilde A')]{450}1l
\putmorphism(-820,-850)(1,0)[\phantom{(B, \tilde A)}``=]{520}1a
\putmorphism(-150,-400)(0,-1)[(\tilde B, \tilde A)`(\tilde B, \tilde A') `(\tilde B, U')]{450}1r
\put(-650,-630){\fbox{$(u,U')$}}
\efig}
$$

\noindent $\bullet$ \quad ($(\frac{u}{u'},U)$) \label{uu',U} 
$$(\frac{u}{u'},U)=
\scalebox{0.86}{
\bfig
 \putmorphism(-150,500)(1,0)[(B,A)`(B,A) `=]{600}1a
\putmorphism(-180,500)(0,-1)[\phantom{Y_2}`(B, \tilde A) `(B,U)]{450}1l
\put(0,50){\fbox{$(u,U)$}}
\putmorphism(-150,-400)(1,0)[(\tilde B, \tilde A)` `=]{500}1a
\putmorphism(-180,50)(0,-1)[\phantom{Y_2}``(u,\tilde A)]{450}1l
\putmorphism(450,50)(0,-1)[\phantom{Y_2}`(\tilde B, \tilde A)`(\tilde B, U)]{450}1r
\putmorphism(450,500)(0,-1)[\phantom{Y_2}`(\tilde B, A) `(u,A)]{450}1r
\putmorphism(450,50)(1,0)[\phantom{(B, \tilde A)}`(\tilde B, A)`=]{620}1a
\putmorphism(1070,50)(0,-1)[\phantom{(B, \tilde A')}``(u',A)]{450}1r
\putmorphism(1070,-400)(0,-1)[(\tilde B', A)`(\tilde B', \tilde A)`(\tilde B', U)]{450}1r
\putmorphism(450,-850)(1,0)[\phantom{(B, \tilde A)}``=]{520}1a
\putmorphism(450,-400)(0,-1)[\phantom{(B, \tilde A)}`(\tilde B', \tilde A) `(U',\tilde A)]{450}1l
\put(600,-630){\fbox{$(u',U)$}}
\efig}
$$

\noindent $\bullet$ \quad ($(k,K)$-l-nat) \label{k,K-l-nat} 
$$
\scalebox{0.86}{
\bfig
 \putmorphism(-150,500)(1,0)[(B,A)`(B', A)`(k, A)]{600}1a
 \putmorphism(450,500)(1,0)[\phantom{A\ot B}`(B', A') `(B', K)]{680}1a
 \putmorphism(-150,50)(1,0)[(B,A)`(B,A')`(B,K)]{600}1a
 \putmorphism(450,50)(1,0)[\phantom{A\ot B}`(B', A') `(k, A')]{680}1a

\putmorphism(-180,500)(0,-1)[\phantom{Y_2}``=]{450}1r
\putmorphism(1100,500)(0,-1)[\phantom{Y_2}``=]{450}1r
\put(350,260){\fbox{$(k,K)$}}

\putmorphism(-150,-400)(1,0)[(\tilde B,A)`(\tilde B,A') `(\tilde B, K)]{640}1a
 \putmorphism(450,-400)(1,0)[\phantom{A'\ot B'}` (\tilde B',A') `(l,A')]{680}1a

\putmorphism(-180,50)(0,-1)[\phantom{Y_2}``(u,A)]{450}1l
\putmorphism(450,50)(0,-1)[\phantom{Y_2}``]{450}1r
\putmorphism(300,50)(0,-1)[\phantom{Y_2}``(u,A')]{450}0r
\putmorphism(1100,50)(0,-1)[\phantom{Y_2}``(v,A')]{450}1r
\put(-20,-180){\fbox{$(u,K)$}}
\put(660,-180){\fbox{$(\omega,A')$}}

\efig}
=
\scalebox{0.86}{
\bfig
 \putmorphism(-150,500)(1,0)[(B,A)`(B', A)`(k, A)]{600}1a
 \putmorphism(450,500)(1,0)[\phantom{A\ot B}`(B', A') `(B', K)]{680}1a

 \putmorphism(-150,50)(1,0)[(\tilde B,A)`(\tilde B',A)`(l, A)]{600}1a
 \putmorphism(450,50)(1,0)[\phantom{A\ot B}`(\tilde B',A') `(\tilde B',K)]{680}1a
\putmorphism(-180,500)(0,-1)[\phantom{Y_2}``(u,A)]{450}1l
\putmorphism(450,500)(0,-1)[\phantom{Y_2}``]{450}1l
\putmorphism(610,500)(0,-1)[\phantom{Y_2}``(v,A)]{450}0l 
\putmorphism(1120,500)(0,-1)[\phantom{Y_3}``(v,A')]{450}1r
\put(-40,270){\fbox{$(\omega,A)$}} 
\put(650,270){\fbox{$(v,K)$}}
\putmorphism(-150,-400)(1,0)[(\tilde B,A)`(\tilde B,A') `(\tilde B, K)]{640}1a
 \putmorphism(450,-400)(1,0)[\phantom{A'\ot B'}` (\tilde B',A') `(l,A')]{680}1a

\putmorphism(-180,50)(0,-1)[\phantom{Y_2}``=]{450}1l
\putmorphism(1120,50)(0,-1)[\phantom{Y_3}``=]{450}1r
\put(300,-200){\fbox{$(l,K)$}}

\efig}
$$

\noindent $\bullet$ \quad ($(k,K)$-r-nat) \label{k,K-r-nat} 
$$
\scalebox{0.86}{
\bfig
 \putmorphism(-150,500)(1,0)[(B,A)`(B', A)`(k, A)]{600}1a
 \putmorphism(450,500)(1,0)[\phantom{A\ot B}`(B', A') `(B', K)]{680}1a
 \putmorphism(-150,50)(1,0)[(B,A)`(B,A')`(B,K)]{600}1a
 \putmorphism(450,50)(1,0)[\phantom{A\ot B}`(B', A') `(k, A')]{680}1a

\putmorphism(-180,500)(0,-1)[\phantom{Y_2}``=]{450}1r
\putmorphism(1100,500)(0,-1)[\phantom{Y_2}``=]{450}1r
\put(350,260){\fbox{$(k,K)$}}

\putmorphism(-180,50)(0,-1)[\phantom{Y_2}``(B,U)]{450}1l
\putmorphism(450,50)(0,-1)[\phantom{Y_2}``]{450}1r
\putmorphism(300,50)(0,-1)[\phantom{Y_2}``(B,V)]{450}0r
\putmorphism(1100,50)(0,-1)[\phantom{Y_2}``(B',V)]{450}1r
\put(-20,-180){\fbox{$(B,\zeta)$}}
\put(660,-180){\fbox{$(k,V)$}}

\putmorphism(-150,-400)(1,0)[(B,\tilde A)`(B,\tilde A') `(B,L)]{640}1a
 \putmorphism(450,-400)(1,0)[\phantom{A'\ot B'}` (B',\tilde A') `(k,\tilde A')]{680}1a
\efig}
=
\scalebox{0.86}{
\bfig
 \putmorphism(-150,500)(1,0)[(B,A)`(B', A)`(k, A)]{600}1a
 \putmorphism(450,500)(1,0)[\phantom{A\ot B}`(B', A') `(B', K)]{680}1a

 \putmorphism(-150,50)(1,0)[(B,\tilde A)`(B',\tilde A)`(k,\tilde A)]{600}1a
 \putmorphism(450,50)(1,0)[\phantom{A\ot B}`(B',\tilde A') `(B',L)]{680}1a

\putmorphism(-180,500)(0,-1)[\phantom{Y_2}``(B,U)]{450}1l
\putmorphism(450,500)(0,-1)[\phantom{Y_2}``]{450}1l
\putmorphism(610,500)(0,-1)[\phantom{Y_2}``(B',U)]{450}0l 
\putmorphism(1120,500)(0,-1)[\phantom{Y_3}``(B',V)]{450}1r
\put(-40,270){\fbox{$(k,U)$}} 
\put(650,270){\fbox{$(B', \zeta)$}}

\putmorphism(-180,50)(0,-1)[\phantom{Y_2}``=]{450}1l
\putmorphism(1120,50)(0,-1)[\phantom{Y_3}``=]{450}1r
\put(300,-200){\fbox{$(k,L)$}}

\putmorphism(-150,-400)(1,0)[(B,\tilde A)`(B,\tilde A') `(B,L)]{640}1a
 \putmorphism(450,-400)(1,0)[\phantom{A'\ot B'}` (B',\tilde A') `(k,\tilde A')]{680}1a
\efig}
$$

\noindent $\bullet$ \quad ($(u,U)$-l-nat) \label{u,U-l-nat} 
$$
\scalebox{0.86}{
\bfig
 \putmorphism(-150,500)(1,0)[(B,A)`(B,A) `=]{600}1a
 \putmorphism(450,500)(1,0)[(B,A)` `(k,A)]{450}1a
\putmorphism(-180,500)(0,-1)[\phantom{Y_2}`(B, \tilde A) `(B,U)]{450}1l
\put(0,50){\fbox{$(u,U)$}}
\putmorphism(-150,-400)(1,0)[(\tilde B, \tilde A)` `=]{500}1a
\putmorphism(-180,50)(0,-1)[\phantom{Y_2}``(u,\tilde A)]{450}1l
\putmorphism(450,50)(0,-1)[\phantom{Y_2}`(\tilde B, \tilde A)`(\tilde B, U)]{450}1l
\putmorphism(450,500)(0,-1)[\phantom{Y_2}`(\tilde B, A) `(u,A)]{450}1l
\put(600,260){\fbox{$(\omega,A)$}}
\putmorphism(450,50)(1,0)[\phantom{(B, \tilde A)}``(l, A)]{500}1a
\putmorphism(1070,50)(0,-1)[\phantom{(B, A')}`(\tilde B', \tilde A)`(\tilde B',U)]{450}1r
\putmorphism(1070,500)(0,-1)[(B', A)`(\tilde B', A)`(v,A)]{450}1r
\putmorphism(450,-400)(1,0)[\phantom{(B, \tilde A)}``(l, \tilde A)]{500}1a
\put(600,-170){\fbox{$(l,U)$}}
\efig}=
\scalebox{0.86}{
\bfig
 \putmorphism(-150,500)(1,0)[(B,A)`(B',A) `(k,A)]{600}1a
 \putmorphism(450,500)(1,0)[\phantom{(B,A)}` `=]{450}1a
\putmorphism(-180,500)(0,-1)[\phantom{Y_2}`(B, \tilde A) `(B,U)]{450}1l
\put(620,50){\fbox{$(v,U)$}}
\putmorphism(-150,-400)(1,0)[(\tilde B, \tilde A)` `(l, \tilde A)]{500}1a
\putmorphism(-180,50)(0,-1)[\phantom{Y_2}``(u,\tilde A)]{450}1l
\putmorphism(450,50)(0,-1)[\phantom{Y_2}`(\tilde B', \tilde A)`(v,\tilde A)]{450}1r
\putmorphism(450,500)(0,-1)[\phantom{Y_2}`(B', \tilde A) `(B',U)]{450}1r
\put(0,260){\fbox{$(k,U)$}}
\putmorphism(-150,50)(1,0)[\phantom{(B, \tilde A)}``(k, \tilde A)]{500}1a
\putmorphism(1070,50)(0,-1)[\phantom{(B, A')}`(\tilde B', \tilde A)`(\tilde B',U)]{450}1r
\putmorphism(1070,500)(0,-1)[(B', A)`(\tilde B', A)`(v,A)]{450}1r
\putmorphism(450,-400)(1,0)[\phantom{(B, \tilde A)}``=]{500}1b
\put(0,-170){\fbox{$(\omega, \tilde A)$}}
\efig}
$$

\noindent $\bullet$ \quad ($(u,U)$-r-nat) \label{u,U-r-nat} 
$$
\scalebox{0.86}{
\bfig
 \putmorphism(-150,500)(1,0)[(B,A)`(B,A) `=]{600}1a
 \putmorphism(550,500)(1,0)[` `(B,K)]{400}1a
\putmorphism(-180,500)(0,-1)[\phantom{Y_2}`(B, \tilde A) `(B,U)]{450}1l
\put(0,50){\fbox{$(u,U)$}}
\putmorphism(-150,-400)(1,0)[(\tilde B, \tilde A)` `=]{500}1a
\putmorphism(-180,50)(0,-1)[\phantom{Y_2}``(u,\tilde A)]{450}1l
\putmorphism(450,50)(0,-1)[\phantom{Y_2}`(\tilde B, \tilde A)`(\tilde B, U)]{450}1l
\putmorphism(450,500)(0,-1)[\phantom{Y_2}`(\tilde B, A) `(u,A)]{450}1l
\put(620,280){\fbox{$(u,K)$}}
\putmorphism(450,50)(1,0)[\phantom{(B, \tilde A)}``(\tilde B,K)]{500}1a
\putmorphism(1070,50)(0,-1)[\phantom{(B, A')}`(\tilde B, \tilde A')`(\tilde B,V)]{450}1r
\putmorphism(1070,500)(0,-1)[(B, A')`(\tilde B, A')`(u,A')]{450}1r
\putmorphism(450,-400)(1,0)[\phantom{(B, \tilde A)}``(\tilde B, L)]{500}1a
\put(620,-170){\fbox{$ (\tilde{B},\zeta)$ } } 
\efig}=
\scalebox{0.86}{
\bfig
 \putmorphism(-150,500)(1,0)[(B,A)`(B,A') `(B,K)]{600}1a
 \putmorphism(450,500)(1,0)[\phantom{(B,A)}` `=]{450}1a
\putmorphism(-180,500)(0,-1)[\phantom{Y_2}`(B, \tilde A) `(B,U)]{450}1l
\put(620,50){\fbox{$(u,V)$}}
\putmorphism(-150,-400)(1,0)[(\tilde B, \tilde A)` `(\tilde B, L)]{500}1a
\putmorphism(-180,50)(0,-1)[\phantom{Y_2}``(u,\tilde A)]{450}1l
\putmorphism(450,50)(0,-1)[\phantom{Y_2}`(\tilde B, \tilde A')`(u,\tilde A')]{450}1r
\putmorphism(450,500)(0,-1)[\phantom{Y_2}`(B, \tilde A') `(B,V)]{450}1r
\put(0,260){\fbox{$(B,\zeta)$}}
\putmorphism(-150,50)(1,0)[\phantom{(B, \tilde A)}``(B,L)]{500}1a
\putmorphism(1070,50)(0,-1)[\phantom{(B, A')}`(\tilde B, \tilde A')`(\tilde B,V)]{450}1r
\putmorphism(1070,500)(0,-1)[(B, A')`(\tilde B, A')`(u,A')]{450}1r
\putmorphism(450,-400)(1,0)[\phantom{(B, \tilde A)}``=]{500}1b
\put(0,-170){\fbox{$(u,L)$}}
\efig}
$$
for any 2-cells 
\begin{equation} \eqlabel{omega-zeta}
\scalebox{0.86}{
\bfig
\putmorphism(-150,50)(1,0)[B` B'`k]{450}1a
\putmorphism(-150,-300)(1,0)[\tilde B`\tilde B' `l]{440}1b
\putmorphism(-170,50)(0,-1)[\phantom{Y_2}``u]{350}1l
\putmorphism(280,50)(0,-1)[\phantom{Y_2}``v]{350}1r
\put(0,-140){\fbox{$\omega$}}
\efig}
\quad\text{and}\quad
\scalebox{0.86}{
\bfig
\putmorphism(-150,50)(1,0)[A` A'`K]{450}1a
\putmorphism(-150,-300)(1,0)[\tilde A`\tilde A' `L]{440}1b
\putmorphism(-170,50)(0,-1)[\phantom{Y_2}``U]{350}1l
\putmorphism(280,50)(0,-1)[\phantom{Y_2}``V]{350}1r
\put(0,-140){\fbox{$\zeta$}}
\efig}
\end{equation}
in $\Bb$, respectively $\Aa$. 
\end{prop}

\begin{proof}
The images of the four types of cells in $\Aa$, which we typically denote as $A,K,U,\zeta$, by the lax double functor 
$\F:\Aa\to\llbracket\Bb,\Cc\rrbracket$ are being denoted by $\F(x)=(-,x)$, for any of such cells $x$ in $\Aa$. 
Then one first sees that $(-,A):\Bb\to\Cc$ is a lax double functor. That $(B,-):\Aa\to\Cc$ is a lax double functor 
follows from the eight axioms of $\F$ as a lax double functor in the following way. Axioms (lx.f.v1) and (lx.f.v2) 
of $\F$ are equalities of vertical lax transformations (v.l.t.). When evaluated at $B$ (this corresponds to the part 1. of $(-,U)$, 
respectively of $(-,1^A)$, being a v.l.t.), they yield axioms (lx.f.v1) and (lx.f.v2) for a lax double functor $(B,-)$. 
The resting six axioms of $\F$ as a lax double functor are equalities of modifications in $\llbracket\Bb,\Cc\rrbracket$, 
and as such we only may evaluate them at $B$. It is evaluating them at $B$ that we cover the resting six axioms for 
$(B,-)$ to be a lax double functor $(B,-)$. 
The origin of the four 2-cells and each of the axioms obtained in part 2. of this Proposition we summarize 
in Table 1. 

\begin{table}[h!]
\begin{center}
\begin{tabular}{ c c } 
New axiom & \hspace{0,8cm} Origin from $\F:\Aa\to\llbracket\Bb,\Cc\rrbracket$ \\ [0.5ex]
\hline
2-cell $(k,K)$ & part 3. of $(-,K)$ being a h.o.t. 
\\ [1ex]   
2-cell $(u,K)$ & part 2. of $(-,K)$ being a h.o.t.  \\ [1ex]
2-cell $(k,U)$ & part 2. of $(-,U)$ being a v.l.t.  \\ [1ex]
2-cell $(u,U)$ & part 3. of $(-,U)$ being a v.l.t.  \\ [1ex]
\hline
$((1_B,K))$  & (h.o.t.-2) of $(-,K)$ \\ [1ex]   
$((k,1_A))$  & (m.ho-vl.-1) of unitor $\F_A: \Id_{(-,A)}\Rrightarrow (-,1_A)$ \\ [1ex]   
$((1^B,K))$  & (h.o.t.-4) of $(-,K)$ \\ [1ex]   
$((u,1_A))$  & (m.ho-vl.-2) of unitor $\F_A: \Id_{(-,A)}\Rrightarrow (-,1_A)$ \\ [1ex]   
$((1_B,U))$  & (v.l.t.-2) of $(-,U)$ \\ [1ex]   
$((k,1^A))$  & (lx.f.v2) of $\F$ (is an equality of v.l.t.) evaluated at $k$ \\ 
           &    ({\em i.e.} part 2. of $(-,1^A)$ being a v.l.t.) \\ [1ex]   
$((1^B,U))$  & (v.l.t.-4) of $(-,U)$ \\ [1ex]   
$((u,1^A))$  & (lx.f.v2) of $\F$ (is an equality of v.l.t.) evaluated at $u$ \\ 
           &    ({\em i.e.} part 3. of $(-,1^A)$ being a v.l.t.) \\ [1ex]   
\hline
$((k'k,K))$  & (h.o.t.-1) of $(-,K)$ \\ [1ex]   
$((k,K'K))$  & (m.ho-vl.-1) of assotiator $\F_{LK}: (-,L)(-,K)\Rrightarrow (-,LK)$ \\ [1ex]   
$((\frac{u}{u'},K))$  & (h.o.t.-3) of $(-,K)$ \\ [1ex]   
$((u,K'K))$  & (m.ho-vl.-2) of assotiator $\F_{LK}: (-,L)(-,K)\Rrightarrow (-,LK)$ \\ [1ex]   
$((k'k,U))$  & (v.l.t.-1) of $(-,U)$ \\ [1ex]   
$((k,\frac{U}{U'}))$  & (lx.f.v1) of $\F$ (is an equality of v.l.t.) evaluated at $k$ \\ 
           &    ({\em i.e.} part 2. of $(-,U)$ being a v.l.t.) \\ [1ex] 
$((u,\frac{U}{U'}))$  & (lx.f.v1) of $\F$ (is an equality of v.l.t.) evaluated at $u$ \\ 
           &    ({\em i.e.} part 3. of $(-,U)$ being a v.l.t.) \\ [1ex] 
$((\frac{u}{u'},U))$  & (v.l.t.-3) of $(-,U)$ \\ [1ex]   
\hline
($(k,K)$-l-nat)  & (h.o.t.-5) of $(-,K)$ \\ [1ex]   
($(k,K)$-r-nat)  & (m.ho-vl.-1) of $(-,\zeta)$ \\ [1ex]   
($(u,U)$-l-nat)  & (v.l.t.-5) of $(-,U)$ \\ [1ex]   
\hline
($(u,U)$-r-nat)  & (m.ho-vl.-2) of $(-,\zeta)$ \\ [1ex]   
\end{tabular}
\caption{Generation of a lax double quasi-functor $\Aa\times\Bb\protect\to\Cc$}
\label{table:1}
\end{center}
\end{table}
\qed\end{proof}

In analogy to \cite[Definition I.4.1]{Gray} we set: 

\begin{defn} \delabel{H dbl}
A {\em lax double quasi-functor} $H: \Aa\times\Bb\to\Cc$ consists of: 
\begin{enumerate} 
\item two families of lax double functors 
$$(-,A):\Bb\to\Cc\quad\text{ and}\quad (B,-):\Aa\to\Cc$$ 
such that $H(A,-)=(-, A), H(-, B)=(B,-)$ and $(-,A)\vert_B=(B,-)\vert_A=(B,A)$, for objects $A\in\Aa, B\in\Bb$, and 
\item four families of 2-cells $(k,K), (u,K), (k,U), (u,U)$ in $\Cc$ for 1h-cells $K$ of $\Aa$ and $k$ of $\Bb$, and 
1v-cells $U$ of $\Aa$ and $u$ of $\Bb$,  
\end{enumerate}
satisfying the conditions listed in part 2. of \prref{char df}. 
\end{defn}

From the data in the above Proposition we may draw several consequences.

\begin{cor} \colabel{obtained oplax trans}
For any 1h-cell $K:A\to A'$, 1v-cell $U:A\to\tilde A$ and 2-cell $\zeta$ in $\Aa$, and 
for any 1h-cell $k:B\to B'$, 1v-cell $u:B\to\tilde B$ and 2-cell $\omega$ in $\Bb$, the following hold: 
\begin{enumerate}
\item $(-,K): (-,A)\to(-,A')$ is a horizontal oplax transformation, 
$(-,U): (-,A)\to(-,\tilde A)$ is a vertical lax transformation, and 
$(-,\zeta)$ is a modification with respect to horizontally oplax and vertically lax transformations, and 
\item $(k,-): (B,-)\to (B', -)$ is a horizontal lax transformation, 
$(u,-): (B,-)\to(\tilde B,-)$ is a vertical oplax transformation, and 
$(\omega,-)$ is a modification with respect to horizontally lax and vertically oplax transformations.
\end{enumerate}
\end{cor}

\begin{proof}
Part 1. highlights the meta-results from the above Proposition: $(-,K), (-,U), (-,\zeta)$ are images of $\F$. 

That $(k,-)$ is a horizontal lax transformation, it follows from: $((k,K'K))$, $((k,1_A))$, $((k,\frac{U}{U'}))$, 
$((k,1^A))$, ($(k,K)$-r-nat) of \prref{char df}. 

That $(u,-)$ is a vertical oplax transformation, it follows from: $((u,K'K))$, $((u,1_A))$, $((u,\frac{U}{U'}))$,
$((u,1^A))$, ($(u,U)$-r-nat). 

That $(\omega,-)$ is a modification in the sense of \deref{mod-hl-vo}, it follows from: ($(k,K)$-l-nat) and ($(u,U)$-l-nat). 
\qed\end{proof}

\subsection{A candidate for a lax Gray type monoidal product} \sslabel{Gray prod}

We may now describe $\Aa\ot\Bb$ by reading off the structure of the image double category $F(\Aa)(\Bb)$  
for any double functor $F:\Aa\to\llbracket\Bb,\Aa\times\Bb\rrbracket$ using the definition of a lax double functor. 
Namely, the result of $\F(\Aa)(\Bb)$ are pairs $(y,x)$ living in the Cartesian product $\Aa\times\Bb$ for any 
0-, 1h-, 1v- or 2-cells $x$ of $\Aa$ and $y$ of $\Bb$. By setting $x\ot y:=(y,x)$ we come to the following definition.  
(Recall that a double category can be seen as a category internal in the category of categories. In this viewpoint, 
we denote the source and target, composition and unit functors by $s,t,c,i$, respectively.) 

\begin{defn} 
Let $\Aa\ot\Bb$ be generated as a double category by the following data: \\
\ul{objects}: $A\ot B$ for objects $A\in\Aa, B\in\Bb$; \\ 
\ul{1h-cells}: $A\ot k, K\ot B$; \\ 
%
\ul{1v-cells}: $A\ot u, U\ot B$ and vertical compositions of such obeying the following rules: 
$$\frac{A\ot u}{A\ot u'}=A\ot \frac{u}{u'}, \quad \frac{U\ot B}{U'\ot B}=\frac{U}{U'}\ot B, \quad A\ot 1^B=1^{A\ot B}=1^A\ot B$$
where $u,u'$ are 1v-cells of $\Bb$ and $U,U'$ 1v-cells of $\Aa$; \\
\ul{2-cells}: $A\ot\omega, \zeta\ot B$:
$$
\scalebox{0.86}{
\bfig
\putmorphism(-150,50)(1,0)[A\ot B`A\ot B'`A\ot k]{600}1a
\putmorphism(-150,-400)(1,0)[A\ot\tilde B`A'\ot\tilde B' `A\ot l]{640}1a
\putmorphism(-180,50)(0,-1)[\phantom{Y_2}``A\ot u]{450}1l
\putmorphism(450,50)(0,-1)[\phantom{Y_2}``A\ot v]{450}1r
\put(0,-180){\fbox{$A\ot \omega$}}
\efig}
\qquad
\scalebox{0.86}{
\bfig
\putmorphism(-150,50)(1,0)[A\ot B`A'\ot B`K\ot B]{600}1a
\putmorphism(-150,-400)(1,0)[\tilde A\ot B`\tilde A'\ot B `L\ot B]{640}1a
\putmorphism(-180,50)(0,-1)[\phantom{Y_2}``U\ot B]{450}1l
\putmorphism(450,50)(0,-1)[\phantom{Y_2}``V\ot B]{450}1r
\put(0,-180){\fbox{$\zeta\ot B$}}
\efig}
$$
where $\omega$ and $\zeta$ are as in \equref{omega-zeta}, 
four (vertically) globular 2-cells from the laxity of double functors $(-,A)$ and $(B,-)$: 
\begin{equation} \eqlabel{4 laxity 2-cells}
(A\ot k')(A\ot k)\stackrel{(A\ot -)_{k'k}}{\Rightarrow} A\ot (k' k), \quad 
(K'\ot B)(K\ot B)\stackrel{(-\ot B)_{K'K}}{\Rightarrow}(K' K)\ot B
\end{equation}
$$1_{A\ot B}\stackrel{(A\ot -)_B}{\Rightarrow} A\ot 1_B, 
\quad 1_{A\ot B}\stackrel{(-\ot B)_A}{\Rightarrow} 1_A\ot B$$
which satisfy associativity and unitality laws, and 
where $k,k'$ are 1h-cells of $\Bb$ and $K, K'$ 1h-cells of $\Aa$, 
and four types of 2-cells coming from the 2-cells of point 2. in \prref{char df}: 
a vertically globular 2-cell $K\ot k: (A'\ot k)(K\ot B) \Rightarrow(K\ot B')(A\ot k)$, 
a horizontally globular 2-cell $U\ot u: \frac{U\ot B}{\tilde A\ot u}\Rightarrow\frac{A\ot u}{U\ot\tilde B}$ 
(so that $1^A\ot 1^B=1_{A\ot B}$), and 
2-cells $K\ot u$ and $U\ot k$, 
subject to the rules induced by the rules of point 2. in \prref{char df} and the following ones: 
$$
\scalebox{0.86}{
\bfig
\putmorphism(-150,500)(1,0)[A\ot B`A\ot B' `A\ot k]{600}1a
 \putmorphism(450,500)(1,0)[\phantom{F(A)}`A\ot B'' `A\ot k' ]{620}1a

 \putmorphism(-150,50)(1,0)[A\ot \tilde B`A\ot \tilde B' ` A\ot l]{600}1a
 \putmorphism(470,50)(1,0)[\phantom{F(A)}`A\ot \tilde B'' `A\ot l']{620}1a

\putmorphism(-180,500)(0,-1)[\phantom{Y_2}``A\ot u]{450}1l
\putmorphism(450,500)(0,-1)[\phantom{Y_2}``]{450}1r
\putmorphism(300,500)(0,-1)[\phantom{Y_2}``A\ot v]{450}0r
\putmorphism(1080,500)(0,-1)[\phantom{Y_2}``A\ot w]{450}1r
\put(-40,280){\fbox{$A\ot\omega$}}
\put(620,280){\fbox{$A\ot\omega'$}}

\putmorphism(-150,-400)(1,0)[A\ot \tilde B`A\ot \tilde B'' `A\ot l'l]{1200}1a

\putmorphism(-180,50)(0,-1)[\phantom{Y_2}``=]{450}1l
\putmorphism(1080,50)(0,-1)[\phantom{Y_3}``=]{450}1r
\put(240,-160){\fbox{$(A\ot -)_{l'l}$}}

\efig}
=
\scalebox{0.86}{
\bfig
\putmorphism(-150,500)(1,0)[A\ot B`A\ot B' `A\ot k]{600}1a
 \putmorphism(450,500)(1,0)[\phantom{F(A)}`A\ot B'' `A\ot k' ]{620}1a
 \putmorphism(-150,50)(1,0)[A\ot B`A\ot B''`A\ot k'k]{1220}1a

\putmorphism(-180,500)(0,-1)[\phantom{Y_2}``=]{450}1r
\putmorphism(1080,500)(0,-1)[\phantom{Y_2}``=]{450}1r
\put(240,290){\fbox{$(A\ot -)_{k'k}$}}

\putmorphism(-150,-400)(1,0)[A\ot \tilde B`A\ot \tilde B'' `A\ot l'l]{1200}1a

\putmorphism(-180,50)(0,-1)[\phantom{Y_2}``A\ot u]{450}1l
\putmorphism(1080,50)(0,-1)[\phantom{Y_3}``A\ot w]{450}1r
\put(300,-180){\fbox{$A\ot\omega'\omega$}} 
\efig}
$$

$$
\scalebox{0.86}{
\bfig
\putmorphism(-150,500)(1,0)[A\ot B`A'\ot B`K\ot B]{600}1a
 \putmorphism(470,500)(1,0)[\phantom{F(A)}`A''\ot B `K'\ot B]{600}1a
 \putmorphism(-150,50)(1,0)[A\ot B`A''\ot B`K'K\ot B]{1220}1a

\putmorphism(-180,500)(0,-1)[\phantom{Y_2}``=]{450}1r
\putmorphism(1080,500)(0,-1)[\phantom{Y_2}``=]{450}1r
\put(240,290){\fbox{$(-\ot B)_{K'K}$}}

\putmorphism(-150,-400)(1,0)[\tilde A\ot B`\tilde A''\ot B `L'L\ot B]{1200}1a

\putmorphism(-180,50)(0,-1)[\phantom{Y_2}``U\ot B]{450}1l
\putmorphism(1080,50)(0,-1)[\phantom{Y_3}``U''\ot B]{450}1r
\put(300,-180){\fbox{$\zeta'\zeta\ot B$}} 
\efig}
=
\scalebox{0.86}{
\bfig
\putmorphism(-150,500)(1,0)[A\ot B`A'\ot B`K\ot B]{600}1a
 \putmorphism(470,500)(1,0)[\phantom{F(A)}`A''\ot B `K'\ot B]{600}1a

 \putmorphism(-150,50)(1,0)[\tilde A\ot B`\tilde A'\ot B`L\ot B]{600}1a
 \putmorphism(470,50)(1,0)[\phantom{F(A)}`\tilde A''\ot B `L' \ot B]{620}1a

\putmorphism(-180,500)(0,-1)[\phantom{Y_2}``U\ot B]{450}1l
\putmorphism(450,500)(0,-1)[\phantom{Y_2}``]{450}1r
\putmorphism(300,500)(0,-1)[\phantom{Y_2}``U'\ot B]{450}0r
\putmorphism(1080,500)(0,-1)[\phantom{Y_2}``U''\ot B]{450}1r
\put(-40,280){\fbox{$\zeta\ot B$}}
\put(620,280){\fbox{$\zeta' \ot B$}}

\putmorphism(-150,-400)(1,0)[\tilde A\ot B`\tilde A''\ot B `L'L\ot B]{1200}1a

\putmorphism(-180,50)(0,-1)[\phantom{Y_2}``=]{450}1l
\putmorphism(1080,50)(0,-1)[\phantom{Y_3}``=]{450}1r
\put(240,-160){\fbox{$(-\ot B)_{L'L}$}}

\efig}
$$

$$A\ot\frac{\omega}{\omega'}=\frac{A\ot\omega}{A\ot\omega'}, \quad \frac{\zeta}{\zeta'}\ot B=\frac{\zeta\ot B}{\zeta'\ot B},$$
$$A\ot\Id_k=\Id_{A\ot k}, \quad \Id_K\ot B=\Id_{K\ot B}, \quad A\ot\Id^u=\Id^{A\ot u}, \quad \Id^U\ot B=\Id^{U\ot B} .$$
The source and target functors $s,t$ on $\Aa\ot\Bb$ are defined as in $\Aa\times\Bb$, the composition functor $c$ is defined by horizontal juxtaposition of 
the corresponding 2-cells, and the unit functor $i$ is defined on generators as follows: 
$$i(A\ot B)=1_{A\ot B}, \, i(A\ot v)=1^A\ot v (=Id^{A\ot v}) \quad \text{and} \quad i(U\ot B)=U\ot 1^B (=Id^{U\ot B}).$$ 
\end{defn}

Since $\Aa\ot\Bb$ is defined by generators and relations on $\Aa\times\Bb$, 
it is clear that there is a lax double quasi-functor $J:\Aa\times\Bb\to\Aa\ot\Bb$ given by $J(-, B)(x)=x\ot B, J(A,-)(y)=A\ot y$ for cells 
$x$ in $\Aa$ and $y$ in $\Bb$ and with unique 2-cells $K\ot k, K\ot u, U\ot k$ and $U\ot u$ in $\Aa\ot\Bb$, where the usual notation 
is used. It turns out that the universal property that $\Aa\ot\Bb$ satisfies is the following: 
for every lax double quasi-functor $H:\Aa\times\Bb\to\Cc$ there is a unique 
strict double functor $\crta H: \Aa\ot\Bb\to\Cc$ such that $H=\crta{H}J$. Moreover, in \ssref{Gray} we will prove a double category isomorphism: 
$$q\x\Lax_{hop}(\Aa\times\Bb,\Cc)\iso \Dbl_{hop}(\Aa\ot\Bb,\Cc).$$ 

\subsection{Why $-\ot-$ fails to be a (lax Gray type) monoidal product} \sslabel{not Gray}

As we announced at the beginning of \seref{mon str}, we will show now why $\Aa\ot\Bb$ does not yield a Gray monoidal product in the 
category of double categories and lax double functors $Dbl^{st}_{lx}$. The fact that we work with lax instead of pseudo or strict 
double functors makes two important steps in the construction of a Gray monoidal product fail. On one hand, 
it is clear that the only kind of ``isomorphism'' between  
double categories $(\Aa\ot\Bb)\ot\Cc$ and $\Aa\ot(\Bb\ot\Cc)$ must be (an invertible) pseudo double functor. On the other hand, 
in \ssref{not-inner} we explained that 
$\llbracket -,-\rrbracket: (Dbl_{lx})^{op}\times Dbl_{lx}\to Dbl_{lx}$ is not a functor, so we can abandon the idea 
to get closedness of $Dbl_{lx}$ in the expected way. 
Nevertheless, one could still ask the question whether $\Aa\ot\Bb$ gives a monoidal product on some category of double categories.  
Let us explore this question.  

\medskip

We find that the main issue that hinders proving monoidality of $\ot$ is its ``laxity'' as it is 
constructed in \ssref{Gray prod}. Concretely: the fact that the 2-cells 
$(k,K)$ in \prref{char df} are not invertible and the laxity of double functors $(-,A)$ and $(B,-)$, 
which then imply that the 2-cells $K\ot k$ and 
$$
(A\ot k')(A\ot k)\stackrel{(A\ot -)_{k'k}}{\Rightarrow} A\ot (k' k), \quad 
(K'\ot B)(K\ot B)\stackrel{(-\ot B)_{K'K}}{\Rightarrow}(K' K)\ot B
$$
from \equref{4 laxity 2-cells} in $\Aa\ot\Bb$ are not invertible. 
We show why this is not good for getting a monoidal product. 

\smallskip

Suppose that we want to define this kind of tensor product for two double functors {\em of any kind} among double categories: $F:\Aa\to\Bb$ and $G:\Cc\to\Dd$. 
In order to define $F\ot G:\Aa\ot\Cc\to\Bb\ot\Dd$ we need to define in particular how $F\ot G$ acts on a 1h-cell of type 
$(A'\ot k)(K\ot C):A\ot C\Rightarrow A'\ot C'$.  
(Observe that we only have a 2-cell $(K\ot C')(A\ot k)\Rightarrow(A'\ot k)(K\ot C)$ in $\Aa\ot\Cc$.) The natural candidate for this definition would be: 
$(F\ot G)((A'\ot k)(K\ot C))=(F(A')\ot G(k))(F(K)\ot G(C)).$

But then with the lax Gray type product $\ot$ there is no way to give a compositor 2-cell for example for: 
$$(F\ot G)\left( (L\ot C'')(A'\ot l)\right) (F\ot G)\left( (K\ot C')(A\ot k)\right) \stackrel{\xi}{\Rightarrow}
(F\ot G)\left( (L\ot C'')(A'\ot l)(K\ot C')(A\ot k)\right).$$ 
Though, in the pseudo or strict Gray type product as in \cite{Gabi} one has: 
$$\left(F(A')\ot G(l)\right)\left(F(K)\ot G(C')\right) \stackrel{*}{\iso} \left(F(K)\ot G(C'')\right)\left(F(A)\ot G(l)\right)$$ 
and also isomorphisms 
$$\left(F(L)\ot G(C'')\right)\left(F(K)\ot G(C')\right)\stackrel{*}{\iso} F(LK)\ot G(C')\quad\text{and}\quad 
\left(F(A)\ot G(l)\right)\left(F(A)\ot G(k)\right)\stackrel{*}{\iso} F(A)\ot G(lk)$$
which permit to construct a 2-cell $\xi$ above and finally a pseudodouble functor structure on $F\ot G$. 

Basically, this type of argument explains why strict or pseudo type Gray 
product $\ot$ permits to obtain a monoidal structure on double categories with pseudo or strict double functors, as it was done in \cite{Gabi} in the 
strict case. And likewise, it explains why the non-invertibility of the 2-cells at places analogous to * in our lax Gray case of $\ot$ prevents us from proving monoidality.

\subsection{What about enrichment?} \sslabel{enrich}

Having constructed double categories $\llbracket\Aa,\Bb\rrbracket=\Lax_{hop}(\Aa,\Bb)$ for every two double categories $\Aa, \Bb$, one may wonder if 
the category $Dbl_{lx}$ of double categories and lax double functors is enriched over double categories or 2-categories. 
For this question we find again that laxity of the double functors in $Dbl_{lx}$ makes problems. 

\smallskip

Firstly, assume that the hom-category for $Dbl_{lx}$ and two double categories $\Aa, \Bb$ is given by $\llbracket\Aa,\Bb\rrbracket=\Lax_{hop}(\Aa,\Bb)$ 
and that the composition on them is induced by the composition of lax double functors (0-cells in the hom-categories), 
{\em horizontal} composition of horizontal oplax and vertical lax transformations and of modifications. Let us denote this kind of composition by 
$\circ$. 
 
This requires to give a horizontal composition among horizontal oplax transformation (and among 1v- and 2-cells). 
However, one can only define the horizontal composition of horizontal oplax transformations if the functors have both lax and colax structure. 
This is the same fact that happens in 2-categories. Namely, we encounter exactly the same issue as in \ssref{not-inner}: take for $H$ there to be 
any lax double functor (as we need in our situation), and you see that one can not define $\delta_{H(\alpha),f}$ (as also colaxity of $H$ is needed). 
If one had pseudodouble functors for 0-cells in the hom-categories, then one can define the horizontal composition of horizontal oplax transformations 
as in \cite[Lemma 3.5]{Fem}:

\begin{lma} \lelabel{horiz comp hot}
Horizontal composition of two horizontal oplax transformations $\alpha: F\Rightarrow G: \Aa\to\Bb$ and $\beta: F\s'\Rightarrow G':\Bb\to\Cc$ 
among pseudodouble functors is well-given by: 
\begin{itemize}
\item for every 0-cell $A$ in $\Aa$ a 1h-cell in $\Cc$:
$$(\beta\comp\alpha)(A)=\big( F\s'F(A)\stackrel{F\s'(\alpha(A))}{\longrightarrow}F\s'G(A) \stackrel{\beta(G(A))}{\longrightarrow} G'G(A) \big),$$ 
\item for every 1v-cell $u:A\to A'$ in $\Aa$ a 2-cell in $\Cc$:
$$(\beta\comp\alpha)^u=
\bfig
\putmorphism(-320,500)(1,0)[F\s'F(A)`F\s'G(A)`F\s'(\alpha(A))]{770}1a
 \putmorphism(500,500)(1,0)[\phantom{F(A)}`G'G(A) `\beta(G(A))]{720}1a

 \putmorphism(-320,50)(1,0)[F\s'F(A')`F\s'G(A')`F\s'(\alpha(A'))]{770}1a
 \putmorphism(520,50)(1,0)[\phantom{F(A)}`G'G(A') `\beta(G(A'))]{730}1a

\putmorphism(-280,500)(0,-1)[\phantom{Y_2}``F\s'F(u)]{450}1l
\putmorphism(450,500)(0,-1)[\phantom{Y_2}``]{450}1r
\putmorphism(300,500)(0,-1)[\phantom{Y_2}``F\s'G(u)]{450}0r
\putmorphism(1180,500)(0,-1)[\phantom{Y_2}``G'G(u)]{450}1r
\put(-160,290){\fbox{$F\s'(\alpha^u)$}}
\put(700,300){\fbox{$(\beta)^{G(u)}$}}
\efig
$$
\item for every 1h-cell $f:A\to B$ in $\Aa$ a 2-cell in $\Cc$:
$$\delta_{\beta\comp\alpha,f}=
\bfig
 \putmorphism(-250,0)(1,0)[F\s'F(A)`F\s'F(B)` F\s'F(f)]{700}1a
 \putmorphism(450,0)(1,0)[\phantom{F\s'F(B)}`F\s'G(B) `F\s'(\alpha(B))]{760}1a

 \putmorphism(-250,-450)(1,0)[F\s'F(A)`F\s'G(A)`F\s'(\alpha(A))]{700}1a
 \putmorphism(480,-450)(1,0)[\phantom{A\ot B}`F\s'G(B) `F\s'G(f)]{740}1a
 \putmorphism(1200,-450)(1,0)[\phantom{A'\ot B'}` G'G(B) `\beta(G(B))]{760}1a

\putmorphism(-280,0)(0,-1)[\phantom{Y_2}``=]{450}1r
\putmorphism(1180,0)(0,-1)[\phantom{Y_2}``=]{450}1r
\put(350,-240){\fbox{$ \delta_{F\s'(\alpha),f}  $}}
\put(1000,-700){\fbox{$\delta_{\beta,G(f)}$}}

 \putmorphism(450,-900)(1,0)[F\s'G(A)` G'G(A) `\beta(G(A))]{760}1a
 \putmorphism(1180,-900)(1,0)[\phantom{A''\ot B'}`G'G(B) ` G'G(f)]{760}1a

\putmorphism(450,-450)(0,-1)[\phantom{Y_2}``=]{450}1l
\putmorphism(1950,-450)(0,-1)[\phantom{Y_2}``=]{450}1r
\efig
$$
where $\delta_{F\s'(\alpha),f}$ is from \ssref{not-inner}. 
\end{itemize}
\end{lma}

Moreover, in order to have some kind of double or 2-functor 
$$\circ: \llbracket\Bb,\Cc\rrbracket \times \llbracket\Aa,\Bb\rrbracket \to \llbracket\Aa,\Cc\rrbracket$$
one needs to explore in particular the interchange law on 1h-cells (recall that the composition of 1h-cells in $\llbracket\Aa,\Bb\rrbracket$ is 
given by the vertical composition of horizontal oplax transformations).

In order not to extend the argumentation too much, we recollect our findings in short. 
Let $F\s'\circ F\stackrel{\beta\circ\alpha}{\Rightarrow} G'\circ G \stackrel{\beta'\circ\alpha'}{\Rightarrow}H'\circ H$ 
be horizontal oplax transformations of pseudodouble functors. One obtains: 
$$[\frac{\alpha}{\alpha'} \hspace{0,2cm} \vert \hspace{0,2cm} \frac{\beta}{\beta'}]^u= 
[F\s'([\alpha^u\vert(\alpha')^u])  \hspace{0,2cm} \vert  \hspace{0,2cm} \beta^{H(u)}  \hspace{0,2cm} \vert  \hspace{0,2cm} \beta'^{H(u)}]$$
and 
$$\frac{[\alpha \hspace{0,2cm} \vert \hspace{0,2cm} \beta]^u}{[\alpha' \hspace{0,2cm} \vert \hspace{0,2cm} \beta']^u}= 
[F\s'(\alpha^u) \hspace{0,2cm} \vert \hspace{0,2cm} \beta^{G(u)} \hspace{0,2cm} \vert \hspace{0,2cm} G'((\alpha')^u)\hspace{0,2cm} \vert \hspace{0,2cm} \beta'^{H(u)}].$$
In order for these two to be equal one would need strictness of $F\s'$ (and thus of all double functors) and strictness of $\beta$ (and all the horizontal transformations). Alternatively, pseudodouble functors and horizontal pseudotransformations would allow to build an invertible modification joining the two 
2-cells above. 

\medskip

These findings mean that in our lax construction of $\llbracket\Aa,\Bb\rrbracket$ as $\Lax_{hop}(\Aa,\Bb)$ we do not have enrichment with hom-categories 
given by $\llbracket\Aa,\Bb\rrbracket$ and horizontal composition $\circ$. Though, the inner-hom in \cite{Gabi} (with strict double functors and 
horizontal pseudotransformations) makes $Dbl$, 
the category of double categories and strict double functors, enriched over $Dbl_{ps}$, in which morphisms are pseudotransformations.

\section{The double categories $\Lax_{hop}(\Aa, \llbracket\Bb,\Cc\rrbracket)$ and \\ $q\x\Lax_{hop}(\Aa\times\Bb,\Cc)$ are isomorphic} \selabel{isom}

Observe that a pair of families of lax functors of 2-categories together with their distributive law, 
which is given by a family of 2-cells $\sigma_{f,g}$ for 1-cells $f,g$, defined in \cite[Definition 3.1]{FMS} present a lax version 
of ``quasi-functors of two variables'' of \cite[Definition I, 4.1]{Gray}. Namely, $\sigma_{f,g}$ from the former precisely corresponds 
to $\gamma_{f,g}$ of the latter, only the functors in \cite{Gray} are strict 2-functors. 
The single condition QF$_23$ of the latter is equivalent to the two conditions (D5) and (D6) of the former. 
In \cite{GG,GPS} the 2-cells $\gamma_{f,g}$ of a quasi-functor of two variables 
were considered to be invertible. Such a quasi-functor of two variables in these references was called ``cubical functor''. 
In Proposition 2.1 and Definition 2.2 of \cite{Fem} we generalized cubical functors to strict double categories and called them {\em cubical double functors}. 

In \prref{char df} and \deref{H dbl} above we generalized cubical double functors to the lax case. (Observe that the corresponding 2-cell 
mentioned in the above paragraph is not invertible, so we do not work here in a cubical setting, and follow Gray's terminology.) 
Thus our \prref{char df} is a generalization to the double categorical setting of \cite[Lemma 4.1]{FMS}.

Morphisms of distributive laws of lax functors from \cite[Definition 4.3]{FMS} are the oplax version of quasi-natural transformations from 
\cite[Definition I, 4.1]{Gray}, which are lax (see I.4.1 and I.3.3 of \cite{Gray}). In this Section we will first 
introduce the corresponding notions to horizontal oplax and vertical lax transformations and their modifications in the lax double quasi-functor setting, 
and then prove that the latter are in 1-1 correspondence with the horizontal oplax and vertical lax transformations between lax double functors 
of the form $\Aa\to\llbracket\Bb,\Cc\rrbracket$ and their modifications.

\subsection{The double category $q\x\Lax_{hop}(\Aa\times\Bb,\Cc)$ } \sslabel{2-cat}

By $q\x\Lax_{hop}(\Aa\times\Bb,\Cc)$ we will denote the double category consisting of lax double quasi-functors, horizontal oplax 
transformations of lax double quasi-functors as 1h-cells, vertical lax transformations as 1v-cells and modifications among the latter two. 
We define its 1- and 2-cells below.

\begin{defn} \delabel{oplax tr cubical}
A horizontal oplax transformation $\theta: (-,-)_1\Rightarrow (-,-)_2$ between lax double quasi-functors $(-,-)_1,(-,-)_2: 
\Aa\times\Bb\to\Cc$ is given by: for each $A\in\Aa$ a horizontal oplax transformation $\theta^A: (-,A)_1\Rightarrow(-,A)_2$ and 
for each $B\in\Bb$ a horizontal oplax transformation $\theta^B: (B,-)_1\Rightarrow(B,-)_2$, both of lax double functors, such that 
$\theta^A_B=\theta^B_A$ and such that 

\medskip
\noindent $(HOT^q_1)$ \vspace{-0,6cm}
$$\scalebox{0.86}{
\bfig
 \putmorphism(450,500)(1,0)[(B', A)_1 `(B', A')_1 `(B', K)_1]{680}1a
 \putmorphism(1140,500)(1,0)[\phantom{A\ot B}`(B', A')_2 ` \theta^{B'}_{A'}]{680}1a

 \putmorphism(-150,50)(1,0)[(B, A)_1 `(B', A)_1`(k, A)_1]{600}1a
 \putmorphism(450,50)(1,0)[\phantom{A\ot B}`(B', A)_2 `\theta^{B'}_A]{680}1a
 \putmorphism(1130,50)(1,0)[\phantom{A\ot B}`(B', A')_2 ` (B', K)_2]{680}1a

\putmorphism(450,500)(0,-1)[\phantom{Y_2}``=]{450}1r
\putmorphism(1750,500)(0,-1)[\phantom{Y_2}``=]{450}1r
\put(1020,270){\fbox{$ \theta^{B'}_K$}}

 \putmorphism(-150,-400)(1,0)[(B, A)_1`(B, A)_2 `\theta^A_B]{640}1a
 \putmorphism(460,-400)(1,0)[\phantom{A'\ot B'}`(B', A)_2 `(k, A)_2]{680}1a

\putmorphism(-180,50)(0,-1)[\phantom{Y_2}``=]{450}1l
\putmorphism(1120,50)(0,-1)[\phantom{Y_3}``=]{450}1r
\put(310,-200){\fbox{$ \theta^A_k$}}

 \putmorphism(1170,-400)(1,0)[\phantom{A\ot B}`(B', A')_2 ` (B', K)_2]{650}1a
\putmorphism(450,-400)(0,-1)[\phantom{Y_2}``=]{450}1r
\putmorphism(1750,-400)(0,-1)[\phantom{Y_2}``=]{450}1r

 \putmorphism(460,-850)(1,0)[(B, A)_2 `(B, A')_2 `(B, K)_2]{650}1a
 \putmorphism(1240,-850)(1,0)[ `(B', A')_2 `(k, A')_2]{600}1a
\put(920,-640){\fbox{$ (k,K)_2$}}
\efig}
=
\scalebox{0.86}{
\bfig

 \putmorphism(-150,450)(1,0)[(B,A)_1`(B',A)_1`(k,A)_1]{600}1a
 \putmorphism(450,450)(1,0)[\phantom{A\ot B}`(B', A')_1 `(B',K)_1]{680}1a

 \putmorphism(-150,0)(1,0)[(B,A)_1`(B,A')_1`(B,K)_1]{600}1a
 \putmorphism(450,0)(1,0)[\phantom{A\ot B}`(B', A')_1 `(k,A')_1]{680}1a
 \putmorphism(1120,0)(1,0)[\phantom{A'\ot B'}`(B', A')_2 `\theta^{A'}_{B'}]{650}1a

\putmorphism(-180,450)(0,-1)[\phantom{Y_2}``=]{450}1r
\putmorphism(1100,450)(0,-1)[\phantom{Y_2}``=]{450}1r
\put(350,210){\fbox{$(k,K)_1$}}
\put(1000,-250){\fbox{$\theta^{A'}_k$}}

 \putmorphism(-150,-450)(1,0)[(B,A)_1`(B,A')_1`(B,K)_1]{600}1a
 \putmorphism(450,-450)(1,0)[\phantom{A''\ot B'}` (B, A')_2 `\theta^{A'}_B]{680}1a
 \putmorphism(1100,-450)(1,0)[\phantom{A''\ot B'}`(B', A')_2 ` (k, A')_2]{660}1a

\putmorphism(450,0)(0,-1)[\phantom{Y_2}``=]{450}1l
\putmorphism(1750,0)(0,-1)[\phantom{Y_2}``=]{450}1r
\putmorphism(-180,-450)(0,-1)[\phantom{Y_2}``=]{450}1r
\putmorphism(1100,-450)(0,-1)[\phantom{Y_2}``=]{450}1r
 \putmorphism(-150,-900)(1,0)[(B,A)_1`(B,A)_2`\theta^B_A]{600}1a
 \putmorphism(580,-900)(1,0)[`(B,A')_2`(B,K)_2]{540}1a
\put(320,-670){\fbox{$\theta^{B}_K$}}
\efig}
$$
for every 1h-cells $K:A\to A'$ and $k:B\to B'$, 

\medskip
\noindent $(HOT^q_2)$ \vspace{-0,6cm}
$$\scalebox{0.86}{
\bfig
\putmorphism(-150,500)(1,0)[(B,A)_1`(B,A')_1`(B,K)_1]{600}1a
 \putmorphism(450,500)(1,0)[\phantom{F(A)}`(B,A')_2 `\theta^{A'}_B]{640}1a

 \putmorphism(-150,50)(1,0)[(\tilde B,A)_1`(\tilde B,A')_1`(\tilde B,K)_1]{600}1a
 \putmorphism(450,50)(1,0)[\phantom{F(A)}`(\tilde B,A')_2 `\theta^{A'}_{\tilde B}]{640}1a

\putmorphism(-180,500)(0,-1)[\phantom{Y_2}``(u,A)_1]{450}1l
\putmorphism(450,500)(0,-1)[\phantom{Y_2}``]{450}1r
\putmorphism(300,500)(0,-1)[\phantom{Y_2}``(u,A')_1]{450}0r
\putmorphism(1100,500)(0,-1)[\phantom{Y_2}``(u,A')_2]{450}1r
\put(-40,280){\fbox{$(u,K)_1$}}
\put(700,280){\fbox{$\theta^{A'}_u$}}

\putmorphism(-150,-400)(1,0)[(\tilde B,A)_1`(\tilde B,A)_2 `\theta^{A}_{\tilde B}]{640}1a
 \putmorphism(490,-400)(1,0)[\phantom{F(A')}` (\tilde B,A')_2 `(\tilde B,K)_2]{640}1a

\putmorphism(-180,50)(0,-1)[\phantom{Y_2}``=]{450}1l
\putmorphism(1120,50)(0,-1)[\phantom{Y_3}``=]{450}1r
\put(320,-200){\fbox{$\theta^{\tilde B}_K$}}
\efig}
\quad=\quad
\scalebox{0.86}{
\bfig
\putmorphism(-150,500)(1,0)[(B,A)_1`(B,A')_1`(B,K)_1]{600}1a
 \putmorphism(450,500)(1,0)[\phantom{F(A)}`(B,A')_2 `\theta^{A'}_B]{640}1a
 \putmorphism(-150,50)(1,0)[(B,A)_1`(B,A)_2`\theta^{A}_B]{600}1a
 \putmorphism(470,50)(1,0)[\phantom{F(A)}`(B,A')_2 `(B,K)_2]{660}1a

\putmorphism(-180,500)(0,-1)[\phantom{Y_2}``=]{450}1r
\putmorphism(1100,500)(0,-1)[\phantom{Y_2}``=]{450}1r
\put(320,280){\fbox{$\theta^B_K$}}

\putmorphism(-150,-400)(1,0)[(\tilde B,A)_1`(\tilde B,A)_2 `\theta^{A}_{\tilde B}]{640}1a
 \putmorphism(490,-400)(1,0)[\phantom{F(A')}` (\tilde B,A')_2 `(\tilde B,K)_2]{640}1a

\putmorphism(-180,50)(0,-1)[\phantom{Y_2}``(u,A)_1]{450}1l
\putmorphism(450,50)(0,-1)[\phantom{Y_2}``]{450}1l
\putmorphism(610,50)(0,-1)[\phantom{Y_2}``(u,A)_2]{450}0l 
\putmorphism(1120,50)(0,-1)[\phantom{Y_3}``(u,A')_2]{450}1r
\put(-40,-180){\fbox{$\theta^A_u$}} 
\put(620,-180){\fbox{$(u,K)_2$}}

\efig}
$$
for every 1h-cell $K:A\to A'$ and 1v-cell $u: B\to \tilde B$, 

\medskip
\noindent $(HOT^q_3)$ \vspace{-0,6cm}
$$\scalebox{0.86}{
\bfig
\putmorphism(-150,500)(1,0)[(B,A)_1`(B',A)_1`(k,A)_1]{600}1a
 \putmorphism(480,500)(1,0)[\phantom{F(A)}`(B',A)_2 `\theta^A_{B'}]{640}1a

 \putmorphism(-150,50)(1,0)[(B,\tilde A)_1`(B',\tilde A)_1`(k,\tilde A)_1]{600}1a
 \putmorphism(480,50)(1,0)[\phantom{F(A)}`(B',\tilde A)_2 `\theta^{\tilde A}_{B'}]{640}1a

\putmorphism(-180,500)(0,-1)[\phantom{Y_2}``(B,U)_1]{450}1l
\putmorphism(450,500)(0,-1)[\phantom{Y_2}``]{450}1r
\putmorphism(300,500)(0,-1)[\phantom{Y_2}``(B',U)_1]{450}0r
\putmorphism(1100,500)(0,-1)[\phantom{Y_2}``(B',U)_2]{450}1r
\put(-40,280){\fbox{$(k,U)_1$}}
\put(700,280){\fbox{$\theta_U^{B'}$}}

\putmorphism(-150,-400)(1,0)[(B,\tilde A)_1`(B,\tilde A)_2 `\theta^{\tilde A}_B]{640}1a
 \putmorphism(480,-400)(1,0)[\phantom{A'\ot B'}` (B',\tilde A)_2 `(k,\tilde A)_2]{680}1a

\putmorphism(-180,50)(0,-1)[\phantom{Y_2}``=]{450}1l
\putmorphism(1120,50)(0,-1)[\phantom{Y_3}``=]{450}1r
\put(320,-200){\fbox{$\theta^{\tilde A}_k$}}
\efig}
\quad=\quad
\scalebox{0.86}{
\bfig
\putmorphism(-150,500)(1,0)[(B,A)_1`(B',A)_1`(k,A)_1]{600}1a
 \putmorphism(480,500)(1,0)[\phantom{F(A)}`(B',A)_2 `\theta^{A}_{B'}]{640}1a
 \putmorphism(-150,50)(1,0)[(B,A)_1`(B,A)_2`\theta^{A}_B]{600}1a
 \putmorphism(480,50)(1,0)[\phantom{F(A)}`(B',A)_2 `(k,A)_2]{680}1a

\putmorphism(-180,500)(0,-1)[\phantom{Y_2}``=]{450}1r
\putmorphism(1100,500)(0,-1)[\phantom{Y_2}``=]{450}1r
\put(320,280){\fbox{$\theta^A_k$}}

\putmorphism(-150,-400)(1,0)[(B,\tilde A)_1`(B,\tilde A)_2 `\theta^{\tilde A}_{B}]{640}1a
 \putmorphism(490,-400)(1,0)[\phantom{F(A')}` (B',\tilde A)_2 `(k,\tilde A)_2]{640}1a

\putmorphism(-180,50)(0,-1)[\phantom{Y_2}``(B,U)_1]{450}1l
\putmorphism(480,50)(0,-1)[\phantom{Y_2}``]{450}1l
\putmorphism(570,50)(0,-1)[\phantom{Y_2}``(B,U)_2]{450}0l 
\putmorphism(1120,50)(0,-1)[\phantom{Y_3}``(B',U)_2]{450}1r
\put(-40,-180){\fbox{$\theta_U^B$}} 
\put(620,-180){\fbox{$(k,U)_2$}}
\efig}
$$
for every 1v-cell $U:A\to \tilde A$ and 1h-cell $k:B\to B'$, and 

\medskip
\noindent $(HOT^q_4)$ \vspace{-0,6cm}
$$\scalebox{0.86}{
\bfig
 \putmorphism(-150,500)(1,0)[(B,A)_1`(B,A)_1 `=]{600}1a
 \putmorphism(450,500)(1,0)[(B,A)_1` `\theta^{A}_B]{450}1a
\putmorphism(-200,500)(0,-1)[\phantom{Y_2}`(B, \tilde A)_1 `(B,U)_1]{450}1l
\put(-50,50){\fbox{$(u,U)_1$}}
\putmorphism(-170,-400)(1,0)[(\tilde B, \tilde A)_1` `=]{500}1a
\putmorphism(-200,50)(0,-1)[\phantom{Y_2}``(u,\tilde A)_1]{450}1l
\putmorphism(450,50)(0,-1)[\phantom{Y_2}`(\tilde B, \tilde A)_1`(\tilde B, U)_1]{450}1l
\putmorphism(450,500)(0,-1)[\phantom{Y_2}`(\tilde B, A)_1 `(u,A)_1]{450}1l
\put(600,260){\fbox{$\theta^A_u$}}
\putmorphism(450,50)(1,0)[\phantom{(B, \tilde A)}``\theta^{A}_{\tilde B}]{500}1a
\putmorphism(1070,50)(0,-1)[\phantom{(B, A')}`(\tilde B, \tilde A)_2`(\tilde B,U)_2]{450}1r
\putmorphism(1070,500)(0,-1)[(B, A)_2`(\tilde B, A)_2`(u,A)_2]{450}1r
\putmorphism(450,-400)(1,0)[\phantom{(B, \tilde A)}``\theta^{\tilde A}_{\tilde B}]{500}1a
\put(600,-170){\fbox{$\theta_U^{\tilde B}$}}
\efig}=
\scalebox{0.86}{
\bfig
 \putmorphism(-150,500)(1,0)[(B,A)_1`(B,A)_2 `\theta^{A}_B]{600}1a
 \putmorphism(450,500)(1,0)[\phantom{(B,A)}` `=]{540}1a
\putmorphism(-180,500)(0,-1)[\phantom{Y_2}`(B, \tilde A)_1 `(B,U)_1]{450}1l
\put(0,280){\fbox{$\theta_U^B$}}
\putmorphism(-180,-400)(1,0)[(\tilde B, \tilde A)_1` `\theta^{\tilde A}_{\tilde B}]{500}1a
\putmorphism(-180,50)(0,-1)[\phantom{Y_2}``(u,\tilde A)_1]{450}1l
\putmorphism(450,50)(0,-1)[\phantom{Y_2}`(\tilde B, \tilde A)_2`(u,\tilde A)_2]{450}1r
\putmorphism(450,500)(0,-1)[\phantom{Y_2}`(B, \tilde A)_2 `(B,U)_2]{450}1r
\put(600,50){\fbox{$(u,U)_2$}}
\putmorphism(-180,50)(1,0)[\phantom{(B, \tilde A)}``\theta^{\tilde A}_B]{500}1a
\putmorphism(1130,50)(0,-1)[\phantom{(B, A')}`(\tilde B, \tilde A)_2`(\tilde B,U)_2]{450}1r
\putmorphism(1130,500)(0,-1)[(B, A)`(\tilde B, A)_2`(u,A)_2]{450}1r
\putmorphism(450,-400)(1,0)[\phantom{(B, \tilde A)}``=]{520}1b
\put(0,-170){\fbox{$\theta^{\tilde A}_u$}}
\efig}
$$
for every 1v-cells $U:A\to \tilde A$ and $u:B\to\tilde B$. 
\end{defn}


\begin{defn} \delabel{lax v tr cubical}
A vertical lax transformation $\theta_0: (-,-)_1\Rightarrow (-,-)_2$ between lax double quasi-functors $(-,-)_1,(-,-)_2: 
\Aa\times\Bb\to\Cc$ is given by: for each $A\in\Aa$ a vertical lax transformation $\theta_0^A: (-,A)_1\Rightarrow(-,A)_2$ and 
for each $B\in\Bb$ a vertical lax transformation $\theta_0^B: (B,-)_1\Rightarrow(B,-)_2$, both of lax double functors, such that 
$(\theta_0^A)_B=(\theta_0^B)_A$ and such that 

\medskip
\noindent $(VLT^q_1)$ \vspace{-0,6cm}

$$
\scalebox{0.8}{
\bfig
 \putmorphism(-150,500)(1,0)[(B,A)_1`(B,A)_1 `=]{600}1a
\putmorphism(-180,500)(0,-1)[\phantom{Y_2}`(B, A)_2 `(\theta_0^A)_B]{450}1l
\put(-10,250){\fbox{$(\theta_0^A)^u$}}
\putmorphism(-150,-400)(1,0)[(\tilde B, A)_2`(\tilde B, A)_2 `=]{640}1a
\putmorphism(-180,50)(0,-1)[\phantom{Y_2}``(u,A)_2]{450}1l
\putmorphism(450,50)(0,-1)[\phantom{Y_2}``(\theta_0^A)_{\tilde B}]{450}1l
\putmorphism(450,500)(0,-1)[\phantom{Y_2}`(\tilde B, A)_1 `(u,A)_1]{450}1r
\putmorphism(-820,50)(1,0)[(B, A)_2``=]{520}1a
\putmorphism(-820,50)(0,-1)[\phantom{(B, \tilde A')}``(B,U)_2]{450}1l
\putmorphism(-820,-400)(0,-1)[(B, \tilde A)_2`(\tilde B, \tilde A)_2`(u,\tilde A)_2]{450}1l
\putmorphism(-820,-850)(1,0)[\phantom{(B, \tilde A)}``=]{520}1a
\putmorphism(-180,-400)(0,-1)[`(\tilde B, \tilde A)_2 `(\tilde B, U)_2]{450}1r
\put(-650,-630){\fbox{$(u,U)_2$}}

\putmorphism(570,50)(1,0)[`(\tilde B, A)_1`=]{520}1a
\putmorphism(1050,50)(0,-1)[\phantom{Y_2}`(\tilde B, \tilde A)_1`(\tilde B, U)_1]{450}1r
\putmorphism(450,-400)(0,-1)[`(\tilde B, \tilde A)_2 `(\tilde B,U)_2]{450}1r
\putmorphism(1050,-400)(0,-1)[`(\tilde B, \tilde A)_2 `(\theta_0^{\tilde B})_{\tilde A}]{450}1r
\putmorphism(570,-850)(1,0)[` `=]{340}1a
\put(600,-200){\fbox{$(\theta_0^{\tilde B})^U$}}
\efig}
=
\scalebox{0.8}{
\bfig
 \putmorphism(-150,500)(1,0)[(B,A)_1`(B,A)_1 `=]{600}1a
\putmorphism(-180,500)(0,-1)[\phantom{Y_2}`(B, A)_2 `(\theta^A_0)_B]{450}1l
\put(0,50){\fbox{$(\theta^B_0)^U$}}
\putmorphism(-180,-400)(1,0)[(B, \tilde A)_2` `=]{500}1a
\putmorphism(-180,50)(0,-1)[\phantom{Y_2}``(B,U)_2]{450}1l
\putmorphism(450,50)(0,-1)[\phantom{Y_2}`(B, \tilde A)_2`(\theta^B_0)_{\tilde A}]{450}1r
\putmorphism(450,500)(0,-1)[\phantom{Y_2}`(B, \tilde A)_1 `(B,U)_1]{450}1r
\putmorphism(450,50)(1,0)[\phantom{(B, \tilde A)}`(B, \tilde A)_1`=]{620}1a
\putmorphism(1070,50)(0,-1)[\phantom{(B, \tilde A')}``(u,\tilde A)_1]{450}1r
\putmorphism(1070,-400)(0,-1)[(\tilde B, \tilde A)_1`(\tilde B, \tilde A)_2`(\theta^{\tilde A}_0)_{\tilde B}]{450}1r
\putmorphism(450,-850)(1,0)[\phantom{(B, \tilde A)}``=]{500}1a
\putmorphism(450,-400)(0,-1)[\phantom{(B, \tilde A)}`(\tilde B, \tilde A) `(u,\tilde A)_2]{450}1l
\put(600,-630){\fbox{$(\theta^{\tilde A}_0)^u$}}

 \putmorphism(1070,500)(1,0)[(B,A)_1`(B,A)_1 `=]{670}1a
\putmorphism(1070,500)(0,-1)[\phantom{Y_2}` `(B,U)_1]{450}1r
\putmorphism(1700,500)(0,-1)[\phantom{Y_2}`(\tilde B, A)_1 `(u,A)_1]{450}1r
\putmorphism(1700,50)(0,-1)[\phantom{Y_2}`(\tilde B, \tilde A)_1 `(\tilde B,U)_1]{450}1r
\putmorphism(1190,-400)(1,0)[` `=]{370}1a
\put(1200,50){\fbox{$(u,U)_1$}}
\efig}
$$
for every 1v-cells $U:A\to \tilde A$ and $u:B\to\tilde B$; 

\medskip
\noindent $(VLT^q_2)$ \vspace{-0,6cm}
$$\scalebox{0.86}{
\bfig
 \putmorphism(-150,500)(1,0)[(B,A)_1`(B,A)_1 `=]{600}1a
 \putmorphism(450,500)(1,0)[\phantom{(B,A)_1}` `(B,K)_1]{450}1a
\putmorphism(-200,500)(0,-1)[\phantom{Y_2}`(B, A)_2 `(\theta_0^A)_B]{450}1l
\put(-50,50){\fbox{$(\theta^A_0)^u$}}
\putmorphism(-170,-400)(1,0)[(\tilde B, A)_2` `=]{500}1a
\putmorphism(-200,50)(0,-1)[\phantom{Y_2}``(u,A)_2]{450}1l
\putmorphism(450,50)(0,-1)[\phantom{Y_2}`(\tilde B, A)_2 `(\theta_0^A)_{\tilde B}]{450}1l
\putmorphism(450,500)(0,-1)[\phantom{Y_2}`(\tilde B, A)_1 `(u,A)_1]{450}1l
\put(600,260){\fbox{$(u,K)_1$}}
\putmorphism(450,50)(1,0)[\phantom{(B, \tilde A)}``(\tilde B,K)_1]{500}1a
\putmorphism(1070,50)(0,-1)[\phantom{(B, A')}`(\tilde B, A')_2`(\theta^{\tilde B}_0)_{A'}]{450}1r
\putmorphism(1070,500)(0,-1)[(B, A')_1`(\tilde B, A')_1`(u,A')_1]{450}1r
\putmorphism(450,-400)(1,0)[\phantom{(B, \tilde A)}``(\tilde B,K)_2]{500}1a
\put(600,-170){\fbox{$(\theta_0^{\tilde B})_K$}}
\efig}=
\scalebox{0.86}{
\bfig
 \putmorphism(-150,500)(1,0)[(B,A)_1`(B,A')_1 `(B,K)_1]{600}1a
 \putmorphism(450,500)(1,0)[\phantom{(B,A)}` `=]{540}1a
\putmorphism(-180,500)(0,-1)[\phantom{Y_2}`(B, A)_2 `(\theta^A_0)_B]{450}1l
\put(0,280){\fbox{$(\theta_0^B)_K$}}
\putmorphism(-180,-400)(1,0)[(\tilde B,A)_2` `(\tilde B,K)_2]{500}1a
\putmorphism(-180,50)(0,-1)[\phantom{Y_2}``(u,A)_2]{450}1l
\putmorphism(450,50)(0,-1)[\phantom{Y_2}`(\tilde B, A')_2`(u, A')_2]{450}1r
\putmorphism(450,500)(0,-1)[\phantom{Y_2}`(B, A')_2 `(\theta^{A'}_0)_B]{450}1r
\put(600,50){\fbox{$(\theta_0^{A'})^u$}}
\putmorphism(-180,50)(1,0)[\phantom{(B, \tilde A)}``(B,K)_2]{500}1a
\putmorphism(1130,50)(0,-1)[\phantom{(B, A')}`(\tilde B, A')_2`(\theta_0^{A'})_{\tilde B}]{450}1r
\putmorphism(1130,500)(0,-1)[(B, A')_1`(\tilde B, A')_1`(u,A')_1]{450}1r
\putmorphism(450,-400)(1,0)[\phantom{(B, \tilde A)}``=]{520}1b
\put(0,-170){\fbox{$(u,K)_2$}}
\efig}
$$
for every 1h-cell $K:A\to A'$ and 1v-cell $u: B\to \tilde B$, 

\medskip
\noindent $(VLT^q_3)$ \vspace{-0,6cm}

$$\scalebox{0.86}{
\bfig
 \putmorphism(-150,500)(1,0)[(B,A)_1`(B,A)_1 `=]{600}1a
 \putmorphism(450,500)(1,0)[\phantom{(B',A)_1}` `(k,A)_1]{450}1a
\putmorphism(-200,500)(0,-1)[\phantom{Y_2}`(B, A)_2 `(\theta_0^B)_A]{450}1l
\put(-40,50){\fbox{$(\theta^B_0)^U$}}
\putmorphism(-170,-400)(1,0)[(B, \tilde A)_2` `=]{500}1a
\putmorphism(-200,50)(0,-1)[\phantom{Y_2}``(B,U)_2]{450}1l
\putmorphism(450,50)(0,-1)[\phantom{Y_2}`(B, \tilde A)_2 `(\theta_0^B)_{\tilde A}]{450}1l 
\putmorphism(450,500)(0,-1)[\phantom{Y_2}`(B, \tilde A)_1 `(B,U)_1]{450}1l
\put(600,260){\fbox{$(k,U)_1$}}
\putmorphism(450,50)(1,0)[\phantom{(B, \tilde A)}``(k,\tilde A)_1]{500}1a %
\putmorphism(1070,50)(0,-1)[\phantom{(B, A')}`(B', \tilde A)_2`(\theta^{\tilde A}_0)_{B'}]{450}1r
\putmorphism(1070,500)(0,-1)[(B, A')_1`(B', \tilde A)_1`(B',U)_1]{450}1r
\putmorphism(450,-400)(1,0)[\phantom{(B, \tilde A)}``(k,\tilde A)_2]{500}1a
\put(600,-170){\fbox{$(\theta_0^{\tilde A})_k$}}
\efig}
\quad=\quad
\scalebox{0.86}{
\bfig
 \putmorphism(-150,500)(1,0)[(B,A)_1`(B',A)_1 `(k,A)_1]{600}1a
 \putmorphism(450,500)(1,0)[\phantom{(B,A)}` `=]{540}1a
\putmorphism(-180,500)(0,-1)[\phantom{Y_2}`(B, A)_2 `(\theta^A_0)_B]{450}1l
\put(0,280){\fbox{$(\theta_0^A)_k$}}
\putmorphism(-180,-400)(1,0)[(\tilde B,A)_2` `(k,\tilde A)_2]{500}1a
\putmorphism(-180,50)(0,-1)[\phantom{Y_2}``(B,U)_2]{450}1l
\putmorphism(450,50)(0,-1)[\phantom{Y_2}`(B', \tilde A)_2`(B', U)_2]{450}1r
\putmorphism(450,500)(0,-1)[\phantom{Y_2}`(B', A)_2 `(\theta^{B'}_0)_A]{450}1r
\put(630,50){\fbox{$(\theta_0^{B'})^U$}}
\putmorphism(-180,50)(1,0)[\phantom{(B, \tilde A)}``(k,A)_2]{500}1a
\putmorphism(1130,50)(0,-1)[\phantom{(B, A')}`(B',\tilde A)_2`(\theta_0^{\tilde A})_{B'}]{450}1r
\putmorphism(1130,500)(0,-1)[(B', A)_1`(B',\tilde A)_1`(B', U)_1]{450}1r
\putmorphism(450,-400)(1,0)[\phantom{(B, \tilde A)}``=]{520}1b
\put(0,-170){\fbox{$(k,U)_2$}}
\efig}
$$
for every 1v-cell $U:A\to \tilde A$ and 1h-cell $k:B\to B'$, and 

\medskip
\noindent $(VLT^q_4)$ \vspace{-0,6cm}
$$\scalebox{0.86}{
\bfig
\putmorphism(-150,500)(1,0)[(B,A)_1`(B',A)_1`(k,A)_1]{600}1a
 \putmorphism(480,500)(1,0)[\phantom{F(A)}`(B',A')_1 `(B',K)_1]{640}1a
 \putmorphism(-150,50)(1,0)[(B,A)_1`(B,A')_1`(B,K)_1]{600}1a
 \putmorphism(470,50)(1,0)[\phantom{F(A)}`(B',A')_1 `(k,A')_1]{660}1a

\putmorphism(-180,500)(0,-1)[\phantom{Y_2}``=]{450}1r
\putmorphism(1100,500)(0,-1)[\phantom{Y_2}``=]{450}1r
\put(320,280){\fbox{$(k,K)_1$}}

\putmorphism(-150,-400)(1,0)[(B,A)_2`(B,A')_2 `(B,K)_2]{640}1a
 \putmorphism(490,-400)(1,0)[\phantom{F(A')}` (B',A')_2 `(k,A')_2]{640}1a

\putmorphism(-180,50)(0,-1)[\phantom{Y_2}``(\theta_0^A)_B]{450}1l %
\putmorphism(450,50)(0,-1)[\phantom{Y_2}``]{450}1l
\putmorphism(610,50)(0,-1)[\phantom{Y_2}``(\theta_0^B)_{A'}]{450}0l 
\putmorphism(1120,50)(0,-1)[\phantom{Y_3}``(\theta_0^{A'})_{B'}]{450}1r
\put(-40,-180){\fbox{$(\theta^B_0)_K$}} 
\put(620,-180){\fbox{$(\theta_0^{A'})_k$}}
\efig}
\quad
=
\quad
\scalebox{0.86}{
\bfig
\putmorphism(-150,500)(1,0)[(B,A)_1`(B',A)_1`(k,A)_1]{600}1a
 \putmorphism(480,500)(1,0)[\phantom{F(A)}`(B',A')_1 `(B',K)_1]{640}1a

 \putmorphism(-150,50)(1,0)[(B,A)_2`(B',A)_2`(k,A)_2]{600}1a
 \putmorphism(450,50)(1,0)[\phantom{F(A)}`(B',K)_2 `(B', K)_2]{640}1a

\putmorphism(-180,500)(0,-1)[\phantom{Y_2}``(\theta_0^A)_B]{450}1l
\putmorphism(450,500)(0,-1)[\phantom{Y_2}``]{450}1r
\putmorphism(300,500)(0,-1)[\phantom{Y_2}``(\theta_0^{B'})_A]{450}0r
\putmorphism(1100,500)(0,-1)[\phantom{Y_2}``(\theta_0^{B'})_{A'}]{450}1r
\put(-40,280){\fbox{$(\theta_0^{A})_k$}}
\put(700,280){\fbox{$(\theta_0^{B'})_K$}}

\putmorphism(-150,-400)(1,0)[(B,A)_2`(B,A')_2 `(B,K)_2]{640}1a
 \putmorphism(490,-400)(1,0)[\phantom{F(A')}` (B',A')_2 `(k,A')_2]{640}1a

\putmorphism(-180,50)(0,-1)[\phantom{Y_2}``=]{450}1l
\putmorphism(1120,50)(0,-1)[\phantom{Y_3}``=]{450}1r
\put(320,-200){\fbox{$(k,K)_2$}}
\efig}
$$
for every 1h-cells $K:A\to A'$ and $k:B\to B'$. 
\end{defn}

\begin{defn} \delabel{modif btw horiz}
Let horizontal oplax transformations $\theta, \theta'$ and vertical lax transformations $\theta_0, \theta'_0$ acting between lax double quasi-functors 
$H_1, H_2, H_3, H_4: \Aa\times\Bb\to\Cc$ be given as in the left diagram below. Denote by $(-,A)_i:\Bb\to\Cc, (B,-)_i:\Aa\to\Cc, i=1,2,3,4$ the pairs of 
lax double functors corresponding to $H_1, H_2, H_3, H_4$, respectively. A modification $\Theta$ (on the left below) is given by a pair of modifications 
$\tau^A, \tau^B$ acting between transformations among lax double functors:
\begin{equation} \eqlabel{q-modif}
\scalebox{0.86}{
\bfig
\putmorphism(-150,50)(1,0)[H_1` H_2` \theta]{400}1a
\putmorphism(-150,-270)(1,0)[H_3 ` H_4 ` \theta' ]{400}1b
\putmorphism(-170,50)(0,-1)[\phantom{Y_2}``\theta_0]{320}1l
\putmorphism(250,50)(0,-1)[\phantom{Y_2}``\theta_0']{320}1r
\put(-30,-140){\fbox{$\tau$}}
\efig}
\qquad\qquad
\scalebox{0.86}{
\bfig
\putmorphism(-180,50)(1,0)[(-,A)_1` (-,A)_2`\theta^A]{550}1a
\putmorphism(-180,-270)(1,0)[(-,A)_3`(-,A)_4 `\theta'^A]{550}1b
\putmorphism(-170,50)(0,-1)[\phantom{Y_2}``\theta_0^A]{320}1l
\putmorphism(350,50)(0,-1)[\phantom{Y_2}``\theta_0'^A]{320}1r
\put(0,-140){\fbox{$\tau^A$}}
\efig}
\qquad\qquad
\scalebox{0.86}{
\bfig
\putmorphism(-180,50)(1,0)[(B,-)_1` (B,-)_2`\theta^B]{550}1a
\putmorphism(-180,-270)(1,0)[(B,-)_3`(B,-)_4 `\theta'^B]{550}1b
\putmorphism(-170,50)(0,-1)[\phantom{Y_2}``\theta_0^B]{320}1l
\putmorphism(350,50)(0,-1)[\phantom{Y_2}``\theta_0'^B]{320}1r
\put(0,-140){\fbox{$\tau^B$}}
\efig}
\end{equation}
such that $\tau^A_B=\tau^B_A$ for every $A\in\Aa, B\in\Bb$. 
\end{defn}

The composition of 1- and 2-cells in $q\x\Lax_{hop}(\Aa\times\Bb,\Cc)$ is given in the analogous way as in 
$\Lax_{hop}(\Aa, \llbracket\Bb,\Cc\rrbracket)$.

\subsection{The 1-1 correspondence between 1h- and 1v-cells}  \sslabel{1-cells}

We proceed to show that the double categories $\Lax_{hop}(\Aa, \llbracket\Bb,\Cc\rrbracket)$ and $q\x\Lax_{hop}(\Aa\times\Bb,\Cc)$ are isomorphic. 
From \prref{char df} and \deref{H dbl} we know that we have a 1-1 correspondence between their corresponding 0-cells. 

Let $F,G: \Aa\to\llbracket\Bb,\Cc\rrbracket$ be two lax double functors and take a horizontal oplax transformation $\alpha:F\Rightarrow G$. 
Set $(-,-)_1$ and $(-,-)_2$ for the two lax double quasi-functors obtained from $F$ and $G$, respectively. 
Evaluating at a 0-cell $A\in\Aa$ we get $\alpha(A): F(A)\to G(A)$ a 1h-cell in $\llbracket\Bb,\Cc\rrbracket$ of the form 
$(-,A)_1\to(-,A)_2$. This 1h-cell is a horizontal oplax transformation between lax double functors, so we have 
the following cells in $\Cc$:
a 1h-cell $\alpha(A)_B: (B,A)_1\to(B,A)_2$, a globular 2-cell $\alpha(A)_k=\delta_{\alpha(A),k}$, and 
a 2-cell $\alpha(A)_u$, for a 0-cell $B$, a 1h-cell $k:B\to B'$, and a 1v-cell $u:B\to\tilde B$ in $\Bb$. 
The 2-cells $\alpha(A)_k$ and $\alpha(A)_u$ in $\Cc$ satisfy the five axioms from \deref{hor nat tr}.

On the other hand, evaluating the horizontal oplax transformation $\alpha:F\Rightarrow G$ at a 1h-cell $K:A\to A'$ in $\Aa$, one obtains a 
globular 2-cell $\alpha_K=\delta_{\alpha,K}: \frac{(-,K)_1}{\alpha(A')}\Rrightarrow \frac{\alpha(A)}{(-,K)_2}$ in 
$\llbracket\Bb,\Cc\rrbracket$, which, by the horizontal restriction of \deref{modif-hv}, is a modification between (the vertical composition of) 
horizontal oplax transformations of lax double functors. 
It has a free slot for 0-cells in $\Bb$, so that after evaluation at some $B\in\Bb$ it yields a 
globular 2-cell $\alpha_K(B)$ in $\Cc$. Finally, evaluating $\alpha$ at a 1v-cell $U:A\to\tilde A$ in $\Aa$, one obtains a 2-cell $\alpha_U$:
\begin{equation}\eqlabel{alfa-U}
\scalebox{0.86}{
\bfig
 \putmorphism(-200,50)(1,0)[(-,A)_1`(-,A)_2`\alpha(A)]{600}1a
\putmorphism(-200,-350)(1,0)[(-,\tilde A)_1`(-,\tilde A)_2 `\alpha(\tilde A)]{600}1a
\putmorphism(-180,50)(0,-1)[\phantom{Y_2}``(-,U)_1]{400}1l
\putmorphism(370,50)(0,-1)[\phantom{Y_2}``(-,U)_2]{400}1r
\put(0,-110){\fbox{$\alpha_U$}} 
\efig}
\end{equation}
in $\llbracket\Bb,\Cc\rrbracket$ (thus a modification in the sense of \deref{modif-hv}) 
with a free slot for 0-cells in $\Bb$ (after evaluation at $B\in\Bb$ it yields a 2-cell $\alpha_U(B)$ in $\Cc$). 

The families of 2-cells $\alpha_K$ and $\alpha_U$ in $\llbracket\Bb,\Cc\rrbracket$ from the horizontal oplax transformation $\alpha:F\Rightarrow G$ 
satisfy the five axioms from \deref{hor nat tr}. Evaluating these five axioms at $B\in\Bb$ one obtains five axioms 
for families of 2-cell $\alpha_K(B)$ and $\alpha_U(B)$ in $\Cc$. The latter axioms mean that $\alpha(-)_B: F(-)(B)\Rightarrow G(-)(B)$, 
obtained by reading $\alpha:F\Rightarrow G$ described above the other way around, that is, evaluating at a 0-cell $B\in\Bb$ and leaving a free slot for cells from $\Aa$, is a horizontal oplax transformation between lax double functors $(B,-)_1\to(B,-)_2$ which act between $\Aa\to\Cc$.
%
Namely, set $\alpha(K)_B=\delta_{\alpha(-)_B,K}:=\alpha_K(B)$ and $\alpha(U)_B:=\alpha_U(B)$, for a 1h-cell $K:A\to A'$ and 1v-cell 
$U:A\to\tilde A$ in $\Aa$. 

Now, we may set $\theta^A:=\alpha(A)$ and $\theta^B:=\alpha(-)_B$ for two horizontal oplax transformations between lax double functors. 
We do have that $\theta^A_B=\theta^B_A$, it remains to check the other four conditions in order for the pairs $(\theta^A, \theta^B)$ for $A\in\Aa, B\in\Bb$ 
to make a horizontal oplax transformation $\theta: (-,-)_1\Rightarrow(-,-)_2$ between lax double quasi-functors.

\begin{prop} \prlabel{equiv-modif}
Let $F,G: \Aa\to\llbracket\Bb,\Cc\rrbracket$ be two lax double functors with the corresponding lax double quasi-functors $(-,-)_1, (-,-)_2:\Aa\times\Bb\to\Cc$. 
For every $A\in\Aa$ and $B\in\Bb$ let $\alpha(A): F(A)\to G(A)$ and $\alpha(-)_B: F(-)(B)\to G(-)(B)$ be horizontal oplax transformations between lax double functors. 
The following are equivalent: 
\begin{enumerate}
\item $\alpha_K: \frac{(-,K)_1}{\alpha(A')}\Rrightarrow \frac{\alpha(A)}{(-,K)_2}$ is a modification 
on the vertical composition of horizontal oplax transformations of lax double functors (with components $(\alpha_K)_B=\delta_{\alpha(-)_B,K}
: \frac{(B,K)_1}{\alpha(A')_B}\Rrightarrow \frac{\alpha(A)_B}{(B,K)_2}$, recall \equref{m-hor}) for every 1h-cell $K: A\to A'$ in $\Aa$, and 
$\alpha_U$ of the form \equref{alfa-U} is a modification in the sense of \deref{modif-hv} for every 1v-cell $U:A\to\tilde A$ in $\Aa$; 
\item the pairs $(\theta^A, \theta^B):=(\alpha(A),\alpha(-)_B)$ for $A\in\Aa, B\in\Bb$ form a horizontal oplax transformation $\theta: (-,-)_1\Rightarrow(-,-)_2$ between lax double quasi-functors. 
\end{enumerate}
\end{prop}

\begin{proof}
In \coref{obtained oplax trans} we saw that $(-,K)_i, i=1,2$ are horizontal oplax transformations. 
From \leref{vert comp hor.ps.tr.} we have that the composite transformations $[\frac{(-,K)_1}{\alpha(A')}]_k$ and $[\frac{\alpha(A)}{(-,K)_2}]_k$ 
evaluated at a 1h-cell $k:B\to B'$ have the following form:
$$\delta_{\frac{(-,K)_1}{\alpha(A')},k}=
\scalebox{0.86}{\bfig

 \putmorphism(-200,-50)(1,0)[(B,A)_1`(B',A)_1` (k,A)_1]{650}1a
 \putmorphism(480,-50)(1,0)[\phantom{F(B)}`(B',A')_1 `(B',K)_1]{700}1a

 \putmorphism(-200,-450)(1,0)[(B,A)_1`(B,A')_1`(B,K)_1]{650}1a
 \putmorphism(460,-450)(1,0)[\phantom{A\ot B}`(B',A')_1 `(k,A')_1]{700}1a
 \putmorphism(1150,-450)(1,0)[\phantom{A'\ot B'}` (B',A')_2 `\alpha(A')_{B'}]{700}1a

\putmorphism(-230,-50)(0,-1)[\phantom{Y_2}``=]{400}1r
\putmorphism(1100,-50)(0,-1)[\phantom{Y_2}``=]{400}1r
\put(300,-240){\fbox{$ \delta_{(-,K)_1,k}  $}}
\put(1000,-660){\fbox{$\delta_{\alpha(A'),k}$}}

 \putmorphism(450,-850)(1,0)[(B,A')_1` (B,A')_2 `\alpha(A')_B]{700}1a
 \putmorphism(1130,-850)(1,0)[\phantom{A''\ot B'}`(B',A')_2 ` (k,A')_2]{700}1a

\putmorphism(450,-450)(0,-1)[\phantom{Y_2}``=]{400}1l
\putmorphism(1800,-450)(0,-1)[\phantom{Y_2}``=]{400}1r
\efig}
$$
and
$$\delta_{\frac{\alpha(A)}{(-,K)_2},k}=
\scalebox{0.86}{\bfig

 \putmorphism(-200,-50)(1,0)[(B,A)_1`(B',A)_1` (k,A)_1]{650}1a
 \putmorphism(480,-50)(1,0)[\phantom{F(B)}`(B',A)_2 `\alpha(A)_{B'}]{700}1a

 \putmorphism(-200,-450)(1,0)[(B,A)_1`(B,A)_2`\alpha(A)_B]{650}1a
 \putmorphism(450,-450)(1,0)[\phantom{A\ot B}`(B',A)_2 `(k,A)_2]{700}1a
 \putmorphism(1140,-450)(1,0)[\phantom{A'\ot B'}` (B',A')_2 `(B',K)_2]{700}1a

\putmorphism(-230,-50)(0,-1)[\phantom{Y_2}``=]{400}1r
\putmorphism(1050,-50)(0,-1)[\phantom{Y_2}``=]{400}1r
\put(300,-240){\fbox{$ \delta_{\alpha(A),k}  $}}
\put(1000,-660){\fbox{$\delta_{(-,K)_2,k}$}}

 \putmorphism(450,-850)(1,0)[(B,A)_2` (B,A')_2 `(B,K)_2]{700}1a
 \putmorphism(1130,-850)(1,0)[\phantom{A''\ot B'}`(B',A')_2 ` (k,A')_2]{700}1a

\putmorphism(450,-450)(0,-1)[\phantom{Y_2}``=]{400}1l
\putmorphism(1780,-450)(0,-1)[\phantom{Y_2}``=]{400}1r
\efig}
$$
Now the first modification condition (m.ho.-1) for $\alpha_K$ reads: 
$$\scalebox{0.8}{
\bfig
 \putmorphism(-200,600)(1,0)[(B,A)_1`(B',A)_1` (k,A)_1]{650}1a
 \putmorphism(480,600)(1,0)[\phantom{F(B)}`(B',A')_1 `(B',K)_1]{700}1a
\putmorphism(-230,600)(0,-1)[\phantom{Y_2}``=]{400}1r
\putmorphism(1100,600)(0,-1)[\phantom{Y_2}``=]{400}1r

 \putmorphism(-200,200)(1,0)[(B,A)_1`(B,A')_1`(B,K)_1]{650}1a
 \putmorphism(460,200)(1,0)[\phantom{A\ot B}`(B',A')_1 `(k,A')_1]{700}1a
 \putmorphism(1150,200)(1,0)[\phantom{A'\ot B'}` (B',A')_2 `\alpha(A')_{B'}]{700}1a
\putmorphism(450,200)(0,-1)[\phantom{Y_2}``=]{400}1l
\putmorphism(1800,200)(0,-1)[\phantom{Y_2}``=]{400}1r

\put(300,400){\fbox{$ \delta_{(-,K)_1,k}  $}}
\put(1000,-10){\fbox{$\delta_{\alpha(A'),k}$}}

\putmorphism(-200,-200)(1,0)[(B,A)_1`\phantom{(B,A')_1}`(B,K)_1]{650}1a
\putmorphism(450,-200)(1,0)[(B,A')_1` (B,A')_2 `\alpha(A')_B]{700}1a
 \putmorphism(1130,-200)(1,0)[\phantom{A''\ot B'}`(B',A')_2 ` (k,A')_2]{700}1a

\putmorphism(-230,-200)(0,-1)[\phantom{Y_2}``=]{400}1r
\putmorphism(1100,-200)(0,-1)[\phantom{Y_2}``=]{400}1r
\putmorphism(-200,-600)(1,0)[(B,A)_1`\phantom{(B,A')_1}`\alpha(A)_B]{650}1a
\putmorphism(450,-600)(1,0)[(B,A')_1` (B,A')_2 `(B,K)_2]{700}1a
\put(260,-430){\fbox{$ \delta_{\alpha(-)_B,K}  $}}
\efig}
=
\scalebox{0.8}{
\bfig
\putmorphism(450,600)(0,-1)[\phantom{Y_2}``=]{400}1l
\putmorphism(1780,600)(0,-1)[\phantom{Y_2}``=]{400}1r
 \putmorphism(460,600)(1,0)[(B',A)_1`(B',A)_1 `(B',K)_1]{700}1a
 \putmorphism(1140,600)(1,0)[\phantom{A'\ot B'}` (B',A')_2 `\alpha(A')_{B'}]{700}1a
\put(1000,410){\fbox{$\delta_{\alpha(-)_{B'},K}$}}

 \putmorphism(-200,200)(1,0)[(B,A)_1`(B',A)_1` (k,A)_1]{650}1a
 \putmorphism(480,200)(1,0)[\phantom{F(B)}`(B',A)_2 `\alpha(A)_{B'}]{670}1a
 \putmorphism(1140,200)(1,0)[\phantom{A'\ot B'}` (B',A')_2 `(B',K)_2]{700}1a
\putmorphism(-230,200)(0,-1)[\phantom{Y_2}``=]{400}1r
\putmorphism(1050,200)(0,-1)[\phantom{Y_2}``=]{400}1r
\put(300,10){\fbox{$ \delta_{\alpha(A),k}  $}}

 \putmorphism(-200,-200)(1,0)[(B,A)_1`(B,A)_2`\alpha(A)_B]{650}1a
 \putmorphism(450,-200)(1,0)[\phantom{A\ot B}`(B',A)_2 `(k,A)_2]{700}1a
 \putmorphism(1140,-200)(1,0)[\phantom{A'\ot B'}` (B',A')_2 `(B',K)_2]{700}1a
\putmorphism(450,-200)(0,-1)[\phantom{Y_2}``=]{400}1l
\putmorphism(1780,-200)(0,-1)[\phantom{Y_2}``=]{400}1r
\put(1000,-410){\fbox{$\delta_{(-,K)_2,k}$}}

 \putmorphism(450,-600)(1,0)[(B,A)_2` (B,A')_2 `(B,K)_2]{700}1a
 \putmorphism(1130,-600)(1,0)[\phantom{A''\ot B'}`(B',A')_2 ` (k,A')_2]{700}1a
\efig}
$$
Recall from \coref{obtained oplax trans} that $\delta_{(-,K)_i,k}=(k,K)_i$ for $i=1,2$ and that we are setting 
$\theta^A:=\alpha(A)$ and $\theta^B:=\alpha(-)_B$, thus $\delta_{\alpha(A),k}=\theta^A_k$ and $\delta_{\alpha(-)_B,K}=\theta^B_K$. 
We have that the above modification condition is precisely $(HOT^q_1)$. 

By \leref{vert comp hor.ps.tr.} we have: 
$$\frac{(-,K)_1}{\alpha(A')}\vert_u=
\scalebox{0.86}{
\bfig
\putmorphism(-150,250)(1,0)[(B,A)_1`(B,A')_1`(B,K)_1]{600}1a
 \putmorphism(450,250)(1,0)[\phantom{F(A)}`(B,A')_2 `\alpha(A')_B]{640}1a

 \putmorphism(-150,-200)(1,0)[(\tilde B,A)_1`(\tilde B,A')_1`(\tilde B,K)_1]{600}1a
 \putmorphism(450,-200)(1,0)[\phantom{F(A)}`(\tilde B,A')_2 `\alpha(A')_{\tilde B}]{640}1a

\putmorphism(-180,250)(0,-1)[\phantom{Y_2}``(u,A)_1]{450}1l
\putmorphism(450,250)(0,-1)[\phantom{Y_2}``]{450}1r
\putmorphism(300,250)(0,-1)[\phantom{Y_2}``(u,A')_1]{450}0r
\putmorphism(1100,250)(0,-1)[\phantom{Y_2}``(u,A')_2]{450}1r
\put(-40,30){\fbox{$(u,K)_1$}}
\put(700,30){\fbox{$\alpha(A')_u$}}
\efig}
$$
and 
$$\frac{\alpha(A)}{(-,K)_2}\vert_u=
\scalebox{0.86}{
\bfig
 \putmorphism(-150,300)(1,0)[(B,A)_1`(B,A)_2`\alpha(A)_B]{600}1a
 \putmorphism(450,300)(1,0)[\phantom{F(A)}`(B,A')_2 `(B,K)_2]{680}1a

\putmorphism(-150,-150)(1,0)[(\tilde B,A)_1`(\tilde B,A)_2 `\alpha(A)_{\tilde B}]{640}1a
 \putmorphism(490,-150)(1,0)[\phantom{F(A')}` (\tilde B,A')_2 `(\tilde B,K)_2]{640}1a

\putmorphism(-180,300)(0,-1)[\phantom{Y_2}``(u,A)_1]{450}1l
\putmorphism(450,300)(0,-1)[\phantom{Y_2}``]{450}1l
\putmorphism(610,300)(0,-1)[\phantom{Y_2}``(u,A)_2]{450}0l 
\putmorphism(1120,300)(0,-1)[\phantom{Y_3}``(u,A')_2]{450}1r
\put(-40,70){\fbox{$\alpha(A)_u$}} 
\put(620,70){\fbox{$(u,K)_2$}}
\efig}
$$
for a 1v-cell $u:B\to \tilde B$. Now the second modification condition (m.ho.-2) for $\alpha_K$ reads:
$$
\scalebox{0.86}{
\bfig
\putmorphism(-150,500)(1,0)[(B,A)_1`(B,A')_1`(B,K)_1]{600}1a
 \putmorphism(450,500)(1,0)[\phantom{F(A)}`(B,A')_2 `\alpha(A')_B]{640}1a

 \putmorphism(-150,50)(1,0)[(\tilde B,A)_1`(\tilde B,A')_1`(\tilde B,K)_1]{600}1a
 \putmorphism(450,50)(1,0)[\phantom{F(A)}`(\tilde B,A')_2 `\alpha(A')_{\tilde B}]{640}1a

\putmorphism(-180,500)(0,-1)[\phantom{Y_2}``(u,A)_1]{450}1l
\putmorphism(450,500)(0,-1)[\phantom{Y_2}``]{450}1r
\putmorphism(300,500)(0,-1)[\phantom{Y_2}``(u,A')_1]{450}0r
\putmorphism(1100,500)(0,-1)[\phantom{Y_2}``(u,A')_2]{450}1r
\put(-40,280){\fbox{$(u,K)_1$}}
\put(700,280){\fbox{$\alpha(A')_u$}}

\putmorphism(-150,-400)(1,0)[F(A')`G(A') `\alpha(A)_{\tilde B}]{640}1a
 \putmorphism(450,-400)(1,0)[\phantom{A'\ot B'}` G(B') `(\tilde B,K)_2]{680}1a

\putmorphism(-180,50)(0,-1)[\phantom{Y_2}``=]{450}1l
\putmorphism(1120,50)(0,-1)[\phantom{Y_3}``=]{450}1r
\put(320,-200){\fbox{$\delta_{\alpha(-)_{\tilde B},K}$}}
\efig}
\quad=\quad
\scalebox{0.86}{
\bfig
\putmorphism(-150,500)(1,0)[(B,A)_1`(B,A')_1`(B,K)_1]{600}1a
 \putmorphism(450,500)(1,0)[\phantom{F(A)}`(B,A')_2 `\alpha(A')_B]{640}1a
 \putmorphism(-150,50)(1,0)[(B,A)_1`(B,A)_2`\alpha(A)_B]{600}1a
 \putmorphism(450,50)(1,0)[\phantom{F(A)}`(B,A')_2 `(B,K)_2]{680}1a

\putmorphism(-180,500)(0,-1)[\phantom{Y_2}``=]{450}1r
\putmorphism(1100,500)(0,-1)[\phantom{Y_2}``=]{450}1r
\put(320,280){\fbox{$\delta_{\alpha(-)_B,K}$}}

\putmorphism(-150,-400)(1,0)[(\tilde B,A)_1`(\tilde B,A)_2 `\alpha(A)_{\tilde B}]{640}1a
 \putmorphism(490,-400)(1,0)[\phantom{F(A')}` (\tilde B,A')_2 `(\tilde B,K)_2]{640}1a

\putmorphism(-180,50)(0,-1)[\phantom{Y_2}``(u,A)_1]{450}1l
\putmorphism(450,50)(0,-1)[\phantom{Y_2}``]{450}1l
\putmorphism(610,50)(0,-1)[\phantom{Y_2}``(u,A)_2]{450}0l 
\putmorphism(1120,50)(0,-1)[\phantom{Y_3}``(u,A')_2]{450}1r
\put(-40,-180){\fbox{$\alpha(A)_u$}} 
\put(620,-180){\fbox{$(u,K)_2$}}
\efig}
$$
Setting $\theta^A_u=\alpha(A)_u$ and $\theta^B_K=\delta_{\alpha(-)_B,K}$ this is precisely $(HOT^q_2)$. 

The two modification conditions $(m.ho-vl.-1)$ and $(m.ho-vl.-2)$ for $\alpha_U$ are: 
$$
\scalebox{0.86}{
\bfig
\putmorphism(-150,500)(1,0)[(B,A)_1`(B',A)_1`(k,A)_1]{600}1a
 \putmorphism(480,500)(1,0)[\phantom{F(A)}`(B',A)_2 `\alpha(A)_{B'}]{640}1a

 \putmorphism(-150,50)(1,0)[(B,\tilde A)_1`(B',\tilde A)_1`(k,\tilde A)_1]{600}1a
 \putmorphism(480,50)(1,0)[\phantom{F(A)}`(B',\tilde A)_2 `\alpha(\tilde A)_{B'}]{640}1a

\putmorphism(-180,500)(0,-1)[\phantom{Y_2}``(B,U)_1]{450}1l
\putmorphism(450,500)(0,-1)[\phantom{Y_2}``]{450}1r
\putmorphism(300,500)(0,-1)[\phantom{Y_2}``(B',U)_1]{450}0r
\putmorphism(1100,500)(0,-1)[\phantom{Y_2}``(B',U)_2]{450}1r
\put(-40,280){\fbox{$(k,U)_1$}}
\put(700,280){\fbox{$\alpha_U(B')$}}

\putmorphism(-150,-400)(1,0)[(B,\tilde A)_1`(B,\tilde A)_2 `\alpha(\tilde A)_B]{640}1a
 \putmorphism(480,-400)(1,0)[\phantom{A'\ot B'}` (B',\tilde A)_2 `(k,\tilde A)_2]{680}1a

\putmorphism(-180,50)(0,-1)[\phantom{Y_2}``=]{450}1l
\putmorphism(1120,50)(0,-1)[\phantom{Y_3}``=]{450}1r
\put(320,-200){\fbox{$\delta_{\alpha(\tilde A),k}$}}
\efig}
\quad=\quad
\scalebox{0.86}{
\bfig
\putmorphism(-150,500)(1,0)[(B,A)_1`(B',A)_1`(k,A)_1]{600}1a
 \putmorphism(480,500)(1,0)[\phantom{F(A)}`(B',A)_2 `\alpha(A)_{B'}]{640}1a
 \putmorphism(-150,50)(1,0)[(B,A)_1`(B,A)_2`\alpha(A)_B]{600}1a
 \putmorphism(480,50)(1,0)[\phantom{F(A)}`(B',A)_2 `(k,A)_2]{680}1a

\putmorphism(-180,500)(0,-1)[\phantom{Y_2}``=]{450}1r
\putmorphism(1100,500)(0,-1)[\phantom{Y_2}``=]{450}1r
\put(320,280){\fbox{$\delta_{\alpha(A),k}$}}

\putmorphism(-150,-400)(1,0)[(B,\tilde A)_1`(B,\tilde A)_2 `\alpha(\tilde A)_{B}]{640}1a
 \putmorphism(490,-400)(1,0)[\phantom{F(A')}` (B',\tilde A)_2 `(k,\tilde A)_2]{640}1a

\putmorphism(-180,50)(0,-1)[\phantom{Y_2}``(B,U)_1]{450}1l
\putmorphism(480,50)(0,-1)[\phantom{Y_2}``]{450}1l
\putmorphism(610,50)(0,-1)[\phantom{Y_2}``(B,U)_2]{450}0l 
\putmorphism(1120,50)(0,-1)[\phantom{Y_3}``(B',U)_2]{450}1r
\put(-40,-180){\fbox{$\alpha_U(B)$}} 
\put(620,-180){\fbox{$(k,U)_2$}}
\efig}
$$
and 
$$
\scalebox{0.86}{
\bfig
 \putmorphism(-150,500)(1,0)[(B,A)`(B,A) `=]{600}1a
 \putmorphism(450,500)(1,0)[(B,A)` `\alpha(A)_B]{450}1a
\putmorphism(-180,500)(0,-1)[\phantom{Y_2}`(B, \tilde A) `(B,U)_1]{450}1l
\put(-50,50){\fbox{$(u,U)_1$}}
\putmorphism(-150,-400)(1,0)[(\tilde B, \tilde A)` `=]{500}1a
\putmorphism(-180,50)(0,-1)[\phantom{Y_2}``(u,\tilde A)_1]{450}1l
\putmorphism(450,50)(0,-1)[\phantom{Y_2}`(\tilde B, \tilde A)`(\tilde B, U)_1]{450}1l
\putmorphism(450,500)(0,-1)[\phantom{Y_2}`(\tilde B, A) `(u,A)_1]{450}1l
\put(600,260){\fbox{$\alpha(A)_u$}}
\putmorphism(450,50)(1,0)[\phantom{(B, \tilde A)}``\alpha(A)_{\tilde B}]{500}1a
\putmorphism(1070,50)(0,-1)[\phantom{(B, A')}`(\tilde B', \tilde A)`(\tilde B,U)_2]{450}1r
\putmorphism(1070,500)(0,-1)[(B', A)`(\tilde B', A)`(u,A)_2]{450}1r
\putmorphism(450,-400)(1,0)[\phantom{(B, \tilde A)}``\alpha(\tilde A)_{\tilde B}]{500}1a
\put(600,-170){\fbox{$\alpha_U(\tilde B)$}}
\efig}=
\scalebox{0.86}{
\bfig
 \putmorphism(-150,500)(1,0)[(B,A)_1`(B,A)_2 `\alpha(A)_B]{600}1a
 \putmorphism(450,500)(1,0)[\phantom{(B,A)}` `=]{520}1a
\putmorphism(-180,500)(0,-1)[\phantom{Y_2}`(B, \tilde A)_1 `(B,U)_1]{450}1l
\put(0,260){\fbox{$\alpha_U(B)$}}
\putmorphism(-180,-400)(1,0)[(\tilde B, \tilde A)` `\alpha(\tilde A)_{\tilde B}]{500}1a
\putmorphism(-180,50)(0,-1)[\phantom{Y_2}``(u,\tilde A)_1]{450}1l
\putmorphism(450,50)(0,-1)[\phantom{Y_2}`(\tilde B', \tilde A)`(u,\tilde A)_2]{450}1r
\putmorphism(450,500)(0,-1)[\phantom{Y_2}`(B', \tilde A) `(B,U)_2]{450}1r
\put(600,50){\fbox{$(u,U)_2$}}
\putmorphism(-180,50)(1,0)[\phantom{(B, \tilde A)}``\alpha(\tilde A)_B]{500}1a
\putmorphism(1110,50)(0,-1)[\phantom{(B, A')}`(\tilde B', \tilde A)`(\tilde B,U)_2]{450}1r
\putmorphism(1110,500)(0,-1)[(B', A)`(\tilde B', A)`(u,A)_2]{450}1r
\putmorphism(450,-400)(1,0)[\phantom{(B, \tilde A)}``=]{500}1b
\put(0,-170){\fbox{$\alpha(\tilde A)_u$}}
\efig}
$$
which by additional identifications $\theta^A_u=\alpha(A)_u$ and $\theta^B_U=\alpha_U(B)$ are $(HOT^q_3)$ and $(HOT^q_4)$. 
\qed\end{proof}

Now we have that $\alpha:F\Rightarrow G$ yields $\theta: (-,-)_1\Rightarrow(-,-)_2$. Before seeing the converse, let us summarize our above findings:

\begin{prop} \prlabel{horiz nat summed up}
A horizontal oplax transformation $\alpha:F\Rightarrow G$ between lax double functors $F,G: \Aa\to\llbracket\Bb,\Cc\rrbracket$ consists of the following data:
\begin{itemize}
\item $\alpha(A): F(A)\to G(A)$ is a horizontal oplax transformation between lax double functors for every $A\in\Aa$;
\item $\alpha_K: \frac{(-,K)_1}{\alpha(A')}\Rrightarrow \frac{\alpha(A)}{(-,K)_2}$ is a (globular) modification for every 1h-cell $K:A\to A'$; 
\item $\alpha_U$ (of the form \equref{alfa-U}) is a modification for every 1v-cell $U:A\to\tilde A$; 
\end{itemize} 
so that $\alpha_K$ and $\alpha_U$ obey five axioms, which (after evaluation at $B\in\Bb$) yield that $\alpha(-)_B: F(-)(B)\to G(-)(B)$ 
is a horizontal oplax transformation between lax double functors for every $B\in\Bb$ (by setting $\alpha(K)_B:=\alpha_K(B)$ and $\alpha(U)_B:=\alpha_U(B)$). 
(Both modifications above are meant in the sense of \deref{modif-hv}.) 
\end{prop}

Now, assuming that $\theta: (-,-)_1\Rightarrow(-,-)_2$ given by pairs $(\theta^A, \theta^B)$ for $A\in\Aa, B\in\Bb$ is a horizontal oplax 
transformation between lax double quasi-functors, whose corresponding lax double functors are $F,G:\Aa\to \llbracket\Bb,\Cc\rrbracket$, we 
define a horizontal oplax transformation $\alpha:F\Rightarrow G$ as follows. For $A\in\Aa$ set $\alpha(A):=\theta^A$, for a 1h-cell 
$K:A\to A'$ set for the desired globular modification $\alpha_K$ to be given by components $\alpha_K(B):=\theta^B_K$, and for a 1v-cell 
$U:A\to\tilde A$ set for the desired modification $\alpha_U$ to be given by components $\alpha_U(B):=\theta^B_U$. Since moreover 
$\alpha(A)_B=\theta^A_B=\theta^B_A$, we have that $\alpha(-)_B:=\theta^B$ is a horizontal oplax transformation of lax double functors 
with $\alpha(K)_B:=\theta^B_K=\alpha_K(B)$ and $\alpha(U)_B:=\theta^B_U=\alpha_U(B)$. Now by \prref{equiv-modif} $\alpha_K$ and 
$\alpha_U$ are modifications. By \prref{horiz nat summed up} we have that $\alpha:F\Rightarrow G$ is indeed a horizontal oplax 
transformation between lax double functors.

\medskip

The two assignments of horizontal oplax transformations are clearly inverse to each other. 

\smallskip

The 1-1 correspondence between 1v-cells works completely analogously as for 1h-cells. (This time modifications $(\alpha_0)^u$ of the type 
\equref{m-hor} are used instead of $\alpha_K$ in the analogous place in the above two Propositions.)

\subsection{The 1-1 correspondence among 2-cells} \sslabel{2-cells}

A modification $\Theta$ in $\Lax_{hop}(\Aa,\llbracket\Bb,\Cc\rrbracket)$, {\em i.e.} a modifications between two horizontally oplax and two 
vertically lax transformations, is given by 2-cells 
$$\scalebox{0.86}{
\bfig
\putmorphism(-180,50)(1,0)[F(A)` G(A)`\alpha(A)_{-}]{550}1a
\putmorphism(-180,-270)(1,0)[F\s'(A)`G'(A) `\beta(A)_{-}]{550}1b
\putmorphism(-170,50)(0,-1)[\phantom{Y_2}``\alpha_0(A)_{-}]{320}1l
\putmorphism(350,50)(0,-1)[\phantom{Y_2}``\beta_0(A)_{-}]{320}1r
\put(0,-140){\fbox{$(\Theta_A)_{-}$}}
\efig}
$$
in $\llbracket\Bb,\Cc\rrbracket$ with free slots in $B\in\Bb$ which satisfy axioms (m.ho-vl.-1) \label{m.ho-vl.-1} and (m.ho-vl.-2)\label{m.ho-vl.-2}.  
Evaluating the latter two axioms in a fixed $B\in\Bb$, and considering the slot occupied by 
0-, 1h- and 1v-cells in $\Aa$ as variable, these two axioms mean that one has a modification $(\Theta_{-})_B$:
$$\scalebox{0.86}{
\bfig
\putmorphism(-180,50)(1,0)[F(A)` G(A)`\alpha(-)_B]{550}1a
\putmorphism(-180,-270)(1,0)[F\s'(A)`G'(A). `\beta(-)_B]{550}1b
\putmorphism(-170,50)(0,-1)[\phantom{Y_2}``\alpha_0(-)_B]{320}1l
\putmorphism(350,50)(0,-1)[\phantom{Y_2}``\beta_0(-)_B]{320}1r
\put(0,-140){\fbox{$(\Theta_{-})_B$}}
\efig}
$$
We now identify $\tau^A=\Theta_A$ and $\tau^B=(\Theta_{-})_B$ and recall the identifications from \ssref{1-cells}: 
$\theta^A=\alpha(A), \theta^B=\alpha(-)_B$, and analogously we have $\theta'^A=\beta(A), \theta'^B=\beta(-)_B$, 
and similarly for the vertical lax transformations: 
$\theta_0^A=\alpha_0(A), \theta_0^B=\alpha_0(-)_B$, and $\theta_0'^A=\beta_0(A), \theta_0'^B=\beta_0(-)_B$.  
Then we clearly have: $\tau^A_B=\tau^B_A$ and hence that $\tau^A$ and $\tau^B$ make a 
modification $\tau$ of horizontal oplax transformations $\theta=(\theta^A, \theta^B)_{\substack{A\in\Aa \\B\in\Bb}}$ and $\theta'=
(\theta'^A, \theta'^B)_{\substack{A\in\Aa\\B\in\Bb}}$ and vertical lax transformations 
$\theta_0=(\theta_0^A, \theta_0^B)_{\substack{A\in\Aa \\B\in\Bb}}$ and $\theta_0'=
(\theta_0'^A, \theta_0'^B)_{\substack{A\in\Aa\\B\in\Bb}}$ between lax double quasi-functors (recall the last two squares in \equref{q-modif}). 

Reading the above characterization of a modification $\Theta$ and how we obtained the modification $\tau$ in the reversed order, 
one finds the converse assignment, and it is clear that these two assignments are inverse to each other. 
 
\medskip

It is directly seen that the assignments that we defined in \ssref{1-cells} and this Subsection determine a strict double functor between double categories 
$\Lax_{hop}(\Aa, \llbracket\Bb,\Cc\rrbracket)$ and $q\x\Lax_{hop}(\Aa\times\Bb,\Cc)$. To see that it is compatible with compositions viewing 
\prref{horiz nat summed up} may be helpful. We conclude that there is an isomorphism of double categories
\begin{equation} \eqlabel{pre-Gray}
q\x\Lax_{hop}(\Aa\times\Bb,\Cc) \iso \Lax_{hop}(\Aa, \llbracket\Bb,\Cc\rrbracket).
\end{equation}

\section{A double functor from $q\x\Lax_{hop}^{st}(\Aa\times\Bb,\Cc)$ to $\Lax_{hop}(\Aa\times\Bb,\Cc)$ } \selabel{isom-spec}

Let $q\x\Lax_{hop}^{st}(\Aa\times\Bb,\Cc)$ be the full double subcategory of $q\x\Lax_{hop}(\Aa\times\Bb,\Cc)$ 
differing only in 0-cells, so that the 2-cells $(u,U)$ of its lax double quasi-functors are trivial. The supra-index {\em ``st''} 
alludes to {\em strict} from the notion of a vertical strict transformation (from the 1v-cells in $\llbracket\Bb,\Cc\rrbracket$). 
(The corresponding full double subcategory of $\Lax_{hop}(\Aa, \llbracket\Bb,\Cc\rrbracket)$ isomorphic to $q\x\Lax_{hop}^{st}(\Aa\times\Bb,\Cc)$ 
in \equref{pre-Gray} we denote by $\Lax_{hop}(\Aa, \llbracket\Bb,\Cc\rrbracket^{st})$. Here $\llbracket\Bb,\Cc\rrbracket^{st}$ 
stands for the variation of $\llbracket\Bb,\Cc\rrbracket$ 
in which the vertical lax transformations are strict.)

We can prove that there is a double functor from $q\x\Lax_{hop}^{st}(\Aa\times\Bb,\Cc)$ 
(and thus also from $\Lax_{hop}(\Aa, \llbracket\Bb,\Cc\rrbracket^{st})$) 
to the double category $\Lax_{hop}(\Aa\times\Bb,\Cc)$, consisting of lax double functors 
on the Cartesian product of double categories, and their corresponding horizontal oplax and vertical lax transformations and modifications. 
We will denote it by 
\begin{equation} \eqlabel{2-functor F}
\F: q\x\Lax_{hop}^{st}(\Aa\times\Bb,\Cc) \to \Lax_{hop}(\Aa\times\Bb,\Cc).
\end{equation}
Moreover, restricting to certain double subcategories of $q\x\Lax_{hop}^{st}(\Aa\times\Bb,\Cc)$ and $\Lax_{hop}(\Aa\times\Bb,\Cc)$ 
we obtain a double equivalence of double categories. 
For this purpose in this Section we will construct a tuple $(\F\s',\G, \kappa, \lambda)$ 
of double equivalence functors and horizontal strict transformations $\kappa: \Id\Rightarrow\G\F'$ and $\lambda: \F'\G\Rightarrow\Id$.  
The results that we obtain will generalize to the context of double categories Theorems 4.10 and 5.3 of \cite{FMS}. 

\subsection{The double functor $\F$ on 0-cells} \sslabel{F on 0}

Let us show that a lax double quasi-functor $H:\Aa\times\Bb\to\Cc$, with lax double functors $H(A,-)=(-, A)$ and $H(-, B)=(B,-)$, 
{\em whose 2-cells $(u,U)$ are identities} determines a lax double functor $P:\Aa\times\Bb\to\Cc$ on the Cartesian product. 

Instead of typing the whole proof, we will indicate the list of its steps. For that purpose recall the notation for the formulaic 
computations $[\alpha\vert\beta]=\beta\alpha$ for the horizontal composition of 2-cells $\alpha$ (first) and $\beta$ (second) from 
the end of the second paragraph of \seref{mon str}.

We set: \\
$P(A,B)=H(A,B), \\
P(K,k)=H(A',k)H(K,B)=(k,A')(B,K)$, for $K:A\to A'$ and $k:B\to B'$, \\
$P(U,u)=\frac{(B,U)}{(u,\tilde A)}$ for 1v-cells $U:A\to\tilde A, u: B\to\tilde B$, and \\
\begin{equation} \eqlabel{P on 2-cells}
P(\alpha,\beta):= 
\scalebox{0.86}{
\bfig

 \putmorphism(-520,50)(1,0)[(B,A)`(B,A')`(B,K)]{600}1a
 \putmorphism(30,50)(1,0)[\phantom{A''\ot B'}` (B', A') `(k, A')]{670}1a

\putmorphism(70,50)(0,-1)[\phantom{Y_2}`(B,\tilde A')`]{450}1r
\putmorphism(50,180)(0,-1)[``(B,U')]{450}0r
\putmorphism(650,50)(0,-1)[\phantom{Y_2}`\phantom{A''\ot B'} `(B',U')]{450}1r
\putmorphism(70,-400)(0,-1)[\phantom{Y_2}``]{450}1r
\putmorphism(70,-490)(0,-1)[\phantom{Y_2}``(u,\tilde A')]{450}0r
\putmorphism(650,-400)(0,-1)[\phantom{Y_2}``(u', \tilde A')]{450}1r

\put(210,-580){\fbox{$(\beta, \tilde A')$}}
\put(210,-210){\fbox{$(k,U')$}}

\putmorphism(-470,50)(0,-1)[\phantom{Y_2}`(B,\tilde A)`(B,U)]{450}1l
\putmorphism(-470,-400)(0,-1)[\phantom{Y_2}`(\tilde B,\tilde A)`(u,\tilde A)]{450}1l
\put(-350,-150){\fbox{$(B,\alpha)$}}
\put(-350,-650){\fbox{$(u,\tilde K)$}}

 \putmorphism(-510,-400)(1,0)[\phantom{A''\ot B'}`\phantom{A''\ot B'}`(B,\tilde K)]{600}1a
 \putmorphism(30,-400)(1,0)[\phantom{A''\ot B'}` \phantom{A''\ot B'} `]{670}1a
 \putmorphism(30,-410)(1,0)[\phantom{A''\ot B'}` (\tilde B,\tilde A')`(k, \tilde A')]{670}0a
 \putmorphism(-530,-850)(1,0)[\phantom{A''\ot B'}`(\tilde B,\tilde A')`(\tilde B,\tilde K)]{600}1b
\putmorphism(30,-850)(1,0)[\phantom{A''\ot B'}` (B',\tilde A') `(\tilde k, \tilde A')]{670}1b
\efig 
}
\end{equation} 
for 2-cells $\alpha$ in $\Aa$ and $\beta$ in $\Bb$ as in \equref{alfa-beta}; \\
for the lax structure $\gamma_{(f',g')(f,g)}: P(f',g')P(f,g)\Rightarrow P(f'f,g'g)$ and $\iota^P: 1_{P(A,B)}\Rightarrow P(1d_{(A,B)})$ 
of $P$ we set: 
$$(k,K):=
\scalebox{0.86}{
\bfig
\putmorphism(210,350)(1,0)[` `(k,A')]{420}1a
\putmorphism(530,350)(1,0)[\phantom{F(A)}` `(B',K')]{450}1a

\putmorphism(-190,10)(1,0)[` `(B,K)]{420}1a
\putmorphism(210,10)(1,0)[` `(B,K')]{420}1a
\putmorphism(530,10)(1,0)[\phantom{F(A)}` `(k,A'')]{450}1a
\putmorphism(900,10)(1,0)[\phantom{F(A)}` `(k',A'')]{450}1a

\putmorphism(210,350)(0,-1)[\phantom{Y_2}``=]{350}1l
\putmorphism(980,350)(0,-1)[\phantom{Y_2}``=]{350}1r

\putmorphism(-180,0)(0,-1)[\phantom{Y_2}``=]{350}1l
\putmorphism(620,0)(0,-1)[\phantom{Y_2}``=]{350}1r
\putmorphism(1340,0)(0,-1)[\phantom{Y_2}``=]{350}1r

\put(460,200){\fbox{$(k,K')$}}

 \putmorphism(-190,-320)(1,0)[``(B,K'K)]{800}1b
 \putmorphism(530,-320)(1,0)[\phantom{F(A)}` `(k'k,A'')]{800}1b

\put(0,-180){\fbox{$(B,-)_{K'K}$}}
\put(760,-180){\fbox{$(-,A'')_{k'k}$}}
\efig}
\quad\text{and}\quad
\iota^P_{(A,B)}:=
\scalebox{0.86}{
\bfig
 \putmorphism(-210,220)(1,0)[(B,A)`\phantom{F(A)} `=]{500}1a
\putmorphism(330,220)(1,0)[(B,A)`(B,A)`=]{500}1a
\putmorphism(-210,220)(0,-1)[\phantom{Y_2}``=]{370}1l
\putmorphism(320,220)(0,-1)[\phantom{Y_2}``=]{370}1l
\putmorphism(810,220)(0,-1)[\phantom{Y_2}``=]{370}1r
 \putmorphism(-210,-150)(1,0)[(B,A)`\phantom{Y_2}`(B,1_A)]{470}1a
 \put(490,50){\fbox{$\iota^B_A$}} 
\putmorphism(320,-150)(1,0)[(B,A)`(B,A) `(1_B,A)]{580}1a
\put(-60,40){\fbox{$\iota^A_B$}}
\efig}
$$
where $\iota^B_A=(B,-)_A$ and $\iota^A_B=(-,A)_B$ of $H$. 

The hexagonal law (lx.f.hex) for $\gamma$ and the unital laws (lx.f.u) can be formulated in the underlying horizontal 2-category, where they 
amount to the same data as in \cite[Theorem 3.2]{FMS}. The same holds for 
preservation of identity 2-cell on a 1h-cell (lx.f.h2). We will discuss vertical naturality of $\iota^P$ with respect to 1v-cells (lx.f.u-nat), the 
naturality of $\gamma$ with respect to 2-cells (lx.f.c-nat), and the vertical functoriality of $P$ with respect to 2-cells (lx.f.h1). $P$ is strictly 
compatible with vertical composition of 1v-cells and vertical identities on objects, rules (lx.f.v1) and (lx.f.v2), 
since $(B,-)$ and $(-,A)$ are strict in the vertical direction.

For (lx.f.c-nat) we take two 2-cells $\alpha, \alpha'$ in $\Aa$ and two 2-cells $\beta, \beta'$ in $\Bb$:
\begin{equation} \eqlabel{alfa-beta}
\scalebox{0.86}{
\bfig
\putmorphism(-150,50)(1,0)[A` A'`K]{450}1a
\putmorphism(-150,-250)(1,0)[\tilde A`\tilde A' `\tilde K]{440}1b
\putmorphism(-170,50)(0,-1)[\phantom{Y_2}``U]{300}1l
\putmorphism(280,50)(0,-1)[\phantom{Y_2}``U']{300}1r
\put(-20,-110){\fbox{$\alpha$}}
\putmorphism(300,50)(1,0)[\phantom{Y}` A''`K']{450}1a
\putmorphism(730,50)(0,-1)[\phantom{Y_2}``U'']{300}1r
\putmorphism(300,-250)(1,0)[\phantom{Y}`\tilde A' `\tilde K]{440}1b
\put(460,-110){\fbox{$\alpha'$}}
\efig}
\quad\text{and}\quad
\scalebox{0.86}{
\bfig
\putmorphism(-150,50)(1,0)[B` B'`k]{450}1a
\putmorphism(-150,-250)(1,0)[\tilde B`\tilde B' `\tilde k]{440}1b
\putmorphism(-170,50)(0,-1)[\phantom{Y_2}``u]{300}1l
\putmorphism(280,50)(0,-1)[\phantom{Y_2}``u']{300}1r
\put(0,-110){\fbox{$\beta$}}
\putmorphism(300,50)(1,0)[\phantom{Y}` B''`k']{450}1a
\putmorphism(730,50)(0,-1)[\phantom{Y_2}``u'']{300}1r
\putmorphism(300,-250)(1,0)[\phantom{Y}`\tilde B' `\tilde k]{440}1b
\put(460,-110){\fbox{$\beta'$}}
\efig}
\end{equation}
and we should prove the equality: 
$$
\scalebox{0.86}{
\bfig

 \putmorphism(-150,450)(1,0)[(B,A')`(B',A')`(k,A')]{600}1a
 \putmorphism(450,450)(1,0)[\phantom{A\ot B}`(B', A'') `(B',K')]{680}1a

 \putmorphism(-750,0)(1,0)[(B,A)`(B,A')`(B,K)]{600}1a
 \putmorphism(-150,0)(1,0)[\phantom{A\ot B}`(B,A'')`(B,K')]{600}1a
 \putmorphism(450,0)(1,0)[\phantom{A\ot B}`(B', A'') `(k,A'')]{680}1a
 \putmorphism(1120,0)(1,0)[\phantom{A'\ot B'}`(B'', A'') `(k', A'')]{720}1a

\putmorphism(-750,0)(0,-1)[\phantom{Y_2}``=]{450}1l
\putmorphism(420,0)(0,-1)[\phantom{Y_2}``=]{420}1l
\putmorphism(1800,0)(0,-1)[\phantom{Y_2}``=]{450}1r

\putmorphism(-180,450)(0,-1)[\phantom{Y_2}``=]{450}1r
\putmorphism(1100,450)(0,-1)[\phantom{Y_2}``=]{450}1r
\put(350,210){\fbox{$(k,K')$}}
\put(-350,-230){\fbox{$(B,-)_{K'K}$}}
\put(950,-230){\fbox{$(-,A'')_{k'k}$}}

 \putmorphism(-750,-450)(1,0)[(B,A)`(B,A'')`(B,K'K)]{1200}1a
 \putmorphism(420,-450)(1,0)[\phantom{A''\ot B'}` (B'', A'') `(k'k, A'')]{1390}1a

\putmorphism(420,-450)(0,-1)[\phantom{Y_2}`(B,\tilde A'')`(B,U'')]{450}1r
\putmorphism(420,-900)(0,-1)[\phantom{Y_2}``(u,\tilde A'')]{450}1r

\putmorphism(-750,-450)(0,-1)[\phantom{Y_2}`(B,\tilde A'')`(B,U)]{450}1l
\putmorphism(-750,-900)(0,-1)[\phantom{Y_2}`(\tilde B,\tilde A)`(u,\tilde A)]{450}1l
\putmorphism(1800,-450)(0,-1)[\phantom{Y_2}`(B'',\tilde A'')`(B'',U'')]{450}1r
\putmorphism(1800,-900)(0,-1)[\phantom{Y_2}`(\tilde B'',\tilde A'')`(u'',\tilde A'')]{450}1r

\put(-350,-650){\fbox{$(B,\alpha'\alpha)$}}
\put(950,-650){\fbox{$(k'k,U'')$}}
\put(-350,-1100){\fbox{$(u,\tilde K'\tilde K)$}}
\put(950,-1100){\fbox{$(\beta'\beta, \tilde A'')$}}

 \putmorphism(-780,-900)(1,0)[\phantom{A''\ot B'}`\phantom{A''\ot B'}`(B,\tilde K'\tilde K)]{1200}1a
 \putmorphism(400,-900)(1,0)[\phantom{A''\ot B'}` \phantom{A''\ot B'} ` (k'k, \tilde A'')]{1390}1a

 \putmorphism(-780,-1350)(1,0)[\phantom{A''\ot B'}`(\tilde B,\tilde A'')`(\tilde B,\tilde K'\tilde K)]{1200}1a
 \putmorphism(400,-1350)(1,0)[\phantom{A''\ot B'}` \phantom{A''\ot B'} `(\tilde k'\tilde k, \tilde A'')]{1390}1a
\efig}
$$
=
$$
\scalebox{0.86}{
\bfig

 \putmorphism(30,-450)(1,0)[(B,A)`(B,A')`(B,K)]{600}1a
 \putmorphism(600,-450)(1,0)[\phantom{A''\ot B'}` (B', A') `(k, A')]{700}1a
 \putmorphism(1280,-450)(1,0)[\phantom{A''\ot B'}` (B', A'') `(B', K')]{800}1a
 \putmorphism(2080,-450)(1,0)[\phantom{A''\ot B'}` (B'', A'') `(k',A'')]{700}1a

\putmorphism(70,-450)(0,-1)[\phantom{Y_2}`(B,\tilde A)`(B,U)]{450}1l
\putmorphism(600,-450)(0,-1)[\phantom{Y_2}`(B,\tilde A')`]{450}1r
\putmorphism(590,-320)(0,-1)[``(B,U')]{450}0r
\putmorphism(590,-900)(0,-1)[\phantom{Y_2}`(\tilde B, \tilde A')`]{450}1r
\putmorphism(590,-790)(0,-1)[``(u,\tilde A')]{450}0r
\putmorphism(1300,-900)(0,-1)[\phantom{Y_2}`(\tilde B',\tilde A')`]{450}1r
\putmorphism(1320,-790)(0,-1)[``(u',\tilde A')]{450}0l
 \putmorphism(560,-1350)(1,0)[\phantom{A''\ot B'}` \phantom{A''\ot B'} `(\tilde k, \tilde A')]{750}1b
 \putmorphism(1270,-1350)(1,0)[\phantom{A''\ot B'}` (\tilde B',\tilde A'') `(\tilde B', \tilde K')]{800}1b
 \putmorphism(2050,-1350)(1,0)[\phantom{A''\ot B'}` \phantom{A''\ot B'} `(\tilde k',\tilde A'')]{700}1a 

\put(180,-650){\fbox{$(B,\alpha)$}}
\put(770,-690){\fbox{$(k,U')$}}
\put(770,-1210){\fbox{$(\beta,\tilde A')$}}

\putmorphism(1300,-450)(0,-1)[\phantom{Y_2}`(B',\tilde A')`]{450}1r
\putmorphism(1320,-320)(0,-1)[``(B',U')]{450}0l 

\putmorphism(70,-900)(0,-1)[\phantom{Y_2}`(\tilde B,\tilde A)`(u,\tilde A)]{450}1l

\putmorphism(2130,-450)(0,-1)[\phantom{Y_2}``]{450}1l 
\putmorphism(2160,-320)(0,-1)[``(B',U'')]{450}0l
 \putmorphism(1280,-900)(1,0)[\phantom{A''\ot B'}` (B',\tilde A'') `(B', \tilde K')]{750}1a
 \putmorphism(2000,-900)(1,0)[\phantom{A''\ot B'}` \phantom{A''\ot B'} `(k', \tilde A'')]{750}1a  

\put(180,-1100){\fbox{$(u,\tilde K)$}}
\putmorphism(2150,-900)(0,-1)[\phantom{Y_2}``]{450}1r 
\putmorphism(2160,-790)(0,-1)[``(u',\tilde A'')]{450}0l

\putmorphism(2730,-450)(0,-1)[\phantom{Y_2}`(B',\tilde A'')`(B'', U'')]{450}1r 
\putmorphism(2730,-900)(0,-1)[\phantom{Y_2}`(\tilde B',\tilde A'')`(u'',\tilde A'')]{450}1r 

 \putmorphism(0,-900)(1,0)[\phantom{A''\ot B'}`\phantom{A''\ot B'}`(B,\tilde K)]{600}1a
 \putmorphism(560,-900)(1,0)[\phantom{A''\ot B'}` \phantom{A''\ot B'} `(k, \tilde A')]{750}1a

 \putmorphism(30,-1350)(1,0)[\phantom{A''\ot B'}`\phantom{A''\ot B'}`(\tilde B,\tilde K)]{600}1a
\put(1520,-690){\fbox{$(B',\alpha')$}}
\put(1520,-1210){\fbox{$(u', \tilde K')$}}

\put(2260,-650){\fbox{$(k', U'')$}}
\put(2260,-1100){\fbox{$(\beta',\tilde A'')$}}

\putmorphism(590,-1350)(0,-1)[\phantom{Y_2}`(\tilde B,\tilde A')`=]{450}1r
\putmorphism(2130,-1350)(0,-1)[\phantom{Y_2}`(\tilde B',\tilde A'')`=]{450}1r 
\putmorphism(550,-1800)(1,0)[\phantom{A''\ot B'}` (\tilde B, \tilde A'') `(\tilde B, \tilde K')]{840}1a
\putmorphism(1360,-1800)(1,0)[\phantom{A''\ot B'}` \phantom{A''\ot B'} `(\tilde k, \tilde A'')]{780}1a
\put(1180,-1600){\fbox{$(\tilde k,\tilde K')$}}

\putmorphism(-10,-1800)(1,0)[(\tilde B, \tilde A)`\phantom{A''\ot B'}`(\tilde B,\tilde K)]{640}1a
\putmorphism(2120,-1800)(1,0)[\phantom{A''\ot B'}` (\tilde B'', \tilde A'') `(\tilde k', \tilde A'')]{750}1a 

\putmorphism(0,-1800)(0,-1)[\phantom{Y_2}``=]{410}1l
\putmorphism(1320,-1800)(0,-1)[\phantom{Y_2}``=]{410}1r
\putmorphism(2840,-1800)(0,-1)[\phantom{Y_2}``=]{410}1r

\put(220,-2000){\fbox{$(\tilde B,-)_{\tilde K'\tilde K}$}}
\put(2140,-2000){\fbox{$(-,\tilde A'')_{\tilde k'\tilde k}$}} 

 \putmorphism(0,-2230)(1,0)[(\tilde B,\tilde A)`(\tilde B,\tilde A'')`(\tilde B,\tilde K'\tilde K)]{1400}1a
 \putmorphism(1370,-2230)(1,0)[\phantom{A''\ot B'}` (\tilde B'', \tilde A'') `(\tilde k'k, \tilde A'')]{1500}1a
\efig}
$$
In order to show that the left hand-side equals the right one, apply the following rules: 
1) naturality of the laxity of $(B,-)$ in $\alpha'\alpha$ and $((k'k,U))$; 
2) ($(k,K)$-r-nat), $((u,K'K))$ and naturality of $(-,A'')$ in $\beta'\beta$, and  
3) ($(k,K)$-l-nat). 

\bigskip

For (lx.f.u-nat) we should prove the equality:
$\frac{\iota^P_{(A,B)}}{P(1^{(U,u)})}=\frac{1^{P(U,u)}}{\iota^P_{(\tilde A,\tilde B)}}$, which translates into:
$$
\scalebox{0.86}{
\bfig
 \putmorphism(-170,650)(1,0)[(B,A)`(B,A)`=]{480}1a
\putmorphism(-170,280)(1,0)[(B,A)`(B,A) `(B,1_A)]{520}1a

 \putmorphism(330,650)(1,0)[\phantom{(B,A)}`(B,A)`=]{500}1a
\putmorphism(330,280)(1,0)[\phantom{(B,A)}`(B,A) `(1_B,A)]{520}1a

\putmorphism(-170,650)(0,-1)[\phantom{Y_2}``=]{350}1l
\putmorphism(300,650)(0,-1)[\phantom{Y_2}``=]{370}1r
\put(-100,480){\fbox{$(B,-)_A$}}
\put(400,480){\fbox{$(-,A)_B$}}
\putmorphism(830,650)(0,-1)[\phantom{Y_2}``=]{350}1r

\putmorphism(400,-70)(1,0)[\phantom{Y_2}`\phantom{Y_2} `(1_{\tilde B},\tilde A)]{380}1a

\putmorphism(-180,280)(0,-1)[\phantom{Y_2}``(B,U)]{350}1l %
 \putmorphism(300,280)(0,-1)[\phantom{Y_2}``]{370}1r

\putmorphism(830,280)(0,-1)[\phantom{Y_2}`(\tilde B,\tilde A)`(B,\tilde U)]{350}1r 

\putmorphism(-150,-70)(1,0)[(B,\tilde A)`(B,\tilde A)`(B,1_{\tilde A})]{500}1a
\put(-100,120){\fbox{$(B,1^U)$}}
\put(400,120){\fbox{$(1_B,U)$}}

\putmorphism(-170,-70)(0,-1)[\phantom{Y_2}`(\tilde B,\tilde A)`(u,\tilde A)]{370}1l
 \putmorphism(300,-70)(0,-1)[\phantom{Y_2}`(\tilde B,\tilde A)`]{370}1l
\putmorphism(830,-70)(0,-1)[\phantom{Y_2}`(\tilde B,\tilde A)`(u,\tilde A) ]{370}1r
\put(-100,-240){\fbox{$(u,1_{\tilde A})$}}
\put(400,-240){\fbox{$(1^u,\tilde A)$}}
\putmorphism(-70,-440)(1,0)[``(\tilde B,1_{\tilde A})]{250}1a
\putmorphism(400,-440)(1,0)[``(1_{\tilde B},\tilde A)]{300}1a
\efig}\quad
=\quad
\scalebox{0.86}{
\bfig
\putmorphism(-170,650)(1,0)[\phantom{(B,A)}`\phantom{(B,A)}`=]{1030}1a
\putmorphism(-170,650)(0,-1)[(B,A)`(B,\tilde A)`(B,U)]{300}1l
\putmorphism(-170,350)(0,-1)[\phantom{Y_2}``(u,\tilde A)]{300}1l

\putmorphism(-170,350)(1,0)[\phantom{(B,A)}`\phantom{(B,A)}`=]{1030}1a

\putmorphism(860,650)(0,-1)[(B, A)`(B,\tilde A)`(B,U)]{300}1r
\putmorphism(860,350)(0,-1)[\phantom{Y_2}``(u,\tilde A)]{300}1r
\put(200,480){\fbox{$1^{(B,U)}$}}
\put(190,180){\fbox{$1^{(u,\tilde A)}$}}

\putmorphism(-170,50)(1,0)[(\tilde B, \tilde A)``=]{450}1a
\putmorphism(300,50)(1,0)[`(\tilde B, \tilde A)`=]{570}1a
\putmorphism(-170,50)(0,-1)[\phantom{Y_2}``=]{300}1l
\putmorphism(280,50)(0,-1)[\phantom{Y_2}``=]{300}1r
\put(-120,-110){\fbox{$(\tilde B, -)_{\tilde A}$}}

\putmorphism(860,50)(0,-1)[\phantom{Y_2}``=]{300}1r
\putmorphism(-200,-250)(1,0)[(\tilde B, \tilde A)` (\tilde B, \tilde A)`(\tilde B, 1_{\tilde A})]{480}1b
\putmorphism(360,-250)(1,0)[\phantom{Y}`(\tilde B,\tilde A) `(1_{\tilde B},\tilde A)]{480}1b
\put(400,-110){\fbox{$(-, \tilde A)_{\tilde B}$}}
\efig}
$$
which is true by: 1) vertical naturality of $(B,-)_A$ and $((u,1_A))$, and 2) by $((1_B,U))$ and vertical naturality of $(-,\tilde A)_B$.

\bigskip 

For 
(lx.f.h1) one takes two vertically composable pairs of 2-cells: 
$$
\scalebox{0.86}{
\bfig
\putmorphism(-150,280)(1,0)[A` \tilde A`K]{420}1a
\putmorphism(-150,-50)(1,0)[A'`\tilde A' `K']{410}1a
\putmorphism(-150,-360)(1,0)[\phantom{Y}` \tilde A''`fK'']{420}1a
\putmorphism(-170,280)(0,-1)[\phantom{Y_2}``U]{330}1l
\putmorphism(250,280)(0,-1)[\phantom{Y_2}``V]{330}1r
\putmorphism(-170,-60)(0,-1)[\phantom{Y_2}`A''`U']{300}1l
\putmorphism(250,-60)(0,-1)[\phantom{Y_2}``V']{300}1r
\put(-20,150){\fbox{$\alpha$}}
\put(-20,-180){\fbox{$\alpha'$}}
\efig}
\quad\text{and}\quad
\scalebox{0.86}{
\bfig
\putmorphism(-150,280)(1,0)[B` \tilde B`k]{420}1a
\putmorphism(-150,-50)(1,0)[B'`\tilde B' `k']{410}1a
\putmorphism(-150,-360)(1,0)[\phantom{Y}` \tilde B''`k'']{420}1a
\putmorphism(-170,280)(0,-1)[\phantom{Y_2}``u]{330}1l
\putmorphism(250,280)(0,-1)[\phantom{Y_2}``v]{330}1r
\putmorphism(-170,-60)(0,-1)[\phantom{Y_2}`B''`u']{300}1l
\putmorphism(250,-60)(0,-1)[\phantom{Y_2}``v']{300}1r
\put(-20,150){\fbox{$\beta$}}
\put(-20,-180){\fbox{$\beta'$}}
\efig}
$$
and one immediately sees that $P(\frac{\alpha}{\alpha'},\frac{\beta}{\beta'})=\frac{P(\alpha,\beta)}{P(\alpha', \beta')}$ by ($(u,U)$-l-nat) and ($(u,U)$-r-nat). 
Thus we obtain functoriality of $P$. 

\begin{rem}
The requirement that the 2-cells $(u,U)$ be identities is needed in the above proof to show that $P$ is strictly functorial 
with respect to 1v- and 2-cells. To prove the naturality of $\gamma$, it is sufficient to assume merely invertibility of $(u,U)$. 
\end{rem}


This finishes the proof that we have a lax double functor $P: \Aa\times\Bb\to\Cc$.


\medskip

Observe that  
if $\iota^A$ and $\iota^B$ are invertible, then so is $\iota^P$ and also $\gamma_{(1_{A},g)(f,1_{B'})}$. This is exactly the same 
as in \cite[Lemma 5.2]{FMS}. When $\iota$'s are invertible the lax double functor in question is called {\em unitary}, whereas the 
lax double functor $P:\Aa\times\Bb\to\Cc$ is called {\em decomposable} when $\gamma_{(1_{A},g)(f,1_{B'})}$ is invertible in {\em loc.cit.}.

We will prove later that the full double subcategory $q\x\Lax_{hop}^{st-u}(\Aa\times\Bb,\Cc)$ of $q\x\Lax_{hop}^{st}(\Aa\times\Bb,\Cc)$ consisting 
of {\em unitary} lax double quasi-functors (in the sense that both $(-,A)$ and $(B,-)$ are unitary) is double equivalent to the full double subcategory 
$\Lax_{hop}^{ud}(\Aa\times\Bb,\Cc)$ of $\Lax_{hop}(\Aa\times\Bb,\Cc)$ consisting of the {\em unitary and decomposable} lax double functors. 
Let 
\begin{equation} \eqlabel{2-functor F'}
\F': q\x\Lax_{hop}^{st-u}(\Aa\times\Bb,\Cc) \to \Lax_{hop}^{ud}(\Aa\times\Bb,\Cc).
\end{equation}
denote the corresponding restriction of $\F$ from \equref{2-functor F}, and 
let $\G$ denote a to-be-defined quasi-inverse double functor for $\F'$.

\subsection{A quasi-inverse $\G$ of $\F'$ and transformations $\kappa, \lambda$ on 0-cells} \sslabel{G on 0}

We will first show that a {\em unitary and decomposable} lax double functor $P:\Aa\times\Bb\to\Cc$ with the structures $\gamma$ and $\iota^P$ 
determines a lax double quasi-functor 
$H: \Aa\times\Bb\to\Cc$. Let $P(A,k):=P(1_A,k), P(A,u):=P(1^A,u), P(A,\beta):=P(\I_A,\beta)$, where $\I_A$ is the identity 2-cell 
for the object $A$, and $k,u,\beta$ as usual, and similarly for $P(-,B)$. It follows: $P(1_A,1_B)=P(A,1_B)=P(1_A,B)$. 

Now set $(-, A)=P(A,-)$ and $(B,-)=P(-, B)$ and $\iota^A_B=\iota^B_A:=\iota^P_{(A,B)}$. Then $(k', A)(k, A)=P(1_A,k')P(1_A, k)$ and similarly 
for $(B,-)$, and we may set: 
$$(-,A)_{k'k}:=\gamma_{(1_A,k')(1_A,k)}\quad\text{and}\quad (B,-)_{K'K}:= \gamma_{(K',1_B)(K,1_B)}$$
and we get that $(-, A)$ and $(B,-)$ are unitary lax double functors.

Observe the form of the 2-cell: 
$$
\scalebox{0.86}{
\bfig
 \putmorphism(200,400)(1,0)[``P(K,k)]{480}1a
\putmorphism(120,70)(1,0)[\phantom{F(A)}` `P(\tilde K,\tilde k)]{560}1b

\putmorphism(220,440)(0,-1)[\phantom{Y_2}``P(U,u)]{400}1l
\putmorphism(680,440)(0,-1)[\phantom{Y_2}``P(V,v)]{400}1r
\put(270,220){\fbox{$P(\alpha,\beta)$}}
\efig}
$$
and notice then whenever either of the two 1h-cells or either of the two 1v-cells in $P(-,-)$ above is identity, the form of 
$P(-,-)$ becomes $(B,x)$ {\em i.e.} $(x,A)$ for the corresponding 1h- or 1v-cell $x$. Then we may further define: 
$$(k,K):=
\scalebox{0.86}{
\bfig
\putmorphism(-190,50)(1,0)[` `(k,A)]{420}1a
\putmorphism(150,50)(1,0)[\phantom{F(A)}` `(B',f)]{450}1a
\putmorphism(-180,50)(0,-1)[\phantom{Y_2}``=]{350}1l
\putmorphism(610,50)(0,-1)[\phantom{Y_2}``=]{350}1r
\put(0,-120){\fbox{$\frac{\gamma_{(K,1_{B'})(1_A,k)}}{\gamma^{-1}_{(1_{A'},k)(K,1_B)}}$}}
 \putmorphism(-190,-290)(1,0)[``(B,K)]{420}1b
 \putmorphism(150,-290)(1,0)[\phantom{F(A)}` `(k,A')]{450}1b
\efig}
$$

$$(u,K):=
\scalebox{0.86}{
\bfig
 \putmorphism(200,200)(1,0)[``(B,K)]{700}1a
\putmorphism(120,-140)(1,0)[\phantom{F(A)}` `(\tilde B,K)]{760}1b

\putmorphism(220,240)(0,-1)[\phantom{Y_2}``(u,A)]{420}1l
\putmorphism(860,240)(0,-1)[\phantom{Y_2}``(u,A')]{420}1r
\put(270,20){\fbox{$P(\Id_K,\Id^u)$}}
\efig}
\qquad\text{and}\qquad
(k,U):=
\scalebox{0.86}{
\bfig
 \putmorphism(200,200)(1,0)[``(k,A)]{700}1a
\putmorphism(120,-140)(1,0)[\phantom{F(A)}` `(k,\tilde A)]{760}1b

\putmorphism(220,240)(0,-1)[\phantom{Y_2}``(B,U)]{420}1l
\putmorphism(880,240)(0,-1)[\phantom{Y_2}``(B',U)]{420}1r
\put(270,20){\fbox{$P(\Id^U,\Id_k)$}}
\efig}
$$
Since $P$ is strict in the vertical direction, we have: 
$\frac{(B,U)}{(u,\tilde A)}=\frac{(u,A)}{(\tilde B,U)}$, so we may define a 2-cell $(u,U)$ (in the desired form) to be identity. 
For the same reason the rules $((\frac{u}{u'},K))$, $((k,\frac{U}{U'}))$ and ($(u,U)$-l-nat), ($(u,U)$-r-nat) hold. The rules $((u,K'K))$ 
and $((k'k,U))$ hold by laxity of $P$. 

Since $P$ as a lax functor when evaluated at an identity 2-cell equals identity, the following rules hold true: 
$((1^B,K))$, $((u,1_A))$, $((1_B,U))$ and $((k,1^A))$. The rules $((1^B,U))$ and $((u,1^A))$ hold since we defined $(u,U)$ to be identity. 

Observe that the naturality of $\gamma$ with respect to 2-cells $\Id_K$ from $\Aa$ and $\beta$ from $\Bb$ reads: 
$$
\scalebox{0.86}{
\bfig
 \putmorphism(-150,500)(1,0)[(B,A)`(B', A)`(k, A)]{600}1a
 \putmorphism(450,500)(1,0)[\phantom{A\ot B}`(B', A') `(B', K)]{600}1a
 \putmorphism(-150,50)(1,0)[P(A,B)`P(A',B')`P(K,k)]{1200}1a

\putmorphism(-180,500)(0,-1)[\phantom{Y_2}``=]{450}1r
\putmorphism(1100,500)(0,-1)[\phantom{Y_2}``=]{450}1r
\put(250,280){\fbox{$\gamma_{(K,1_{B'})(1_A,k)}$}}

\putmorphism(-180,50)(0,-1)[\phantom{Y_2}``P(1^A,u)]{450}1l
\putmorphism(1100,50)(0,-1)[\phantom{Y_2}``P(1^{A'},v)]{450}1r
\put(280,-160){\fbox{$P(\Id_K,\beta)$}}

\putmorphism(-150,-400)(1,0)[P(A,\tilde B)`P(A',\tilde B') `P(K,\tilde k)]{1200}1a
\efig}
=
\scalebox{0.86}{
\bfig
 \putmorphism(-150,500)(1,0)[(B,A)`(B', A)`(k, A)]{600}1a
 \putmorphism(450,500)(1,0)[\phantom{A\ot B}`(B', A') `(B', K)]{600}1a

 \putmorphism(-170,50)(1,0)[(\tilde B,A)`(\tilde B', A)`(\tilde k, A)]{630}1a
 \putmorphism(450,50)(1,0)[\phantom{A\ot B}`(\tilde B', A') `(\tilde B',K)]{680}1a

\putmorphism(-200,500)(0,-1)[\phantom{Y_2}``P(1^A,u)]{450}1l
\putmorphism(450,500)(0,-1)[\phantom{Y_2}``]{450}1l
\putmorphism(630,500)(0,-1)[\phantom{Y_2}``P(1^A,v)]{450}0l 
\putmorphism(1120,500)(0,-1)[\phantom{Y_3}``P(1^{A'},v)]{450}1r
\put(-180,270){\fbox{$P(\I_A,\beta)$}} 
\put(600,270){\fbox{$P(\Id_K,\Id^v)$}}

\putmorphism(-180,50)(0,-1)[\phantom{Y_2}``=]{450}1l
\putmorphism(1120,50)(0,-1)[\phantom{Y_3}``=]{450}1r
\put(220,-160){\fbox{$\gamma_{(K,1_{\tilde B'})(1_A,\tilde k)}$}}

\putmorphism(-150,-400)(1,0)[P(A,\tilde B)`P(A',\tilde B') `P(K,\tilde k)]{1280}1a
\efig}
$$
Then applying this naturality of $\gamma_{(K,1_{\tilde B'})(1_A,\tilde k)}$ and analogously of $\gamma^{-1}_{(1_{A'},k)(K,1_B)}$ 
one obtains that the rule ($(k,K)$-l-nat) holds. Analogously, the naturality of $\gamma$ with respect to 2-cells $\alpha$ from $\Aa$ and $\Id_k$ from $\Bb$ 
are used to prove the rule ($(k,K)$-r-nat). 

All the resting rules from \prref{char df} concern only the horizontal structures and are already shown to hold in \cite[Theorem 5.3]{FMS}. 
We conclude that the unitary and decomposable lax double functor $P:\Aa\times\Bb\to\Cc$ indeed determines a unitary lax double 
quasi-functor $H: \Aa\times\Bb\to\Cc$. Thus $\G$, with $\G(P)=H$, is well-defined on 0-cells.

\bigskip 

Let us now start to define two horizontal strict transformations: $\kappa: \Id\Rightarrow\G\F'$ and $\lambda: \F'\G\Rightarrow\Id$. 
They will be isomorphisms if we restrict to unitary double lax quasi-functors in $q\x\Lax_{hop}(\Aa\times\Bb,\Cc)$ and to 
unitary and decomposable lax double functors in $\Lax_{hop}^{st}(\Aa\times\Bb,\Cc)$. This will yield that $\F'$ in \equref{2-functor F'} 
is a double equivalence, as announced. 

\medskip

Let $H$ be a lax double quasi-functor and $\F'(H)=P$. 
Observe that $\G\F'(H\vert_{(B,K)})=P(K,B)=(B,K)(1_B,A')$, and similarly $\G\F'(H\vert_{(k,A)})=P(A,k)=(B,1_A)(k,A)$. For 
a fixed $A\in\Aa$ we proceed to define a horizontal oplax transformation $\chi^A$ (and similarly and independently for a fixed $B\in\Bb$ 
a horizontal oplax transformation $\chi^B$). 
We set $\chi^A(B)=\chi^A_B:=1_{(B,A)}$; for a globular 2-cell $\chi^A_k: \chi^A(B')(k,A)\Rightarrow\G\F(k,A)\chi^A(B)$ which is actually 
$\chi^A_k: (k,A)\Rightarrow(B,1_A)(k,A)$, we set $\chi^A_k:=[\iota^B_A\vert\Id_{(k,A)}]$, and $(\chi^A)^u:=(u,1^A)=\Id^{(u,A)}$ by $((u,1^A))$, 
with notations as usual. We obtain that $(-,A)$ is naturally isomorphic to $P(A,-)$ if $\iota^B_A$ is 
invertible, through the horizontal oplax transformation $\chi^A$. To prove the property (h.o.t.-5) of \deref{hor nat tr} 
in this double categorical setting the interchange law is used. Analogously, one proves a natural isomorphism 
$(B,-)\iso P(-,B)$ through $\chi^B$, if $\iota^A_B$ is invertible. 
Thus when $H$ is unitary ({\em i.e.} $\iota^A$ and $\iota^B$ are invertible), $\chi^A$ and $\chi^B$ are isomorphisms. 
It is easily seen that $\chi^A_B=\chi^B_A$.

Observe further that $\G\F'(k,K)$ is a 2-cell whose both source and target 1h-cells are compositions of four 1h-cells, and not of two 1h-cells 
as in the case of $(k,K)$. To express $\G\F'(k,K)$, one uses that $\gamma^{-1}_{(1_{A'},k)(K,1_B)}=
[\Id_{(B,K)}\vert\iota^{A'}_B\vert\iota^B_{A'}\vert\Id_{(k,A')}]$ 
by \cite[Lemma 5.2]{FMS}. Moreover, $\G\F'(u,K)$ and $\G\F'(k,U)$ are 2-cells whose source and target 1h-cells are compositions of two 1h-cells, 
and not a single 1h-cell as in the case of $(u,K)$ and $(k,U)$, respectively. It is easily seen and it is proved in \cite[Theorem 5.3]{FMS} 
that $\chi^A$ and $\chi^B$ obey $(HOT^q_1)$.  
The axioms $(HOT^q_2)-(HOT^q_4)$ for $\chi^A$ and $\chi^B$ hold almost trivially: $\chi^A_u, \chi^B_U$ are identities 
and in properties $(HOT^q_2)$ and $(HOT^q_3)$ use the interchange law to move the 2-cells $\iota^\bullet_\bullet$.  

In this way we have defined a 0-component $\kappa^H$ at a 0-cell $H$ in $q\x\Lax_{hop}^{st}(\Aa\times\Bb,\Cc)$ 
of $\kappa: \Id\Rightarrow\G\F'$, 
which will be a horizontal strict transformation. 

To define 0-component of a horizontal strict transformation $\lambda: \F'\G\Rightarrow\Id$, we see that $\F'\G(P)(A,B)=P(A,B)$ for a 
unitary decomposable double lax functor $P:\Aa\times\Bb\to\Cc$, so we may set $\lambda^P$ to be identity between 0-cells. 
Similarly, as $\F'\G(P)(U,u)=\frac{(B,U)}{(u,\tilde A)}=\frac{P(1^B,U)}{P(u, 1^{\tilde A})}=P(U,u)$ we may set $\lambda^P$ to be identity 
also between 1v-cells. Though, $\F'\G(P)(K,k)=(k,A')(B,K)=P(1_{A'}, k)P(K,1_B)$, then we set $\lambda^P$ on a 1h-cell $(K,k)$ to be 
$\gamma_{(1_{A},k)(K,1_{B'})}: P(1_{A'}, k)P(K,1_B) \Rightarrow P(K,k)$. Such defined $\lambda^P$ is indeed a horizontal 
oplax transformation of double lax functors: properties (h.o.t.-1) and (h.o.t.-2) are proved in \cite[Theorem 5.3]{FMS}, properties (h.o.t.-3) and 
(h.o.t.-4) hold since 1v-components of $\lambda^P$ are identities, and (h.o.t.-5) holds by naturality of $\gamma$. 

In the next two Subsections we will finalize the proof that $\kappa$ and $\lambda$ are horizontal strict transformations. 
Observe from above 
that restricting to the full double subcategories $q\x\Lax_{hop}^{st-u}(\Aa\times\Bb,\Cc)$ and $\Lax_{hop}^{ud}(\Aa\times\Bb,\Cc)$ we indeed 
obtain a double equivalence. Namely, in these double subcategories $\kappa$ is an isomorphism since so are $\chi^A$'s and $\chi^B$'s, and 
$\lambda$ is an isomorphism since the 0-component of $\lambda$ is defined to be $\gamma_{(1_{A},k)(K,1_{B'})}$ on 1h-cells $(K,k)$.

\subsection{$(\F,\G,\kappa,\lambda)$ on 1h- and 1v-cells}

We first give the definitions of $\F$ and $\G$ on 1h-cells. 

To define $\F$ on 1h-cells, 
let a horizontal oplax transformation between lax double quasi-functors $(-,-)_1, (-,-)_2$ with images $P,P'$ be given via a pair of families  
$\theta^A, \theta^B, A\in\Aa, B\in\Bb$ of horizontal oplax transformations. We define $\Sigma: P\Rightarrow P'$ by $\Sigma(A,B):=\theta^A_B$, 
$\Sigma_{(U,u)}:=\frac{(\theta^B)^U}{(\theta^{\tilde A})^u}$ and $\delta_{\Sigma,(K,k)}=\Sigma_{(K,k)}:=
\frac{[\Id_{(B,K)_1}\vert\theta^{A'}_k]}{[\theta^B_K\vert\Id_{(k,A')_2}]}$. The properties (h.o.t.-1) and (h.o.t.-2) of \deref{hor nat tr} for 
$\Sigma$ are proved in \cite[Theorem 5.3]{FMS}. The properties (h.o.t.-3) and (h.o.t.-4) follow by the same properties for $\theta^B_{\frac{U}{V}}$ and 
$\theta^{\tilde{\tilde A}}_{\frac{u}{v}}$, $(HOT^q_4)$ and since by assumption 
the vertically globular 2-cells $(u,V)_1, (u,V)_2$ are identities. (This includes the proof for $\Sigma_{(1^A,1^B)}=\Id_{\Sigma(A,B)}$.)

To prove the property (h.o.t.-5) of \deref{hor nat tr} for $\Sigma$ one uses: property ($(u,U)$-l-nat) of $(-,-)_1$ and property (h.o.t.-5) for $\theta^{\tilde A'}$, 
then simultaneously $(HOT^q_2)$ and $(HOT^q_3)$, and finally property (h.o.t.-5) for $\theta^B$ and ($(u,U)$-l-nat) of $(-,-)_2$. Then $\Sigma: P\Rightarrow P'$ 
is indeed a horizontal oplax transformation of lax double functors. 

By inspecting $\F(\frac{\theta}{\theta'})$ and $\frac{\F(\theta)}{\F(\theta')}$ on $(A,B), (U,u)$ and $(K,k)$ and the action of $\F$ on the idenity 
one sees that $\F$ is a strict double functor. 

For $\G$, let $P,P':\Aa\times\Bb\to\Cc$ be unitary and decomposable lax double functors and $\Sigma: P\Rightarrow P'$ a horizontal oplax transformation 
between them. Then $\G(\Sigma)$ is a horizontal oplax transformation between lax double quasi-functors given by the families of 
$\theta^A:=\Sigma(A,-), \theta^B:=\Sigma(-,B)$ for 
$A\in\Aa, B\in\Bb$. It's clearly $\theta^A_B=\theta^B_A$ and the condition $(HOT^q_1)$ is proved to hold in \cite[Theorem 5.3]{FMS}. 
Conditions $(HOT^q_2)$ and $(HOT^q_3)$ hold by the property (h.o.t.-5) of $\Sigma$ from \deref{hor nat tr} with $a=(\Id_K,\Id^u)$ and 
$a=(\Id^U,\Id_k)$, respectively, while $(HOT^q_4)$ holds by the property (h.o.t.-3) thereof: both sides in $(HOT^q_4)$ equal $\Sigma_{(U,u)}$. 
It is directly computed that $\G$ is a strict double functor. 

\bigskip 

We now define $\F$ and $\G$ on 1v-cells.  

Given a vertical lax transformation $\theta_0$ between lax double quasi-functors $(-,-)_1, (-,-)_3$ with images $P,\tilde P$, with a pair of families  
$\theta_0^A, \theta_0^B, A\in\Aa, B\in\Bb$ of vertical lax transformations, we define $\Sigma_0: P\Rightarrow\tilde P$ by $\Sigma_0(A,B):=(\theta_0^A)_B
=(\theta_0^B)_A, \Sigma_0^{(U,u)}=[\frac{(\theta_0^B)^U}{\Id^{\tilde P(\tilde A,u)}}\vert\frac{\Id^{P(U,B)}}{(\theta_0^{\tilde A})^u}]$ 
and $(\Sigma_0)_{(K,k)}:=[(\theta_0^B)_K\vert (\theta_0^{A'})_k]$. Analogously as in the case for 1h-cells it is proved that $\Sigma_0: P\Rightarrow\tilde P$  
is a vertical lax transformation between lax double functors.

For $P,\tilde P:\Aa\times\Bb\to\Cc$ and $\Sigma_0: P\Rightarrow \tilde P$ a vertical lax transformation 
we define $\G(\Sigma_0)$ as a vertical lax transformation given by the families:  
$\theta^A_0:=\Sigma_0(A,-), \theta_0^B:=\Sigma_0(-,B)$ for $A\in\Aa, B\in\Bb$. The proof is analogous as in the case for 1h-cells and as before we have 
$\G\F'( (\theta^A_0)_B, (\theta_0^B)_A)_\bullet = ((\theta^A_0)_B, (\theta_0^B)_A)_\bullet$, where $\bullet$ stands for the indexing over 
$A\in\Aa, B\in\Bb$.

\bigskip

Let us now study $\kappa: \Id\Rightarrow\G\F'$ at a 1h-cell component, a horizontal oplax transformation between lax double quasi-functors 
given by a family $(\theta^A,\theta^B), A\in\Aa, B\in\Bb$. To prove that $\kappa$ is a horizontal strict transformation, 
on one hand, we should show the identity $\kappa^{H'}\comp (\theta^A,\theta^B)=
\G\F'(\theta^A,\theta^B)\comp\kappa^H$. (Then the axioms (h.o.t.-1) \label{(h.o.t.-1)} and (h.o.t.-2) \label{(h.o.t.-2)} hold trivially.) 
This means that both $\chi^A\comp\theta^A=\G\F'(\theta^A)\comp\chi^A$ (with abuse of notation 
by writing $\G\F'(\theta^A)$ which is easily understood) and the analogous identity for $B$ must hold. We check only the first identity. 
At a 0-component $B$ we have that $\chi^A_B$ is identity and it is easily seen that $\theta^A_B=\G\F'(\theta^A_B)$. At a 1h-component $k$: 
$\G\F'(\theta^A_k)=\frac{[\Id_{(B,1_A)}\vert\theta^{A}_k]}{[\theta^B_{1_A}\vert\Id_{(k,A)}]}$ and observe that by the property (h.o.t.-4) of 
\deref{hor nat tr} we have that $\theta^B_{1^A}$ is identity. Recall that the $\chi^A_B$'s are identities by construction, so 
the compositions in the identity $\chi^A_k\comp\theta^A_k=\G\F'(\theta^A_k)\comp\chi^A_k$ make sense, and the identity is shown to 
hold by the interchange law. 
Finally, at a 1v-component $u$ we have that $\chi^A_u$ is identity, on one hand, and observe that $\G\F'(\theta^A_u)=
\frac{\theta^B_{1_A}}{\theta^A_u}$, on the other. But $\theta^B_{1^A}$ is identity, so we get indeed $\chi^A_u\comp\theta^A_u=
\G\F'(\theta^A_u)\comp\chi^A_u$, as desired. 

On the other hand, we should define $\kappa: \Id\Rightarrow\G\F'$ at a 1v-cell component and check if the axioms 
(h.o.t.-3) \label{(h.o.t.-3)} - (h.o.t.-5) \label{(h.o.t.-5)} hold. 
To define $\kappa: \Id\Rightarrow\G\F'$ at a 1v-cell component $\theta_0=(\theta^A_0, \theta_0^B)_\bullet$ we should 
define a modification among four lax double quasi-functors (on the left below), which is given by two modifications of four lax double functors, as shown                                                      :
$$
\scalebox{0.86}{
\bfig
\putmorphism(-180,50)(1,0)[H` \G\F'(H)` \kappa^H]{450}1a
\putmorphism(-180,-270)(1,0)[\tilde H ` \G\F'(\tilde{H}) ` \kappa^{\tilde H}]{450}1b
\putmorphism(-170,50)(0,-1)[\phantom{Y_2}``\theta_0]{320}1l
\putmorphism(250,50)(0,-1)[\phantom{Y_2}``\G\F'(\theta_0)]{320}1r
\put(-60,-140){\fbox{$\kappa^{\theta_0}$}}
\efig}
\qquad\qquad
\scalebox{0.86}{
\bfig
\putmorphism(-180,50)(1,0)[(-,A)` (-,A)`\chi^A]{550}1a
\putmorphism(-180,-270)(1,0)[(-,A)_3`(-,A)_3 `\tilde\chi^A]{550}1b
\putmorphism(-170,50)(0,-1)[\phantom{Y_2}``\theta_0^A]{320}1l
\putmorphism(350,50)(0,-1)[\phantom{Y_2}``\G\F'(\theta_0^A)]{320}1r
\put(0,-140){\fbox{$\kappa^{\theta_0^A}$}}
\efig}
\qquad\qquad
\scalebox{0.86}{
\bfig
\putmorphism(-180,50)(1,0)[(B,-)` (B,-)`\chi^B]{550}1a
\putmorphism(-180,-270)(1,0)[(B,-)_3`(B,-)_3 `\tilde\chi^B]{550}1b
\putmorphism(-170,50)(0,-1)[\phantom{Y_2}``\theta_0^B]{320}1l
\putmorphism(350,50)(0,-1)[\phantom{Y_2}``\G\F'(\theta_0^B)]{320}1r
\put(0,-140){\fbox{$\kappa^{\theta_0^B}$}}
\efig}
$$
where $(\kappa^{\theta_0^A})_B=(\kappa^{\theta_0^B})_A$ for every $A\in\Aa, B\in\Bb$. When evaluating the latter two to-be-defined 
modifications at 0-cells $B\in\Bb$, respectively $A\in\Aa$, we see that by construction the 1h-cells appearing in the obtained 2-cells 
$(\kappa^{\theta_0^A})_B=(\kappa^{\theta_0^B})_A$ are identities: $\chi^A_B=\chi^B_A=1_{(B,A)}$ and $\tilde\chi^A_B=\tilde\chi^B_A=1_{(B,A)_3}$. 
Moreover, we also have that $\G\F'( (\theta^A_0)_B, (\theta_0^B)_A)_\bullet = ((\theta^A_0)_B, (\theta_0^B)_A)_\bullet$. Then we may set 
$(\kappa^{\theta_0^A})_B:=\Id^{(\theta^A_0)_B}$ and $(\kappa^{\theta_0^B})_A:=\Id^{(\theta^B_0)_A}$, which actually are the same 2-cell for every 
$A\in\Aa, B\in\Bb$. 

We next verify if this 2-cell obeys the axioms (m.ho-vl.-1) \label{(m.ho-vl.-1)} and (m.ho-vl.-2) \label{(m.ho-vl.-2)}. As we showed above that 
$\kappa$ at a 1h-cell component gives an identity (globular) 2-cell (recall the fact responsible for $\kappa$ being a horizontal strict transformation), 
in the axiom (m.ho-vl.-1) the 2-cells corresponding to $\delta_{\alpha,f}$ and $\Theta_A$ (and similarly for $B$) there are now identities, 
so the axiom trivially holds. 
As for the axiom (m.ho-vl.-2)\label{(m.ho-vl.-2)}, the 2-cells corresponding to $\alpha^u$ and $\beta^u$ there are now $(\chi^A)^u$ and 
$(\tilde\chi^A)^u$, which are by definition $\Id^{(u,A)}$ and $\Id^{(u,A)_3}$, respectively. Then we again have that the axiom holds. 

This terminates the definition of $\kappa$ at a 1v-cell component $\theta_0$. 

\medskip

It rests to check that the axioms (h.o.t.-3) \label{(h.o.t.-3)} - (h.o.t.-5) \label{(h.o.t.-5)} hold for $\kappa$ as a horizontal strict transformation. 
As  $\kappa$ at a 1v-cell component is given by identity 2-cells in both variables, $(\kappa^{\theta_0^A})_B=(\kappa^{\theta_0^B})_A=\Id^{(\theta^A_0)_B}$, 
and as we saw further above $\kappa$ at a 1h-cell component is also given by identity 2-cells in both variables, these three resting axioms are trivially 
fulfilled (although we still haven't defined $\F$ and $\G'$ on 2-cells, {\em i.e.} on the respective modifications). 

\bigskip

To prove that $\lambda$ is a horizontal strict transformation, to the proof in \cite[Theorem 5.3]{FMS} we need to add: 1) the check that  
$\lambda^{P'}_{(U,u)}\comp\F'\G(\Sigma_{(U,u)})=\Sigma_{(U,u)}\comp\lambda^P_{(U,u)}$, 
for an oplax transformation of double lax functors $\Sigma:P\Rightarrow P'$, 
2) define a modification $\lambda^{\Sigma_0}$ corresponding to $\lambda$ at a 1v-vell component $\Sigma_0$, and 3) verify the 
axioms (h.o.t.-3) \label{(h.o.t.-3)} - (h.o.t.-5) \label{(h.o.t.-5)} for $\lambda$ to be a horizontal strict transformation. 
We leave the points 2) and 3) to the reader. To finish 1), recall that the 1v-components of $\lambda^P$ are identities, so 
it remains to check if $\F'\G(\Sigma_{(U,u)})=\Sigma_{(U,u)}$ holds. We find: $\F'\G(\Sigma_{(U,u)})=
\frac{\G(\Sigma_{(U,u)})^B_U}{\G(\Sigma_{(U,u)})^{\tilde A}_u}=\frac{\Sigma(U,1^B)}{\Sigma(1^{\tilde A},u)}=\Sigma(U,u)$, the latter 
identity holds by the property (h.o.t.-3) of $\Sigma$ being a horizontal oplax transformation of double lax functors. Thus we proved the desired equality.

\subsection{$\F$ and $\G$ on 2-cells}

We start by defining $\F$ on modifications. Let a modification $\tau=(\tau^A, \tau^B)_{A\in\Aa, B\in\Bb}$ in $q\x\Lax_{hop}^{st}(\Aa\times\Bb,\Cc)$ 
be given, recall \equref{q-modif}. We define $\F(\tau)$ by $\F(\tau)_{(A,B)}:=\tau^A_B=\tau^B_A$. It is directly checked that this is a modification 
in $\Lax_{hop}(\Aa\times\Bb,\Cc)$. 

Conversely, given a modification $\Theta$ in $\Lax_{hop}(\Aa\times\Bb,\Cc)$ between horizontal oplax and vertical lax transformations of lax double functors, 
define $\tau^A:=\Theta(A,-)$ and $\tau^B:=\Theta(-,B)$. It is directly seen that they give modifications in the sense of \deref{modif-hv}, 
and it is clearly $\tau^A_B=\tau^B_A$, so we obtain a modification $\G(\Theta)=(\tau^A, \tau^B)_{A\in\Aa, B\in\Bb}$ of horizontal oplax 
and vertical lax transformations of lax double quasi-functors.



\bigskip

To summarize, in this Section we have proved the following results:

\begin{prop} \prlabel{F}
With notations as at the beginning of \seref{isom-spec} there is a double functor 
$$\F: q\x\Lax_{hop}^{st}(\Aa\times\Bb,\Cc) \to \Lax_{hop}(\Aa\times\Bb,\Cc).$$
\end{prop}

\begin{thm} \thlabel{equiv-fun}
With notations as explained above \equref{2-functor F'}, the double functor $\F$ restricts to a double equivalence functor 
$$\F': q\x\Lax_{hop}^{st-u}(\Aa\times\Bb,\Cc) \to \Lax_{hop}^{ud}(\Aa\times\Bb,\Cc)$$
with quasi-inverse $\G$. 
\end{thm}

This Theorem is a double categorical version of \cite[Theorem 5.3]{FMS}. We can straighten a bit its formulation by passing to pseudo 
(quasi-) functors, as they are lax unitary (quasi-) functors with all $\gamma$'s invertible. The notation $\Lax$ changes then to $\Ps$, 
the supraindex $u$ becomes superfluous, but also $d$ in the right hand-side. 
Moreover, observe that in the 0-cells of $q\x\Ps_{hop}^{st}(\Aa\times\Bb,\Cc)$ in the left, the 2-cells $(k,K)$ of quasi pseudofunctors 
are invertible. Then the double equivalence functor $\F'$ restricts to a double equivalence 
$$\F'': q\x\Ps_{hop}^{st}(\Aa\times\Bb,\Cc) \to \Ps_{hop}(\Aa\times\Bb,\Cc).$$
Observe that by \cite[Proposition 6.2]{FMS} if the 2-cells $(k,K)$ of a lax unitary quasi-functor are invertible, as it is the case in pseudo quasi-functors, 
then the properties $((1_B,K))$ and $((k,1_A))$ in \prref{char df} are redundant. Namely, in the underlying horizontal 2-category of $\Aa\ot\Bb$ one can perform the computation carried out in the proof of \cite[Proposition 6.2]{FMS} and pull the result back to the double category. 
It comes down to a series of ``tricks'': add an identity 2-cell in the form of the unity axiom of the lax functor structure of $(-,A')$, 
add $\frac{(1_B,K)}{(1_B,K)^{-1}}=\Id$ between one lag of unity and multiplicativity of the functor, use the naturality of $(1_B,K)^{-1}$ with respect 
to the multiplicativity 2-cell $(-,A')_{1_B1_B}$ (this is $((k'k,K))$), use again the unity axiom of the lax functor structure and finally 
$\frac{(1_B,K)}{(1_B,K)^{-1}}=\Id$.

\section{Applications}

After proving our main results in Sections 3 and 4 we dedicate this last Section to some specific cases. 
We will also 
prove the universal property of $\Aa\ot\Bb$ and discuss monads in double categories.

\subsection{``(Un)currying'' functor}

At the beginning of \seref{isom-spec} we commented that the double category isomorphism \equref{pre-Gray} restricts to 
a double category isomorphism $q\x\Lax_{hop}^{st}(\Aa\times\Bb,\Cc)\iso\Lax_{hop}(\Aa, \llbracket\Bb,\Cc\rrbracket^{st})$. 
Composing this with $\F$ we obtain a double functor: 
\begin{equation} \eqlabel{uncurry}
\Lax_{hop}(\Aa, \llbracket\Bb,\Cc\rrbracket^{st})\to \Lax_{hop}(\Aa\times\Bb,\Cc)
\end{equation}
that is a double categorical version of the ``uncurrying'' double functor $J$ at the end of Section 4 of \cite{FMS}. ($J$ was implicitly constructed in \cite{Nik}.)

In \equref{2-functor F'} we moreover restricted to {\em unitary} lax double (quasi) functors. On the left hand-side therein (and in the last Theorem above) unitarity of a lax double quasi-functor $H$ refers to the unitarity of both $(-,A)$ and $(B,-)$ lax double functors 
comprising $H$. In the isomorphism $q\x\Lax_{hop}^{st}(\Aa\times\Bb,\Cc)\iso\Lax_{hop}(\Aa, \llbracket\Bb,\Cc\rrbracket^{st})$ 
unitarity of $(-,A)$ corresponds to the unitarity of 0-cells in $\llbracket\Bb,\Cc\rrbracket^{st}$, while unitarity of $(B,-)$ 
corresponds to the unitarity of 0-cells in $\Lax_{hop}(\Aa, \llbracket\Bb,\Cc\rrbracket^{st})$. Then 
the isomorphism $q\x\Lax_{hop}^{st}(\Aa\times\Bb,\Cc)\iso\Lax_{hop}(\Aa, \llbracket\Bb,\Cc\rrbracket^{st})$ restricts further to a 
double category isomorphism  
\begin{equation} \eqlabel{corch}
q\x\Lax_{hop}^{st-u}(\Aa\times\Bb,\Cc)\iso\Lax_{hop}^u(\Aa, \llbracket\Bb,\Cc\rrbracket^{st-u})
\end{equation}
where $\llbracket\Bb,\Cc\rrbracket^{st-u}$ denotes the double category of 0: unitary lax double functors $\Bb\to\Cc$, 
1h: horizontal oplax transformations, 1v: vertical strict transformations, and 2: modifications, 
and $\Lax_{hop}^u(\Aa, \llbracket\Bb,\Cc\rrbracket^{st-u})$ is the double category of 0: unitary lax double functors 
$\Aa\to\llbracket\Bb,\Cc\rrbracket^{st-u}$, 1h: horizontal oplax transformations between them, 1v: vertical strict transformations 
and modifications between the latter two. Joining \equref{corch} and \thref{equiv-fun} yields 
$$\Lax_{hop}^{ud}(\Aa\times\Bb,\Cc)\simeq\Lax_{hop}^u(\Aa, \llbracket\Bb,\Cc\rrbracket^{st-u})$$
which presents a setting in which the uncurrying double functor \equref{uncurry} restricts to a double category equivalence, {\em i.e.} 
in which a ``currying'' functor exists.

\subsection{The universal property of $\ot$} \sslabel{Gray}

At the end of \ssref{Gray prod} we announced a universal property of $\Aa\ot\Bb$ by which it strictifies lax double quasi-functors. 
We prove it here and upgrade it to an isomorphism of double categories.

\begin{prop}
There is an isomorphism of double categories 
$$q\x\Lax_{hop}(\Aa\times\Bb,\Cc)\iso \Dbl_{hop}(\Aa\ot\Bb,\Cc)$$ 
where the right hand-side is the double category of strict double functors, horizontal oplax transformations as 1h-cells, 
vertical lax transformations as 1v-cells, and modifications. 
\end{prop}

\begin{proof}
For $H\in q\x\Lax_{hop}(\Aa\times\Bb,\Cc)$ define $\crta{H}:\Aa\ot\Bb\to\Cc$ by $\crta{H}(A\ot y):=H(A,y)=(y,A)$ and $\crta{H}(x\ot B):=H(x,B)=(B,x)$ 
for all four types of cells $x$ in $\Aa$ and $y$ in $\Bb$, extend $\crta H$ to a strict double functor (in particular, $\crta H(1_{A\ot B})
=1_{H(A,B)}=1_{(B,A)}$) and define $\crta H((A\ot-)_B):=(-,A)_B, \crta H((A\ot-)_{k'k}):=(-,A)_{k'k}$, and similarly for the other entry, 
as well as for the 2-cells $K\ot k, K\ot u, U\ot k$ and $U\ot u$. 

Conversely, take $G\in\Dbl_{hop}(\Aa\ot\Bb,\Cc)$, define $\crta{(-,A)}:\Bb\to\Cc$ by $\crta{(y,A)}:=G(A\ot y)$, two globular 2-cells: 
$\crta{(-,A)}_{k'k}: G(A\ot k')G(A\ot k)=G((A\ot k')(A\ot k))\stackrel{G((A\ot-)_{k'k})}{\Rightarrow}G(A\ot k'k)$ and 
$\crta{(-,A)}_B:= G(A\ot -)_B$, and analogously for $\crta{(B,-)}:\Aa\to\Cc$. 
Then it is easily and directly proved that $\crta{(-,A)}$ and $\crta{(B,-)}$ are lax double functors. Define the 2-cells 
$\crta{(k,K)}, \crta{(u,K)}, \crta{(k,U)}$ and $\crta{(u,U)}$ in the obvious way, then the laws from \prref{char df} for 
$\crta{(-,A)}$ and $\crta{(B,-)}$ to make a lax double quasi-functor 
pass {\em mutatis mutandi} 
from the defining relations of $\Aa\ot\Bb$, since $G$ is a strict double functor.

Given a horizontal oplax transformation $\theta=(\theta^A, \theta^B)_{A\in\Aa, \\B\in\Bb}$ between lax double quasi-functors $H\Rightarrow H'$ 
we define a horizontal oplax transformation $\Sigma:\crta{H}\Rightarrow\crta{H'}$ by setting: 
$\Sigma(A\ot B)=\theta^A_B=\theta^B_A, \Sigma_{A\ot k}=\theta^A_k, \Sigma_{K\ot B}=\theta^B_K$ and $\Sigma_{A\ot u}=\theta^A_u, 
\Sigma_{U\ot B}=\theta^B_U$. 
To check the property (h.o.t.-5) of \deref{hor nat tr} for $\Sigma$ one should check it for ten types of 2-cells $a$ in $\Aa\ot\Bb$: for $a$ being 
$A\ot \omega$ or $\zeta\ot B$ the property (h.o.t.-5) for $\Sigma$ holds since $\theta^A$ respectively $\theta^B$ is a horizontal oplax transformation, 
while for $a$ being $K\ot k, K\ot u, U\ot k$ and $U\ot u$ the property (h.o.t.-5) for $\Sigma$ holds by the properties $HOT^q_1-HOT^q_4$, respectively, 
and for 2-cells of the type 
\equref{4 laxity 2-cells} it holds by $HOT^q_1$. The properties (h.o.t.-1)-(h.o.t.-4)  
of \deref{hor nat tr} for $\Sigma$ hold by the same properties for $\theta^A$ and $\theta^B$. For the converse, provided a 
horizontal oplax transformation of strict double functors $\tilde\Sigma:G\Rightarrow\crta G'$, 
define $\zeta=(\zeta^A, \zeta^B)_{A\in\Aa, \\B\in\Bb}$ in the obvious (converse) way. 

The definition and correspondence on vertical lax transformations is analogous as for horizontal oplax ones. 
Given a vertical lax transformation $\theta_0=(\theta_0^A, \theta_0^B)_{A\in\Aa, \\B\in\Bb}$ between lax double quasi-functors $H\Rightarrow H_0$ 
one constructs a vertical lax transformation $\Sigma_0:\crta{H}\Rightarrow\crta{H_0}$ on strict double functors $\Aa\ot\Bb\to\Cc$.

Given a modification $\tau=(\tau^A, \tau^B)_{A\in\Aa, \\B\in\Bb}$ as in \equref{q-modif} 
we define a modification $\aaa$: 
$$
\scalebox{0.86}{
\bfig
\putmorphism(-150,50)(1,0)[\crta{H}` \crta{H'}` \Sigma]{400}1a
\putmorphism(-150,-270)(1,0)[\crta{H_0} ` \crta{H_0'} ` \Sigma' ]{400}1b
\putmorphism(-170,50)(0,-1)[\phantom{Y_2}``\Sigma_0]{320}1l
\putmorphism(250,50)(0,-1)[\phantom{Y_2}``\Sigma_0']{320}1r
\put(-30,-140){\fbox{$\Theta$}}
\efig}
$$
by $\Theta(A\ot B)=\tau^A_B=\tau^B_A$. It is immediate to see that $\Theta$ is well-defined. For the converse, formulate the obvious (converse) definition.

On all the four levels of cells it is clear that one has a 1-1 correspondence, so that one obtains an isomorphism of double categories, as claimed. 
\qed\end{proof}

Joining the isomorphism from the above Proposition and \equref{pre-Gray} we obtain that there is an isomorphism of double categories: 
\begin{equation} \eqlabel{adj}
\Dbl_{hop}(\Aa\ot\Bb,\Cc)\iso \Lax_{hop}(\Aa, \llbracket\Bb,\Cc\rrbracket).
\end{equation}
This is 
a strictification result for lax double functors $\Aa\to \llbracket\Bb,\Cc\rrbracket$. 

Forgetting the vertical direction in the above double category isomorphism, 
{\em i.e.} restricting to the horizontal 2-categories of $\Aa,\Bb, \Cc$, we recover \cite[Proposition 2.9]{Nik} (more precisely (78) 
in Corollary 2.12 of {\em loc. cit.}, as we work with horizontal oplax transformations rather than their lax counterparts). 
Namely, the underlying horizontal 2-category of our tensor product $\Aa\ot\Bb$ is precisely the author's $\A\Del_{cmp}\B$ constructed in 
Section 2.8 for 2-categories 
$\A$ and $\B$ seen as the horizontal 2-categories of $\Aa$ and $\Bb$, respectively: $\HH(\Aa\ot\Bb)=\HH(\Aa)\Del_{cmp}\HH(\Bb)$. 

\begin{rem}
The reader may have noticed that the order of $\Aa$ and $\Bb$ in \equref{adj} is the same on both sides, whereas it appears swapped in 
(78) of \cite[Corollary 2.12]{Nik}. However, our result is in accordance with Gray's \cite[Theorem I.4.14]{Gray} for the oplax version 
of transformations, while the order in Proposition 2.9 and (143) in Section 4.1 of \cite{Nik} appears swapped with respect to Gray's 
\cite[Theorem I.4.9]{Gray} in the lax case. 
\end{rem}

By the 1-1 correspondence at the level of 0-cells in the double category isomorphism  \equref{adj} 
we conclude that there is an isomorphism of sets
$$
Dbl_{st}(\Aa\ot\Bb, \Cc)\iso Dbl_{lx}(\Aa, \llbracket\Bb,\Cc\rrbracket)
$$
where $Dbl_{st}$ is the category of double categories and strict double functors.

\subsection{Monads in double categories}

Many authors had observed that although various algebraic structures appear as monads in suitable bicategories, 
the corresponding morphisms are not morphisms of monads, considered as 1-cells in the bicategory $\Mnd(\B)$ of monads in a bicategory $\B$ 
from \cite{S}. 
A very well-known example is the bicategory $\B=Span(\C)$ of spans over a category $\C$ with pullbacks, introduced in \cite{Be}. 
It is immediate to see that monads in 
$Span(\C)$ are categories internal in $\C$. (As a matter of fact, in \cite[Section 5.4.3]{Be} categories internal in $\C$ are defined this way.) 
However, although monads in $Span(\C)$ are internal categories in $\C$, morphisms of monads in $Span(\C)$ are not morphisms in $\Cat(\C)$, 
the category of internal categories in $\C$. 

To remediate this inconsistency, in \cite[Example 2.1]{FGK} a pseudodouble 
category $\Span(\C)$ of spans in $\C$ was introduced whose horizontal bicategory is precisely the bicategory $Span(\C)$. Moreover, in 
\cite[Definition 2.4]{FGK} the authors introduced a pseudodouble category $\Mnd(\Dd)$ of monads in a pseudodouble category $\Dd$, so that when 
$\Dd=\Span(\C)$, the vertical 1-cells in $\Mnd(\Span(\C))$ are morphisms of internal categories in $\C$ (see \cite[Example 2.6]{FGK}). The 
construction of $\Mnd(\Dd)$ enhanced also other examples of the described inconsistency for bicategories $\B$ that could be upgraded into a 
double category $\Dd(\B)$. 

This explains why the authors defined a monad in a double category $\Dd$ as a monad in the horizontal 2-category $\HH(\Dd)$ of $\Dd$. On 
the other hand, B\'enabou observed in \cite{Be} that a lax functor $*\to\K$ from the trivial 2-category to a 2-category $\K$ is nothing but a 
monad in $\K$: the lax functor structure corresponds to the multiplication and the unit of the monad.  
It is straightforwardly seen that the analogous holds for monads in a double category $\Dd$: the only new thing now is that we have the identity 
1v-cell on the unique 0-cell, which is strictly preserved by a lax double functor, so no new data is obtained. 
Let now $*$ denote the trivial double category, then we may write:

\begin{prop} 
A lax double functor $*\to\Dd$ is a monad in $\Dd$. 
\end{prop}

Moreover, a 0-cell in $q\x\Lax_{hop}(*\times *,\Dd)$ is then given by two monads in $\Dd$, and the only surviving 2-cell 
(and laws) in the characterization \prref{char df} is the one of type $(k,K)=(id_*,id_*)$ and the rules $((1_B,K))$, $((k,1_A))$, 
$((k'k,K))$ and $((k,K'K))$, which correspond to monad-monad distributive laws. 

So far we have proved the 1-1 correspondence at the levels of 0-cells in the Proposition below. Since 1h-cells in $\Mnd(\Dd)$ correspond to 
horizontal {\em lax} transformations between lax double quasi-functors from the trivial double category, we are led to the double category $\Lax_{hop}^*(*,\Dd)$ 
from the end of \seref{definitions} (and correspondingly to the double category $q\x\Lax_{hop}^*(*\times *,\Dd)$, with the obvious meaning). To the axioms 
for cells in $\Lax_{hop}^*(*,\Dd)$ we will refer by the same labels as for those in $\Lax_{hop}^*(*,\Dd)$, to avoid introduction of additional notation. 
This said, observe that one can consider two distinct versions of the double category $\Mnd(\Dd)$, depending on whether the distributive laws in its 1h-cells 
are taken as lax or oplax.

\begin{prop}
The following two pairs of double categories are isomorphic:
$$\Lax_{hop}^*(*,\Dd)\iso\Mnd(\Dd)$$ 
and 
$$q\x\Lax_{hop}^*(*\times *,\Dd)\iso\Mnd(\Mnd(\Dd)).$$
\end{prop}

\begin{proof}
For a 1h-cell $\theta$ in $\Lax_{hop}^*(*,\Dd)$ we find the following. It is $(\theta)^{id_*}=\Id_{\theta_0(*)}$ (by (h.o.t.-4)\label{(h.o.t.-4)}). 
Additionally, since the only 2-cell in the double category $*$ is the trivial one, the axioms (h.o.t.-3) \label{(h.o.t.-3)} and (h.o.t.-5) \label{(h.o.t.-5)} 
hold trivially. On the other hand, $(\theta)_{id_*}$ is a non-trivial 2-cell such that (h.o.t.-1) \label{(h.o.t.-1)} and (h.o.t.-2) \label{(h.o.t.-2)} 
mean that a 1h-cell $\theta$ in $\Lax_{hop}^*(*,\Dd)$ (given thus only by $(\theta)_{id_*}$) is a monad-monad distributive law, {\em i.e.} 
a 1h-cell in $\Mnd(\Dd)$. 

For a 1v-cell $\theta_0$ in $\Lax_{hop}^*(*,\Dd)$ the situation is similar. One has $(\theta_0)^{id_*}=\Id^{\theta_0(*)}$ (by (v.l.t.-4)\label{(v.l.t.-4)}). 
The only non-trivial laws now are (v.l.t.-1) \label{(v.l.t.-1)} and (v.l.t.-2)\label{(v.l.t.-2)}. They involve a non-trivial 2-cell 
$(\theta_0)_{id_*}$ and they precisely mean that a 1v-cell $\theta_0$ in $\Lax_{hop}^*(*,\Dd)$ (given thus by $(\theta_0)_{id_*}$) is a 1v-cell in $\Mnd(\Dd)$. 

For a modification in $\Lax_{hop}^*(*,\Dd)$, which is given by a 2-cell 
$\scalebox{0.86}{
\bfig
\putmorphism(-180,50)(1,0)[F(A)` G(A)`\theta(*)]{550}1a
\putmorphism(-180,-270)(1,0)[F\s'(A)`G'(A) `\tilde\theta(*)]{550}1b
\putmorphism(-170,50)(0,-1)[\phantom{Y_2}``\theta_0(*)]{320}1l
\putmorphism(350,50)(0,-1)[\phantom{Y_2}``\theta_0'(*)]{320}1r
\put(0,-140){\fbox{$\Theta_*$}}
\efig}$, 
by triviality of $(\theta)^{id_*}$ and $(\theta_0)^{id_*}$ the axiom (m.ho-vl.-2) \label{m.ho-vl.-2} is trivial. The other axiom (m.ho-vl.-1) \label{m.ho-vl.-1} 
is the only possible identity relating $\Theta_*, (\theta)_{id_*}$ and $(\theta_0)_{id_*}$, and it 
precisely means that $\Theta_*$ is a 2-cell in $\Mnd(\Dd)$. 

\medskip

The inspection for $q\x\Lax_{hop}^*(*\times *,\Dd)$ goes similarly. A 1h-cell there is now a pair of 1h-cells in $\Lax_{hop}^*(*,\Dd)$ which 
relate according to $(HOT^q_1)$ $\label{(HOT^q_1)}$, as the other axioms are trivial now. These data precisely define a 1h-cell in $\Mnd(\Mnd(\Dd)).$ 
The situation for 1v-cells is symmetric, now the only non-trivial axiom is $(VLT^q_4)$ $\label{(VLT^q_4)}$. 

A modification in $q\x\Lax_{hop}^*(*\times *,\Dd)$ is given by a pair of modifications in $\Lax_{hop}^*(*,\Dd)$: 
$$\scalebox{0.86}{
\bfig
\putmorphism(-180,50)(1,0)[\A` \A'`\theta^1]{550}1a
\putmorphism(-180,-270)(1,0)[\tilde\A`\tilde{\A'} `\tilde\theta^1]{550}1b
\putmorphism(-170,50)(0,-1)[\phantom{Y_2}``\theta_0^1]{320}1l
\putmorphism(350,50)(0,-1)[\phantom{Y_2}``\theta_0^{1'}]{320}1r
\put(0,-140){\fbox{$\tau^1$}}
\efig}
\qquad\qquad
\scalebox{0.86}{
\bfig
\putmorphism(-180,50)(1,0)[\A` \A'`\theta^2]{550}1a
\putmorphism(-180,-270)(1,0)[\tilde\A`\tilde{\A'} `\tilde\theta^2]{550}1b
\putmorphism(-170,50)(0,-1)[\phantom{Y_2}``\theta_0^2]{320}1l
\putmorphism(350,50)(0,-1)[\phantom{Y_2}``\theta_0^{2'}]{320}1r
\put(0,-140){\fbox{$\tau^2$}}
\efig}
$$
which satisfy: $\tau^1_*=\tau^2_*$. This means that the 2-cell $\tau^1_*=\tau^2_*$ obeys two identities: one relating 
$\tau^1_*, (\theta^1)_{id_*}$ and $(\theta_0^1)_{id_*}$, and another relating $\tau^2_*, (\theta^2)_{id_*}$ and $(\theta_0^2)_{id_*}$. 
This gives a 2-cell in $\Mnd(\Mnd(\Dd)).$  
\qed\end{proof}

Recall that $\Mnd$ is an endofunctor on the category $2\x\Cat$ of 2-categories and 2-functors 
(it sends a 2-category $\K$ to the 2-category $\Mnd(\K)$), and that there is a natural transformation $\Comp:\Mnd\Mnd\to\Mnd$, 
which evaluated at $\K$ sends a distributive law in $\K$ to the induced composite monad, \cite{S}. 

Analogously, we can see $\Mnd$ as an endofunctor on the category $Dbl_{st}$ of double categories and strict double functors 
by the construction from \cite{FGK}. Moreover, we can consider a natural transformation $\Comp:\Mnd\Mnd\to\Mnd$, 
which evaluated at $\Dd$ sends a distributive law in $\Dd$ to the induced composite monad in $\Dd$. 

In the present setting, being $\Aa=\Bb=*$, note that we can write the double functor $\F$ from \equref{2-functor F} as 
$\F: q\x\Lax_{hop}(*\times *,\Dd)\to\Lax_{hop}(*,\Dd)$. Then one has that the following diagram commutes 
\begin{equation} \eqlabel{dist-mnd diag}
\scalebox{0.84}{
\bfig
\putmorphism(-430,330)(1,0)[q\x\Lax_{hop}(*\times *,\Dd)`\Lax_{hop}(*,\Dd)`\F]{1530}1a
\putmorphism(-480,320)(0,-1)[``\iso]{320}1l
\putmorphism(-400,0)(1,0)[\Mnd(\Mnd(\Dd))`\Mnd(\Dd).`\Comp(\Dd)]{1530}1b
\putmorphism(1060,330)(0,-1)[``\iso]{330}1r
\efig}
\end{equation}
Moreover, it indicates that the general double functor $\F$ can be seen as a sort of generalization of the double functor $\Comp(\Dd)$. 
Either of the two horizontal arrows in this diagram corresponds to the double categorification of Beck's result, that given a monad-monad 
distributive law between monads $T$ and $S$ on a category (given by a natural transformation $\phi: ST\Rightarrow TS$), 
$TS$ is a monad. In terms of our double functor $\F$, the 2-cell $\phi: ST\Rightarrow TS$ corresponds to the 2-cell 
$\gamma_{(id_*,id_*)(id_*,id_*)}$ from \ssref{F on 0}. 

\bigskip


To further develop applications to monads, there are some prospects, along the lines mentioned in the 2-categorical setting at the end of \cite {FMS}. 
We highlight the cases when $\Dd$ is the double category $\D\x\Mat$ of matrices in a category $\D$ with coproducts (see \cite{Power}), 
and $\Span(\C)$ of spans in a category $\C$ with pullbacks, as at the beginning of this Subsection. 
The vertical categories of the double categories $\D\x\Mat$ and $\Span(\C)$ seen as internal categories are then 
$\D\x\Cat$, the category of categories enriched over $\D$, and $\Cat(\C)$, the category of categories internal in $\C$. 
The double functor $\F$ in \equref{dist-mnd diag} then may yield some kind of product on the categories $\D\x\Cat$ and $\Cat(\C)$, 
under certain assumptions.

\bigskip

\bigskip

{\bf Acknowledgments.} 
The author was supported by the Science Fund of the Republic of Serbia, Grant No. 7749891, Graphical Languages - GWORDS. 
My deep and humble gratitude to the referee for her or his commitment and devotion to indicate me how the presentation of the results of this research 
and the results themselves could be taken to a higher level.

\end{document}